\documentclass{amsart}
\newtheorem{theorem}{Theorem}[section]
\newtheorem{lemma}[theorem]{Lemma}

\usepackage{graphicx}

\newtheorem{definition}[theorem]{Definition}
\newtheorem{corollary}[theorem]{Corollary}
\newtheorem{example}[theorem]{Example}
\newtheorem{remark}[theorem]{Remark}

\numberwithin{equation}{section}

\def\Xint#1{\mathchoice
{\XXint\displaystyle\textstyle{#1}}%
{\XXint\textstyle\scriptstyle{#1}}%
{\XXint\scriptstyle\scriptscriptstyle{#1}}%
{\XXint\scriptscriptstyle\scriptscriptstyle{#1}}%
\!\int}
\def\XXint#1#2#3{{\setbox0=\hbox{$#1{#2#3}{\int}$ }
\vcenter{\hbox{$#2#3$ }}\kern-.6\wd0}}

\def\dashint{\Xint-}

\begin{document}

\title[Wiener-type Criterion for the Removability of the Singularity]{Wiener-type Criterion for the Removability of the Fundamental Singularity for the Heat Equation and its Consequences }
\author{Ugur G. Abdulla}
\address{Analysis \& PDE Unit, Okinawa Institute of Science and Technology, Okinawa, Japan}
%\curraddr{Department of Mathematics and Statistics,
%Case Western Reserve University, Cleveland, Ohio 43403}
%\email{xyz@math.university.edu}


%    General info

%\date{January 1, 2001 and, in revised form, June 22, 2001.}



\begin{abstract}
We prove the necessary and sufficient condition for the removability of the fundamental singularity, and equivalently for the unique solvability of the singular Dirichlet problem for the heat equation. In the measure-theoretical context, the criterion determines whether the $h$-parabolic measure of the singularity point is null or positive. From the probabilistic point of view, the criterion presents an asymptotic law for conditional Brownian motion. In {\it U.G. Abdulla, J Math Phys, 65, 121503 (2024)} the Kolmogorov-Petrovsky-type test was established. Here, we prove a new Wiener-type criterion for the "geometric" characterization of the removability of the fundamental singularity for arbitrary open sets in terms of the fine-topological thinness of the complementary set near the singularity point. In the special case when the boundary of the open set is locally represented by a graph, the minimal thinness criterion for the removability of the singularity is expressed in terms of the minimal regularity of the boundary manifold at the singularity point.
\end{abstract}

\maketitle
\tableofcontents
%{\it Dedicated to Norbert Wiener, on the occassion of $100^{th}$ anniversary of his seminal papers in {\it J Math Phys, 3 (1924).}} 

\section{Introduction and Historical Remarks}
The major problem in the Analysis of PDEs is understanding the nature of singularities in solutions to the PDEs that reflect natural phenomena. 
It would be convenient to make some remarks on the analysis of singularities for the Laplace and heat equations, as well as more general second-order elliptic and parabolic PDEs.
The solvability, in some generalized sense, of the classical DP in a bounded open set $E \subset \mathbb{R}^N$, with prescribed data on $\partial E$, is realized within the class of resolutive boundary functions, identified by Perron's method and its Wiener \cite{Wiener1,Wiener2} and Brelot \cite{Brelot} refinements.
Such a method is referred to as the PWB method, and the corresponding solutions are PWB solutions. Paralleling the theory of PWB solutions, the DP for the heat equation in an arbitrary open set is solvable within the class of resolutive boundary functions. We refer to \cite{watson,Doob} for an account of the theory. Wiener, in his pioneering works \cite{Wiener1,Wiener2}, proved a necessary and sufficient condition for the finite boundary point $x_o \in \partial E$ to be regular in terms of the ``thinness" of the complementary set in the neighborhood of $x_o$. If the boundary of the domain is a graph in a neighborhood of $x_0$, the Wiener criterion is entirely geometrical. A key advance made in Wiener's work was the introduction of the concept of capacity - a sub-additive set function dictated by the Laplacian for the accurate measuring of the thinness of the complementary set in the neighborhood of $x_0$ for the boundary regularity of a harmonic function. Formalized through the powerful Choquet capacitability theorem \cite{Choquet}, the concept of capacity became a standard tool for the characterization of singularities for the elliptic and parabolic equations. The question of removability of isolated singularities for the linear second-order elliptic and parabolic PDEs was settled in \cite{SerrinW}, and in \cite{Aronson,Aronson1,Aronson2}. 
The Wiener criterion for the boundary continuity of harmonic functions became a canonical result, driving the boundary regularity theory for the elliptic and parabolic PDEs. 
In 1935, Petrovsky proved a geometric necessary and sufficient condition for the regularity of the characteristic top boundary point for the heat equation in the domain of revolution \cite{Petrovski} (see also \cite{AbdullaKP}). In the same paper, he also presented an elegant solution of the Kolmogorov problem (see Section~\ref{prelude}, Problem~$\tilde{\mathcal{A}}_0$) for the special domain of revolution (see also \cite{Abdulla6}). The results formed the so-called Kolmogorov-Petrovsky test for the asymptotic behavior of the standard Brownian path as $t\downarrow 0$ and $t\uparrow+\infty$, and opened a path for the deep connection between the regularity theory of elliptic and parabolic PDEs and asymptotic properties of the associated Markov processes \cite{ito2}. The geometric iterated logarithm test for the regularity of the boundary point for an arbitrary open set with respect to the heat equation is proved in \cite{AbdullaSiam}. 
Paralleling the Wiener regularity theory, Wiener's criterion for the regularity of the finite boundary point for the heat equation was formulated in \cite{lanconelli} along with the proof of the irregularity assertion. The problem was accomplished in \cite{evansgariepy}, where the long-awaited regularity assertion was proved. As in its elliptic counterpart, the concept of {\it heat capacity} was a key concept to extend the Wiener regularity theory to the case of the heat equation \cite{watson}. However, the major technical difficulty in doing so was related to the nature of the singularities of the fundamental solutions of the Laplace and heat equations. The former is an isolated singularity for the spherical level sets of the fundamental solution, while the latter is a non-isolated singularity point for the level sets of the fundamental solution of the heat equation. To complete the Wiener regularity result at finite boundary points for the heat equation, the major technical advance of the paper \cite{evansgariepy} was a proof of an elegant boundary Harnack estimate near the non-isolated singularity point of the level sets of the fundamental solution to the backward heat equation. The result of \cite{evansgariepy} was extended to the class of linear second-order divergence-form parabolic PDEs with $C^1$-Dini-continuous coefficients in \cite{garofalo1,garofalo2}.  

In \cite{LSW} it is proved that the Wiener test for the regularity
of finite boundary points concerning a second-order divergence form uniformly elliptic operator with bounded measurable coefficients
coincides with the classical Wiener test for the boundary regularity of harmonic functions. The Wiener test for the regularity of finite boundary points for linear degenerate elliptic equations is proved in \cite{FJK}. 
The Wiener test for the regularity of finite boundary points for quasilinear 
elliptic equations was settled due to \cite{Mazya3,gariepy,kilpi}. Nonlinear potential theory of degenerate elliptic PDEs was developed along the same lines as classical potential theory
for the Laplace operator, for which we refer to \cite{HKM}. As in the classical case, the well-posedness and boundary regularity in arbitrary irregular domains for nonlinear degenerate and singular parabolic PDEs are more challenging than in their elliptic counterparts. In \cite{A4,A7,A14} the existence, boundary regularity, uniqueness, and $L_1$-contraction estimates are proved for the nonlinear degenerate and singular parabolic PDEs in general non-cylindrical domains under minimal regularity conditions on the boundary manifolds. Precisely, assuming that the boundary is locally a continuous graph, the well-posedness and boundary regularity results are established under the minimal regularity assumption on the local modulus of left-lower semicontinuity of the boundary manifold (see Section~\ref{geometrictest}). 

To solve the DP in an unbounded open set, Brelot introduced the idea of compactifying $\mathbb{R}^N$ into $\mathbb{R}^N\cup\{\infty\}$, where $\infty$ is the point at $\infty$ of $\mathbb{R}^N$ \cite{brelot1}. The PWB method is extended to the compactified framework, thus providing a powerful existence and uniqueness result for the DP in arbitrary open sets in the class of resolutive boundary functions. The new concept of regularity of $\infty$ was introduced in \cite{Abdulla1} for the classical DP, and in \cite{Abdulla4} for its parabolic counterpart. The DP with a bounded Borel measurable boundary function has one and only one or infinitely many solutions without prescribing the boundary value at $\infty$. The point at $\infty$ is called regular if there is a unique solution, and it is called irregular otherwise. Equivalently, in the measure-theoretical context, the new concept of regularity or irregularity of $\infty$ is introduced according to whether the harmonic measure of $\infty$ is null or positive. In \cite{Abdulla1}, the Wiener criterion for the regularity of $\infty$ in the classical DP for the Laplace equation in an open set $E\subset\mathbb{R}^N$ with $N\geq 3$ is proved. In \cite{Abdulla4}, it is proved that the Wiener criterion at $\infty$ for the linear second-order divergence form elliptic PDEs with bounded measurable coefficients coincides with the Wiener criterion at $\infty$ for the Laplacian operator. The Wiener criterion at $\infty$ for the heat equation is proved in \cite{Abdulla2}. Remarkably, the Kolmogorov problem (see Section~\ref{prelude}, Problem~$\tilde{\mathcal{A}}_0$) is a particular case of the problem of uniqueness of the bounded solution of the parabolic Dirichlet problem in an arbitrary open set in $\mathbb{R}^{N+1}$ without prescribing the limit of the solution at $\infty$. Hence, the Wiener criterion at $\infty$ proved in \cite{Abdulla2} presents a full solution to the Kolmogorov problem. 

The new Wiener criterion at $\infty$ for the elliptic and parabolic PDEs broadly extends the role of the Wiener regularity theory in mathematics. The Wiener test at $\infty$ arises as a global characterization of uniqueness in boundary value problems in arbitrary unbounded open sets. From a topological point of view, the Wiener test at $\infty$ arises as a thinness criterion at $\infty$ in fine topology. In a probabilistic context, the Wiener test at $\infty$ characterizes asymptotic laws for the Markov processes whose generator is a given differential operator. The counterpart of the new Wiener test at a finite boundary point leads to uniqueness in a Dirichlet problem for a class of unbounded functions growing at a certain rate near the boundary point; a criterion for the removability of singularities and/or for unique continuation at the finite boundary point: let $E\subset \mathbb{R}^{N}, N\geq 3$ be an open set, and $x_0\in E$ be a finite boundary point. Consider a singular Dirichlet problem for the linear second-order uniformly elliptic PDE with bounded measurable coefficients in a class $O(|x-x_0|^{2-N}$ as $x\to x_0$. In \cite{Abdulla4} it is proved that the Wiener test at $x_0$ is a necessary and sufficient condition for the unique solvability of the singular Dirichlet problem, and equivalently for the removability of the fundamental singularity at $x_0$. In a recent paper \cite{Abdulla5}, an appropriate 2D analog of this result is established. 
Let $E\subset\mathbb{R}^2$ be a Greenian open set, and $x_0\in\partial E$ be a boundary point (finite or $\infty$). Consider a singular Dirichlet problem for the linear second-order
uniformly elliptic operator with bounded measurable coefficients in the class $O(\log|x-x_0|)$ if $x_0$ is finite, and in a class of functions with logarithmic growth, if $x_0=\infty$. 
In \cite{Abdulla5} it is proved that the Wiener criterion at $x_0$ is a necessary and sufficient condition for the unique solvability of the singular Dirichlet problem, and equivalently for the removability of the logarithmic singularity. Precisely, in \cite{Abdulla5}, the concept of $\log$-regularity (or $\log$-irregularity) of the boundary point (finite or $\infty$) is introduced according as the $\log$-harmonic measure of it is null or positive. The removability of the logarithmic singularity is expressed in terms of the Wiener criterion for the $\log$-regularity of $x_0$. 

The goal of this paper is to establish a necessary and sufficient condition for the removability of the fundamental singularity, and equivalently for the unique solvability of the singular PDP in arbitrary open sets. In a recent paper \cite{AbdullaI}, we proved the Kolmogorov-Petrovsky-type test for a special case of domains with a boundary formed by the surface of revolution around the time axis. In this paper, we introduce the concept of $h$-capacity and prove a Wiener-type criterion for the removability of fundamental singularity for arbitrary open sets expressed in terms of the fine-topological thinness of the complementary set near the singularity point. In the special case when the boundary of the open set is locally represented by a graph, we express the minimal thinness criteria for the removability of the singularity in terms of the minimal regularity of the boundary manifold at the singularity point.
\section{Prelude: Kolmogorov-type Problem}\label{prelude}
Consider the fundamental solution of the heat equation:
\begin{equation}\label{fs}
F(x,t)=
\left\{
\begin{array}{l}
(4\pi t)^{-\frac{N}{2}} e^{-\frac{|x|^2}{4t}}, \ x\in\mathbb{R}^N, t>0,\\
0, \ x\in\mathbb{R}^N, x\neq 0, t=0.
\end{array}\right.
\end{equation}
It is a distributional solution of the Cauchy Problem
\[ \mathcal{H}F:=F_t-\Delta F=0 \ \text{in} \ \mathbb{R}^{N+1}_+, \ F=\delta \ \text{on} \ \mathbb{R}^N \times \{t=0\} \]
where $\delta$ is a unit measure with support at $x=0$. For any fixed $\gamma\in \mathbb{R}^{N}$, let
\[ h(x,t)=F(x-\gamma,t) \]
be a fundamental solution with a pole at $\mathcal{O}:=(\gamma,0)$. 
The singularity of $h$ at $\mathcal{O}$ represents the natural phenomenon of the space-time distribution of the unit energy initially blown up at a single point. The fundamental singularity is non-removable 
for the heat equation in  $\mathbb{R}^{N+1}_+$. In particular, the Cauchy Problem for the heat equation in  $\mathbb{R}^{N+1}_+$ has infinitely many solutions in class
$O(h)$. 

The goal of this paper is to reveal the criterion for the removability of the fundamental singularity for open subsets of $\mathbb{R}^{N+1}_+$.
Let $\Omega\subset\mathbb{R}^{N+1}_+$ be an arbitrary open set and $\partial\Omega\cap\{t=0\}=\{\mathcal{O}\}$. Let $g:\partial\Omega\to \mathbb{R}$ be a boundary function such that $g/h$ is a bounded Borel measurable.
 Consider a {\it singular} {\bf parabolic Dirichlet problem}(PDP):
\begin{equation}\label{sDproblem}
\mathcal{H}u=0 \quad\text{in} \ \Omega, \  u=g \quad\text{on} \ \partial\Omega\setminus\{\mathcal{O}\};  \   u=O(h) \quad\text{at} \ \mathcal{O}.
\end{equation}
Solution of the singular PDP is understood in Perron's sense [see Section~\ref{E:1:1}, and formulae \eqref{Eq:W:1:7add} and \eqref{Eq:W:1:8add}]. Furthermore, the expression
"prescribing the behavior of $u/h$ at $\mathcal{O}$" is understood in the sense of requiring relations \eqref{Eq:W:1:7add} and \eqref{Eq:W:1:8add}) at the boundary point $\mathcal{O}$ .
Without prescribing the behavior of $u/h$ at $\mathcal{O}$, there exists either one and only one or infinitely many solutions of PDP (see Section~\ref{E:1:2} and formula \eqref{Eq:W:1:8h-1}). 
The main goal of this paper is to find a necessary and sufficient condition for open sets $\Omega$ for the uniqueness of the solution to the PDP without prescribing $u/h$ at $\mathcal{O}$. Note that we are not assuming that the open set $\Omega$ is connected. Furthermore, without loss of generality, we assume that
\begin{equation}\label{assonomega}
\mathcal{O}\in \partial \Lambda(z,\Omega), \ \quad\text{for some} \ \Lambda(z,\Omega)\subset \Omega,
\end{equation}
where $\Lambda(z,\Omega)$ for $z\in \Omega$ is a set of points $w\in \Omega\setminus \{z\}$ which can be joined to $z$ by a continuous curve in $\Omega$ along which $t$ is strictly increasing as the line is described from $w$ to $z$. Without this assumption within each connected component $\Lambda(z,\Omega)$, the solution is uniquely defined from boundary values on the parabolic boundary of this component. This makes the problem uniquely solvable in the whole $\Omega$ without prescribing $u/h$ at $\mathcal{O}$. 

The problem of removability vs. non-removability of the fundamental singularity is equivalent to the question of the uniqueness of the solution to PDP \eqref{sDproblem}
without prescribing the behavior of $u/h$ at $\mathcal{O}$. Note that we are not requiring the existence of the limit of $u/h$ at $\mathcal{O}$. Assume that $\lim_{z\to\mathcal{O}, z\in\partial\Omega\setminus\{\mathcal{O}\}}g/h$ exists. We prove in this case that the removability of the fundamental singularity at $\mathcal{O}$ is equivalent to the existence of a unique solution of the singular PDP \eqref{sDproblem} such that
\[ \lim_{z\to\mathcal{O}, z\in\Omega} \frac{u}{h} = \lim_{z\to\mathcal{O}, z\in\partial\Omega\setminus\{\mathcal{O}\}} \frac{g}{h}. \]
Otherwise speaking, for a unique solution $u$, $u/h$ will pick up the limit value at the singularity point $\mathcal{O}$ without being required (see Section~\ref{E:1:3} and Definition~\ref{defbr}).

The following procedure provides a key problem in testing the removability of fundamental singularity.
Let 
\[ \Omega_n:=\Omega\cap \{t>n^{-1}\}, \ n=1,2,... \]
and $u_n$ be a unique solution of the parabolic Dirichlet problem 
\begin{equation}\label{u_n}
\mathcal{H} u = 0, \quad\text{in} \ \Omega_n; \ u |_{\partial\Omega_n\cap \{t>n^{-1}\}}=0; \ u  |_{\partial\Omega_n\cap \{t=n^{-1}\}}=h(x,n^{-1}).
\end{equation}
From the maximum principle, it follows that
\begin{equation}\label{mp}
 0\leq u_{n+1}(x,t)\leq u_n(x,t) \leq h(x,t), \quad\text{on} \ \Omega_n. 
 \end{equation}
Therefore, there exists a limit function
\begin{equation}\label{ulim}
u_*(x,t)=\lim_{n\to +\infty}u_n(x,t), \ (x,t)\in \Omega,
\end{equation}
which satisfies \eqref{sDproblem}, and
\[ 0\leq u_*(x,t) \leq h(x,t), \ (x,t) \in \Omega.\]
The following is the key problem.

{\bf Problem} $\mathcal{A}_\gamma$: {\it Is $u_* \equiv 0$ in $\Omega$ or $u_*(x,t)\not\equiv 0$ in $\Omega$?  Equivalently, is the fundamental singularity at $\mathcal{O}$ removable or nonremovable for $\Omega$?}

Note that the Problems $\mathcal{A}_\gamma$ are equivalent for various $\gamma$. Indeed, the heat equation is translation invariant, so that solving the problem when the pole of the fundamental solutions is $(0,0)$ is equivalent to
solve the problem when the pole is at $(\gamma, 0)$ for any $\gamma\in \mathbb{R}^N$.

Next, we formulate the equivalent problem in $\mathbb{R}^{N+1}_-$. In that context, we are going to consider one-point Alexandrov compactification: $\mathbb{R}^{N+1}_-\rightarrow \mathbb{R}^{N+1}_-\cup \{\infty\}$. 
For any fixed finite $\gamma\in\mathbb{R}^N$, consider a function
\begin{equation}\label{afs}
\tilde{h}(x,t)= e^{\langle x,\gamma\rangle + |\gamma|^2 t}.
\end{equation}
It is a positive solution of the heat equation in $\mathbb{R}^{N+1}_-$. If $\gamma=0$, then $\tilde{h}\equiv 1$, while in the case when $\gamma\neq 0$, it is
an unbounded solution with singularity at $\infty$. The singularity is not removable for the heat equation in $\mathbb{R}^{N+1}_-$. We aim to reveal the criterion for the removability of the fundamental singularity for the open subsets of $\mathbb{R}^{N+1}_-$.
Let $\tilde{\Omega}\subset\mathbb{R}^{N+1}_-$ be an arbitrary open set, and $g:\partial\tilde{\Omega}\to \mathbb{R}$ be a boundary function, such that $g/\tilde{h}$ is a bounded Borel measurable. Consider a {\it singular} {\bf parabolic Dirichlet problem}(PDP):
\begin{equation}\label{sDproblem-}
\mathcal{H}u=0 \quad\text{in} \ \tilde{\Omega}, \  u=g \quad\text{on} \ \partial\tilde{\Omega}; \   u=O(\tilde{h}) \quad\text{at} \ \infty. 
\end{equation}
Solution of the singular PDP is understood in Perron's sense [see Section~\ref{E:1:1}, and formulae \eqref{Eq:W:1:7addtilde} and \eqref{Eq:W:1:8addtilde}]. Furthermore, the expression
"prescribing the behavior of $u/\tilde{h}$ at $\infty$" is understood in the sense of requiring relations \eqref{Eq:W:1:7addtilde} and \eqref{Eq:W:1:8addtilde}) at the boundary point $\infty$. 
Without prescribing the behavior of $u/\tilde{h}$ at $\infty$ there exists either one and only one or infinitely many solutions of PDP (see Section~\ref{E:1:2} and formula \eqref{Eq:W:1:8h-1tilde}). The main goal of this paper is to find a necessary and sufficient condition for open sets $\tilde{\Omega}$ for the uniqueness of the solution to the PDP without prescribing $u/\tilde{h}$ at $\infty$. We are not assuming that the open set $\tilde{\Omega}$ is connected. Furthermore, without loss of generality, we assume that
\begin{equation}\label{assontildeomega}
\infty\in \partial \Lambda(z,\tilde{\Omega}), \ \quad\text{for some} \ \Lambda(z,\tilde{\Omega})\subset \tilde{\Omega},
\end{equation}
Without this assumption within each connected component $\Lambda(z,\tilde{\Omega})$, the solution is uniquely defined from boundary values on the parabolic boundary of this component. This makes the problem uniquely solvable in the whole $\tilde{\Omega}$ without prescribing $u/\tilde{h}$ at $\infty$. 

The problem of removability vs. non-removability of the fundamental singularity is equivalent to the question of the uniqueness of the solution to PDP \eqref{sDproblem-}
without prescribing the behavior of $u/\tilde{h}$ at $\infty$. Note that we are not requiring the existence of the limit of $u/\tilde{h}$ at $\infty$. Assume that $\lim_{z\to\infty, z\in\partial\tilde{\Omega}\}}g/\tilde{h}$ exists. We prove in this case that the removability of the fundamental singularity at $\infty$ is equivalent to the existence of a unique solution of the singular PDP \eqref{sDproblem-} such that
\[ \lim_{z\to\infty, z\in\tilde{\Omega}} \frac{u}{\tilde{h}} = \lim_{z\to\infty, z\in\partial\tilde{\Omega}\}} \frac{g}{\tilde{h}}. \]
Otherwise speaking, for a unique solution $u$, $u/\tilde{h}$ will pick up the limit value at the singularity point $\infty$ without being required (see Section~\ref{E:1:3} and Definition~\ref{defbrtilde}).
Similar to its counterpart \eqref{sDproblem}, the key problem to test the removability of the singularity at $\infty$ is formulated as follows:

Let  
\[ \tilde{\Omega}_n:=\tilde{\Omega}\cap \{t>-n\}, \ n=1,2,... \]
and $\tilde{u}_n$ be a unique solution of the parabolic Dirichlet problem 
\begin{equation}\label{u_n-}
\mathcal{H} u = 0, \quad\text{in} \ \tilde{\Omega}_n; \ u |_{\partial\tilde{\Omega}_n\cap \{t>-n\}}=0; \ u  |_{\partial\tilde{\Omega}_n\cap \{t=-n^\}}=\tilde{h}(x,-n).
\end{equation}
From the maximum principle, it follows that
\begin{equation}\label{mp-}
 0\leq \tilde{u}_{n+1}(x,t)\leq \tilde{u}_n(x,t) \leq \tilde{h}(x,t), \quad\text{on} \ \tilde{\Omega}_n. 
 \end{equation}
Therefore, there exists a limit function
\begin{equation}\label{u*=lim-}
\tilde{u}_*(x,t)=\lim_{n\to +\infty}\tilde{u}_n(x,t), \ (x,t)\in \tilde{\Omega},
\end{equation}
which satisfies \eqref{sDproblem-}, and
\[ 0\leq \tilde{u}_*(x,t) \leq \tilde{h}(x,t), \ (x,t) \in \tilde{\Omega}.\]
The following is the key problem.

{\bf Problem} $\tilde{\mathcal{A}}_\gamma$: Is $\tilde{u}_* \equiv 0$ in $\tilde{\Omega}$ or $\tilde{u}_*(x,t)\not\equiv 0$ in $\tilde{\Omega}$?  Equivalently, is a fundamental singularity at $\infty$ removable or nonremovable for $\tilde{\Omega}$?

\begin{remark}
Problem~$\tilde{\mathcal{A}}_\gamma$ can be formulated in $\mathbb{R}^{N+1}$ without any change. Indeed, the parabolic Dirichlet problem for the heat equation is uniquely solvable in any open subset  $\mathbb{R}^{N+1}_+$ in a class $O(\tilde{h})$. Therefore, given an arbitrary open set $\Omega\subset \mathbb{R}^{N+1}$, the solution of the parabolic Dirichlet problem in $\Omega$ can be constructed as a unique continuation of the solution to the parabolic Dirichlet problem in $\Omega_-=\Omega\cap \mathbb{R}^{N+1}_-$. Moreover, the latter is independent of the boundary values assigned on $\partial\Omega_-\cap\{t=0\}$, since $\partial\Omega_-\cap\{t=0\}$ is a parabolic measure null subset of $\partial\Omega_-$. This implies that the Problem~$\tilde{\mathcal{A}}_\gamma$ is equivalent for $\Omega$ and $\Omega_-$. Otherwise speaking, the fundamental singularity at $\infty$ is removable for $\Omega \subset\mathbb{R}^{N+1}$ if and only if it is removable for $\Omega_-$. 
\end{remark}
%\begin{remark}
%Due to the translation invariance of the heat equation, the solution of the Problem~${\cal A}$ is independent of $\gamma$. However, the solution of the Problem~$\tilde{\mathcal A}}_\gamma$ depends on $\gamma$. 
%\end{remark}

The only problem in the family of formulated problems that is solved is the Problem~$\tilde{\mathcal{A}}_\gamma$ when $\gamma=0$ (or Problem~$\tilde{\mathcal{A}}_0$). The Problem~$\tilde{\mathcal {A}}_0$ was formulated by Kolmogorov in 1928 in the seminar on probability theory at Moscow State University, in the particular case with
$\Omega=\{|x|<f(t), -\infty<t<0\}\subset \mathbb{R}^2$, with $f\in C(-\infty,0]$ such that $f(-\infty)=+\infty, f \uparrow +\infty, (-t)^{-\frac{1}{2}} f \uparrow +\infty$ as $t \downarrow -\infty$. Kolmogorov's motivation for posing this problem was a connection to the probabilistic problem of finding the asymptotic behavior at infinity of the standard Brownian path. Let $\{\xi(t): t\geq 0, P_\bullet\}$ be a standard 1-dimensional Brownian motion and $P_\bullet({\bf B})$ is the probability of the event ${\bf B}$ as a function of the starting point $\xi(0)$. Consider  the event
\[ {\bf B}=\Big [ \quad\text{the set} \ \{t\geq 0: \xi(t)>f(-t)\} \quad\text{clusters to} \ +\infty \Big ]  \]
Blumenthal's 01 law implies that $P_0({\bf B})=0$ or $1$; $f(-t)$ is said to belong to the lower class if this probability is 1 and to the upper class otherwise. 
Remarkably, the Kolmogorov Problem's solution $u_*$ is $=0$ or $>0$ according to whether $f(-t)$ is in the lower or upper class, respectively. The Kolmogorov Problem in a one-dimensional setting was solved by Petrovsky in 1935, and the celebrated result is called the Kolmogorov-Petrovski test in the probabilistic literature \cite{Petrovski} (see also \cite{Abdulla6}). 

The full solution of the Kolmogorov Problem for arbitrary open sets $\Omega$ (or Problem~$\tilde{\mathcal{A}}_0$) is presented in \cite{Abdulla2}. A new concept of regularity or irregularity of $\infty$ is introduced according to whether the parabolic measure of $\infty$ is null or positive, and the necessary and sufficient condition for the Problem~$\tilde{\mathcal{A}}_0$ is proved in terms of the Wiener-type criterion for the regularity of $\infty$. 

In the probabilistic context, the formulated problems $\mathcal{A}_\gamma$ and $\tilde{\mathcal {A}}_\gamma$ are generalizations of the Kolmogorov problem to establish asymptotic laws for the $h$-{\it Brownian processes} \cite{Doob} (see Section~\ref{hBrownian}).

%%%%%%
%In this paper, we formulate and prove a new Wiener-type criterion for the "geometric" characterization of the removability of the fundamental singularity for arbitrary open sets in terms of the fine-topological thinness of the complementary set near the singularity point. We shall follow the notations of \cite{AbdullaI} and freely use the formulations of the problems, results, and techniques employed there. 

\section{Overture: Wiener-type Criterion} %for the Removability of the Fundamental Singularity
Next, we introduce a necessary terminology and formulate the Wiener-type criteria for the {\bf Problem}s $\mathcal{A}_\gamma$ and $\tilde{\mathcal{A}}_\gamma$.

We write typical points $z$ and $w\in \mathbb{R}^{N+1}$ as $z=(x,t), w=(y,\tau), x,y\in\mathbb{R}^N, t,\tau\in \mathbb{R}$. For any compact set $K\subset\mathbb{R}^{N+1}_+$, denote by $\mathcal{M}_K$ the set of all nonnegative Radon measures on $\mathbb{R}^{N+1}_+$, with support in $K$.  For $\mu\in\mathcal{M}_K$ introduce an $h$-{\it potential of} $\mu$:
\begin{equation}\label{potential}
\mathbb{P}^h_\mu(z)=\int\limits_{\mathbb{R}^{N+1}_+}\frac{F(z-w)}{h(z)h_*(w)}d\mu(w), \ z\in\mathbb{R}^{N+1}_+.
\end{equation}
Here 
\[ h_*(z):=\Big (\frac{\pi}{t}\Big )^{\frac{N}{2}} e^{\frac{|x-\gamma|^2}{4t}}, \ z=(x,t)\in\mathbb{R}^{N+1}_+\]
be a fundamental solution of the adjoined heat equation
\begin{equation}\label{*heat}
\mathcal{H}_*u:=-u_t-\Delta u=0
\end{equation}
with singularity at $\mathcal{O}$. 
Recall that for $z=(x,t), w=(y,\tau)$ we have
\begin{equation}\label{fkernel}
F(z-w)=
\left\{
\begin{array}{l}
(4\pi (t-\tau))^{-\frac{N}{2}} e^{-\frac{|x-y|^2}{4(t-\tau)}} \  t>\tau,\\
0 \ \ \ \ \ \ \ \ \ \ \ \ \ \ \ \ \ \ \ \ \ \ \ \ \ \ \ \ t\leq \tau.
\end{array}\right.
\end{equation}
Introduce the key concept of $h$-capacity.\\
\begin{definition}\label{hcapacity}
For a compact set  $K\subset \mathbb{R}^{N+1}_+$, the $h$-capacity of $K$ is
\begin{equation}\label{ChK}
 C_h(K)\equiv \max \{\mu(K): \ \mu\in \mathcal{M}_K, \ \mathbb{P}^h_\mu\leq 1 \ \text{in} \ \mathbb{R}^{N+1}_+\}.
 \end{equation}
 \end{definition}
There exists a unique measure $\lambda\in\mathcal{M}_K$, called an $h$-{\it capacitary measure} of $K$ such that (Lemma~\ref{hpotentialvsreduction1})
\[ \lambda(K)=C_h(K)\]
\begin{definition}\label{hheatball}
The $h$-{\it heat ball} of center $\overline{z}=(\gamma,1)$ and radius $4c(1+4c)^{-1}$ is
\begin{gather} B(\bar{z},c):=\Big \{w\in\mathbb{R}^{N+1}_+: \frac{F(\bar{z}-w)}{h(\bar{z})h_*(w)}>(4\pi c)^{-\frac{N}{2}}\Big \} \nonumber\\
 \equiv \Big \{z\in\mathbb{R}^{N+1}_+: \ |x-\gamma|^2<2Nt(1-t)\log \frac{4ct}{1-t}, \ \frac{1}{1+4c}<t<1 \Big \}\nonumber
 \end{gather}
 \end{definition}
Consider a closed $h$-{\it heat shell}
\begin{gather}
D(\bar{z},c)=\Big \{w\in\mathbb{R}^{N+1}_+: (2\pi c)^{-\frac{N}{2}}\geq \frac{F(\bar{z}-w)}{h(\bar{z})h_*(w)}\geq(4\pi c)^{-\frac{N}{2}}\Big \} \cup\{\bar{z}\} \nonumber\\
\equiv \overline{B(\overline{z},c)\setminus B(\overline{z},c/2)}; \nonumber
\end{gather}
This is the set of points between the level surfaces of the $h_*$-{\it parabolic} function $\frac{F(\bar{z}-\cdot)}{h_*(\cdot)}$ (see Section~\ref{averaging}). 
Given open set $\Omega\subset \mathbb{R}^{N+1}_+$ define (see Figure 1)
\[ E_n\equiv \Omega^c\cap D(\overline{z},2^n).\]
The solution of the {\bf Problem} $\mathcal{A}_\gamma$ for arbitrary open set $\Omega\subset\mathbb{R}^{N+1}_+$ reads:
 \begin{theorem}\label{wienertest}
$u_*\equiv 0$ or $u_*>0$, that is to say, the fundamental singularity is removable or non-removable
according as the series
\begin{equation}\label{wtest}
\sum\limits_n2^{-\frac{nN}{2}}C_h(E_n)
\end{equation} 
diverges or converges. An equivalent characterization is valid if the series \eqref{wtest} is replaced with
\[ \sum\limits_n\lambda^{-n}C_h\Big ( \Big \{ \lambda^{-n+1}\geq \frac{F(\bar{z}-w)}{h(\bar{z})h_*(w)} \geq\lambda^{-n}\Big \}\Big )  \ \text{for any} \ \lambda>1\]
\end{theorem}  
\begin{figure}[ht]\label{fig1}
    \centering
    \includegraphics[width=0.9\textwidth, trim=1cm 8cm 1cm 7cm, clip]{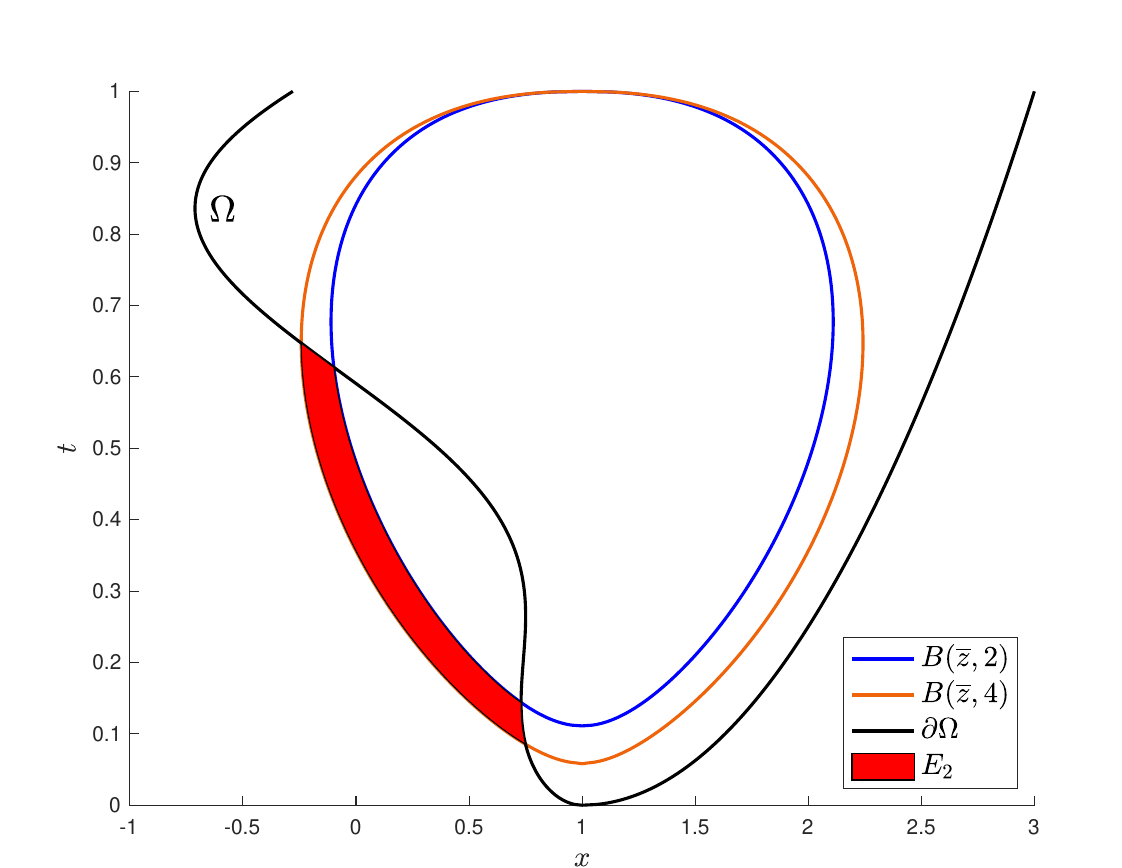}
    \caption{Sketch of $h$-heat balls $B(\bar{z},2), B(\bar{z},4)$, $h$-heat shell $D(\bar{z},4)$ and $E_2$ for a given open set $\Omega$ for $N=1$ and $\gamma = 1$}
\end{figure}

Next, we introduce a similar terminology for the formulation of the solution to the Problem~$\tilde{\mathcal{A}}_\gamma$. For any compact set $K\subset\mathbb{R}^{N+1}_-$, we adopt the same notation $\mathcal{M}_K$ for a set of nonnegative Radon measures on $\mathbb{R}^{N+1}_-$ with support in $K$. We introduce a $\tilde{h}$-{\it potential} of $\mu$ as follows:
\begin{equation}\label{tildehpotential}
\mathbb{P}^{\tilde{h}}_\mu(z)=\int\limits_{\mathbb{R}^{N+1}_-}\frac{F(z-w)}{\tilde{h}(z)\tilde{h}_*(w)}d\mu(w), \ z\in\mathbb{R}^{N+1}_-,
\end{equation}
where
\[ \tilde{h}_*(y,\tau):= e^{-\langle y,\gamma\rangle - |\gamma|^2 \tau}, \ (y,\tau)\in \mathbb{R}^{N+1}_- \]
be a solution of the adjoined heat equation \eqref{*heat} with singulaity at $\infty$. \\
\begin{definition}\label{tildehheat ball}
For a compact set $K\subset \mathbb{R}^{N+1}_-$, the $\tilde{h}$-capacity of $K$ is
\begin{equation}\label{ChK-}
C_{\tilde{h}}(K)\equiv \max \{\mu(K): \ \mu\in \mathcal{M}_K, \ \mathbb{P}^{\tilde{h}}_\mu\leq 1 \ \text{in} \ \mathbb{R}^{N+1}_-\}
\end{equation}
\end{definition}
There exists a unique measure $\tilde{\lambda}\in\mathcal{M}_K$, called an $\tilde{h}$-{\it capacitary measure} of $K$ such that (Lemma~\ref{hpotentialvsreduction1})
\[ \tilde{\lambda}(K)=C_{\tilde{h}}(K)\]
\begin{definition}\label{tildehheatball}
The $\tilde{h}$-{\it heat ball} of center $\overline{w}=(\gamma/2,-1/4)$ and radius $c$ is
\begin{gather}
\tilde{B}(\overline{w},c):=\Big \{w\in\mathbb{R}^{N+1}_-: \frac{F(\bar{w}-w)}{\tilde{h}(\bar{w})\tilde{h}_*(w)}>(4\pi c)^{-\frac{N}{2}}\Big \} \nonumber\\
\equiv \Big \{w=(y,\tau)\in\mathbb{R}^{N+1}_-: \ |y+2\tau\gamma|^2<-\frac{N}{2}(1+4\tau)\log \Big (\frac{-4c}{1+4\tau}\Big ), \ -\frac{1}{4}-c<\tau<-\frac{1}{4}\Big \}. \nonumber
\end{gather}
\end{definition}
Consider a closed $\tilde{h}$-{\it heat shell}
\begin{gather}
 \tilde{D}(\overline{w},c)=\Big \{w\in\mathbb{R}^{N+1}_-: (2\pi c)^{-\frac{N}{2}}\geq\frac{F(\bar{w}-w)}{\tilde{h}(\bar{w})\tilde{h}_*(w)}\geq (4\pi c)^{-\frac{N}{2}}\Big \}\cup \{\bar{w}\}\nonumber\\ \equiv \overline{\tilde{B}(\overline{w},c)\setminus \tilde{B}(\overline{w},c/2)}; \nonumber
 \end{gather}
this is the set of points between the level surfaces of the $\tilde{h}_*$-{\it parabolic} function $\frac{F(\bar{z}-\cdot)}{\tilde{h}_*(\cdot)}$ (see Section~\ref{averaging}).
Given open set $\tilde{\Omega}\subset \mathbb{R}^{N+1}_-$ define (see Figure 2)
\[ \tilde{E}_n\equiv \tilde{\Omega}^c\cap \tilde{D}(\overline{w},2^n).\]
%Let $\tilde{\mu}_n$ be a capacitary measure of $\tilde{E}_n$, i.e. $\tilde{\mu}(\tilde{E}_n)=C_{\tilde{h}}(\tilde{E}_n)$.
\begin{figure}[ht]\label{fig2}
    \centering
    \includegraphics[width=0.9\textwidth, trim=1cm 8cm 1cm 7cm, clip]{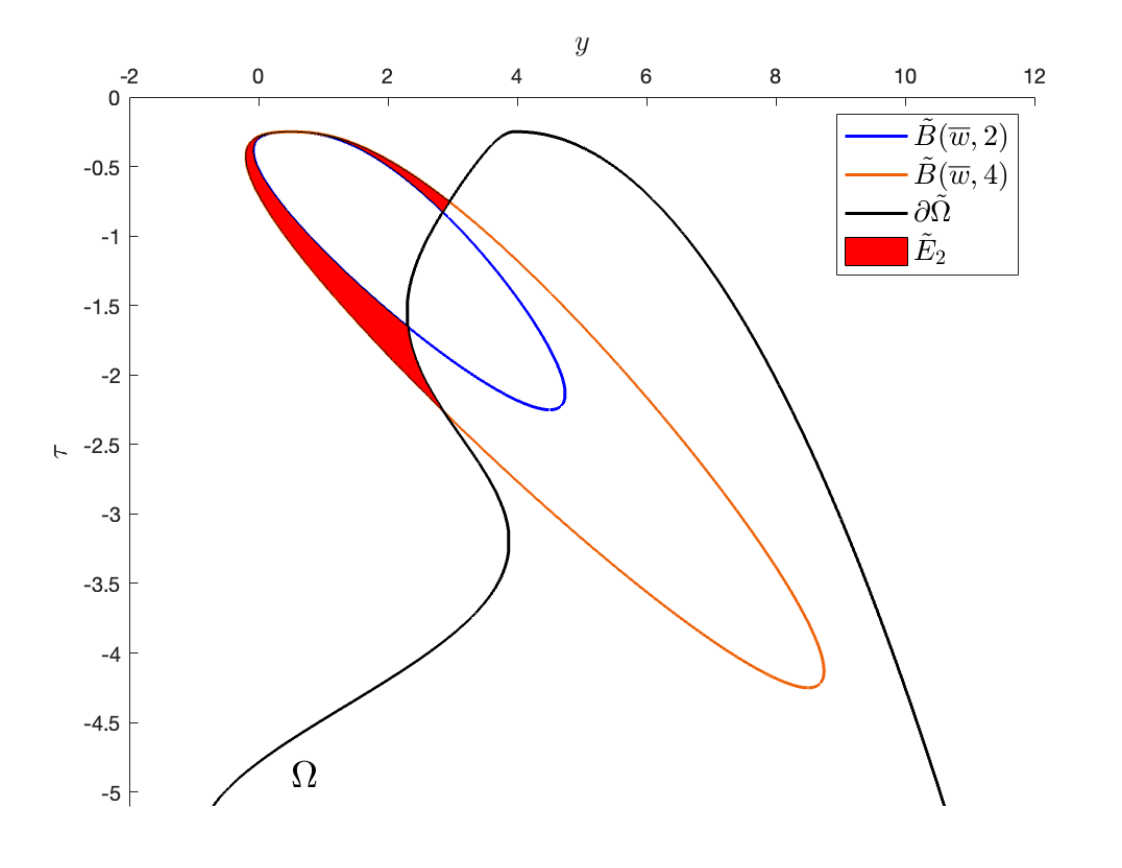}
    \caption{Sketch of $\tilde{h}$-heat balls $\tilde{B}(\bar{w},2), \tilde{B}(\bar{w},4)$, $\tilde{h}$-heat shell $\tilde{D}(\bar{w},4)$ and $\tilde{E}_2$ for a given open set $\tilde{\Omega}$ for $N=1$ and $\gamma = 1$}
\end{figure}

The solution of the {\bf Problem} $\tilde{\mathcal{A}}_\gamma$ for arbitrary open set $\tilde{\Omega}\subset\mathbb{R}^{N+1}_-$ reads:
 \begin{theorem}\label{wienertest2}
$\tilde{u}_*\equiv 0$ or $\tilde{u}_*>0$, that is to say the fundamental singularity at infinity is removable or non-removable
according as the series
\begin{equation}\label{wtest2}
\sum\limits_n 2^{-\frac{nN}{2}} C_{\tilde{h}}(\tilde{E}_n)
\end{equation} 
diverges or converges. An equivalent characterization is valid if the series \eqref{wtest2} is replaced with
\[ \sum\limits_n\lambda^{-n}C_{\tilde{h}}\Big ( \Big \{ \lambda^{-n+1}\geq \frac{F(\bar{w}-w)}{\tilde{h}(\bar{w})\tilde{h}_*(w)}\geq\lambda^{-n}\Big \} \Big )  \ \text{for any} \ \lambda>1\]

\end{theorem}

Remarkably, in the particular case $\gamma = 0$, the criterion \eqref{wtest2} is equivalent to the Wiener criterion at $\infty$ proved in \cite{Abdulla2}. 
In this particular case, $\tilde{h}\equiv 1$, and the concept $\tilde{h}$-capacity from Definition 1.3 coincides with the classical concept of 
thermal capacity \cite{watson}. The $\tilde{h}$-{\it heat ball} in Definition 1.4 coincides with the heat ball
\begin{equation}\label{heatball} 
\tilde{B}(\overline{w},c)=\{w: \ F(\overline{z}-w)> (4\pi c)^{-N/2} \}
\end{equation}
of center $\overline{w}=(0,-1/4)$ and radius $c$; the $\tilde{h}$-{\it heat shell}  $\tilde{D}(\overline{z},c)$ is a closed set of points between the level surfaces $F=(4\pi c)^{-N/2}$ and 
$F=(2\pi c)^{-N/2}$ of the fundamental solution of the backward heat equation with pole at $\overline{w}$:
\[ \tilde{D}(\overline{w},c)=\{w: (2\pi c)^{-N/2}\geq F(\overline{w}-w) \geq (4\pi c)^{-N/2}\} \cup \{\overline{w}\} \]
The divergence of the series \eqref{wtest2} is equivalent to that of the following series
\begin{equation}\label{A2wtest2}
 \sum\limits_n 2^{-\frac{nN}{2}} cap (E_n) 
 \end{equation}
with $E_n=\Omega^c\cap \tilde{D}(\overline{w},2^n)$. In \cite{Abdulla2} it is proved that the divergence of the series \eqref{A2wtest2} is a necessary and sufficient condition for the uniqueness of a bounded solution to the parabolic Dirichlet problem in an arbitrary open set $\tilde{\Omega}\subset\mathbb{R}^{N+1}$. The Problem~$\tilde{\mathcal {A}}_0$ was formulated by Kolmogorov in 1928 in the seminar on probability theory at Moscow State University, in the particular case with
$\Omega=\{|x|<f(t), -\infty<t<0\}\subset \mathbb{R}^2$. Kolmogorov's motivation for posing this problem was a connection to the probabilistic problem of finding the asymptotic behavior at infinity of the standard Brownian path. The Kolmogorov Problem in a one-dimensional setting was solved by Petrovsky in 1935, and the celebrated result is called the Kolmogorov-Petrovski test in the probabilistic literature \cite{Petrovski}. 

The full solution of the Kolmogorov Problem for arbitrary open sets $\Omega$ (or Problem~$\tilde{\mathcal{A}}_0$) is presented in \cite{Abdulla2}. A new concept of regularity or irregularity of $\infty$ is introduced according to whether the parabolic measure of $\infty$ is null or positive, and the necessary and sufficient condition for the Problem~$\tilde{\mathcal{A}}_0$ is proved in terms of the Wiener-type criterion for the regularity of $\infty$. 

In the probabilistic context, the formulated problems $\mathcal{A}_\gamma$ and $\tilde{\mathcal {A}}_\gamma$ are generalizations of the Kolmogorov problem to establish asymptotic laws for the $h$-{\it Brownian processes} \cite{Doob} (see Section~\ref{hBrownian}). \\

\section{Formulation of Problems}\label{description of results,historical remarks}

Being a generalization of the Kolmogorov problem, the Problems ${\mathcal A}_\gamma$, $\tilde{{\mathcal A}}_\gamma$, and their solution expressed in Theorems~\ref{wienertest} and ~\ref{wienertest2} have far-reaching measure-theoretical, topological, and probabilistic implications in Analysis, PDEs, and Potential Theory. The goal of this section is to formulate five outstanding problems equivalent to the Problems ${\mathcal A}_\gamma$ and $\tilde{{\mathcal A}}_\gamma$. Since the problems ${\mathcal A}_\gamma$ and $\tilde{{\mathcal A}}_\gamma$ are equivalent via the Appell transformation, without loss of generality, we are going to mainly focus on the formulations in the framework of the Problem ${\mathcal A}_\gamma$. The equivalent formulations are pursued in the framework of the Problem $\tilde{{\mathcal A}}_\gamma$ by replacing the triple $(\mathbb{R}^{N+1}_+, \Omega, h)$ with singularity point at $\mathcal{O}$ through the triple $(\mathbb{R}^{N+1}_-, \tilde{\Omega}, \tilde{h})$ with the singularity point at $\infty$ respectively.

%{\large 
\subsection{Unique Solvability of the Singular Parabolic Dirichlet Problem}\label{E:1:1}
%%%%%%%%%%%%
 Consider a {\it singular} {\bf parabolic Dirichlet problem}(PDP) \eqref{sDproblem}.
The solution of the PDP can be constructed by Perron's method (or the method by Perron, Wiener, Brelot, and Bauer)\cite{Doob,watson}. Let us introduce some necessary terminology. 

We will often write a typical point $z\in\mathbb{R}^{N+1}$ as $z=(x,t), x\in \mathbb{R}^N, t\in \mathbb{R}$. A smooth solution of the heat equation is called a {\it parabolic function}. A bounded open set $U\subset \mathbb{R}^{N+1}$ is {\it regular} if for each continuous $f: \partial U \to\mathbb{R}$ there exists one (and only one) parabolic function $H^U_f: U\to\mathbb{R}$, such that
\[ \lim\limits_{\substack{z\to w \\ w\in \partial U}}H_f^U(z)=f(w), \ w\in \partial U.      \] 
 A function $u$ is called a {\it superparabolic} in $\Omega$ if it satisfies the following conditions:
\begin{enumerate}
\item $-\infty<u\le+\infty$, $u < +\infty$ on a dense subset of $\Omega$; 
\item $u$ is lower semicontinuous (l.s.c.);
\item for each regular open set $U \subset \Omega$ and each parabolic function $v\in C(\bar{U})$,
the inequality $u\ge v$ on $\partial U$ implies $u \ge v$ in $U$.    
\end{enumerate}
A function $u$ is called a {\it subparabolic} if $-u$ is a superparabolic. 

A function $u=v/h$ is called a $h$-{\it parabolic}, $h$-{\it superparabolic}, or $h$-{\it subparabolic} in $\Omega$ if $v$ is {\it parabolic}, {\it superparabolic}, or {\it subparabolic} \cite{Doob}.

We use the notation
$\mathcal{S}_h(\Omega)$ for a class of all $h$-superparabolic functions in $\Omega$.
Similarly, the class of all $h$-subparabolic functions
in $\Omega$ is $-\mathcal{S}_h(\Omega)$.

Given boundary function $f$ on $\partial \Omega$, consider a {\bf $h$-parabolic Dirichlet
problem ($h$-PDP): find $h$-parabolic function $u$ in $\Omega$ such that }
\begin{equation}\label{Eq:W:1:6}
u=f \quad\text{on}\ \partial \Omega
\end{equation}
It is easy to see that $h$-parabolic function $u=\frac{v}{h}$ is a bounded solution of the
$h$-PDP if and only if $v$ is a solution of the PDP \eqref{sDproblem}.

Assuming for a moment
that $f\in C(\partial \Omega)$, the {\it generalized upper (or lower) solution} of the $h$-PDP is defined as
\begin{equation}\label{Eq:W:1:7}
^h \bar{H}_f^{\Omega} \equiv \inf \{u\in \mathcal{S}_h(\Omega): \liminf_{z\to w, z\in \Omega} u \ge f(w) \ \quad\text{for all}\ w \in \partial\Omega \}
\end{equation}
\begin{equation}\label{Eq:W:1:8}
^h \b{H}_f^{\Omega} \equiv \sup \{u\in -\mathcal{S}_h(\Omega): \limsup_{z\to w, z\in \Omega} u \le f(w) \ \quad\text{for all}\ w \in \partial \Omega \}
\end{equation}
%%%%%new%%%%%%%
The class of functions defined in (\ref{Eq:W:1:7}) (or in (\ref{Eq:W:1:8})) is called {\it upper class} (or {\it lower class}) of the $h$-PDP.
According to classical potential theory \cite{Doob}, $f$ is a {\it $h$-resolutive boundary function} in the sense that
\[ ^h\bar{H}_f^{\Omega} \equiv \ ^h\b{H}_f^{\Omega} \equiv \ ^hH_f^{\Omega}. \]
The indicator function of any Borel measurable boundary subset, and equivalently any bounded Borel measurable boundary function, is resolutive. Being $h$-parabolic in $\Omega$, $^hH_f^{\Omega}$ is called a generalized solution of the $h$-PDP for $f$. 
The generalized solution is unique by construction. It is essential to note that
the construction of the generalized solution
is accomplished by prescribing the behavior of
the solution at $\mathcal{O}$. 

Equivalently, we can define a generalized solution of the PDP \eqref{sDproblem}:

\begin{definition}\label{def1}
Let $g:\partial\Omega\to \mathbb{R}$ be a boundary function, such that $g/h$ is a bounded Borel measurable. Then $g$ is called a resolutive boundary function for the PDP  \eqref{sDproblem}, if $f=g/h$ (extended to $\mathcal{O}$) is $h$-resolutive for the $h$-PDP. The function
\begin{equation}\label{newPWB}
H_g^\Omega := h \ ^hH_f^{\Omega}
\end{equation}
is called a generalized solution of the PDP \eqref{sDproblem}.
\end{definition}
In particular, for the resolutive boundary function $g$, we have
\[ \overline{H}_g^{\Omega} \equiv \ \underline{H}_g^{\Omega} \equiv \ H_g^{\Omega}, \]
where
\begin{equation}\label{Eq:W:1:7add}
\overline{H}_g^{\Omega} \equiv \inf \{u\in \mathcal{S}(\Omega): \liminf_{z\to w, z\in \Omega} \frac{u}{h} \ge \limsup_{z\to w, z\in \partial\Omega}\frac{g}{h} \ \quad\text{for all}\ w \in \partial\Omega \}
\end{equation}
\begin{equation}\label{Eq:W:1:8add}
\underline{H}_g^{\Omega} \equiv \sup \{u\in -\mathcal{S}(\Omega): \limsup_{z\to w, z\in \Omega} \frac{u}{h} \le  \liminf_{z\to w, z\in \Omega}\frac{g}{h} \ \quad\text{for all}\ w \in \partial \Omega \}
\end{equation}

Again, note that the unique solution $H_g^\Omega$ of the PDP \eqref{sDproblem} is constructed by prescribing the behavior 
of the ratio $H_g^{\Omega} /h$ at $\mathcal{O}$.

The elegant theory, 
while identifying a class of unique solvability, 
leaves the following questions open:
\begin{itemize}
\item Would a unique solution of the $h$-DP still exist if its 
limit at $\mathcal{O}$ were not specified? That is, could it be 
that the solutions would pick up the ``boundary value" $f(\mathcal{O})$
without being required? Equivalently, would a unique solution of the PDP \eqref{sDproblem} still exist if the limit of the ratio of solution to $h$ at $\mathcal{O}$ is not prescribed?  In particular, is the fundamental singularity at $\mathcal{O}$
removable?  
\item What if the boundary datum $f$ (or $g/h$) on $\partial \Omega$, while being continuous at $\partial\Omega\setminus\{\mathcal{O}\}$, 
does not have a limit at $\mathcal{O}$, for example, it 
exhibits bounded oscillations. Is the $h$-PDP (or PDP \eqref{sDproblem}) uniquely solvable?
\end{itemize}
\begin{example}\label{example1}
 Let $\Omega=\mathbb{R}^{N+1}_+$.
It is easy to see that the boundary of $\mathbb{R}^{N+1}_+$ is $h$-resolutive and the only possible solutions of the $h$-PDP in $\mathbb{R}^{N+1}_+$ are constants. Precisely, the unique solution of the $h$-PDP is identical with the constant $f(\mathcal{O})$. Indeed, for arbitrary $\epsilon >0$, the function 
\[ u(\cdot)=f(\mathcal{O})+\frac{\epsilon}{h(\cdot)} \ \ \Big (\quad\text{or}\ v(\cdot)=f(\mathcal{O})-\frac{\epsilon}{h(\cdot)} \Big ) \]
is in the upper class (or lower class) for $h$-PDP in $\mathbb{R}^{N+1}_+$ for $f$. Hence,
\[ f(\mathcal{O})-\frac{\epsilon}{h(\cdot)}\leq \ ^h \b{H}_f^{\mathbb{R}^{N+1}_+}(\cdot) \leq \ ^h \bar{H}_f^{\mathbb{R}^{N+1}_+}(\cdot) \leq f(\mathcal{O})+\frac{\epsilon}{h(\cdot)} \]
Since $\epsilon>0$ is arbitrary, the assertion follows. Equivalently, all possible solutions of the PDP \eqref{sDproblem}in $\mathbb{R}^{N+1}_+$ are constant multiples of $h$, and the unique solution is identified by prescribing the ratio $u/h$ at $\mathcal{O}$. 
\end{example}
Example~\ref{example1} demonstrates that if $\Omega= \mathbb{R}^{N+1}_+$, the answer is negative and 
arbitrary constant $C$ is a solution of the $h$-PDP, $Ch$ is a solution of the PDP \eqref{sDproblem}, and the fundamental singularity at $\mathcal{O}$ is not removable. 

The positive answer to these fundamental questions is possible if 
$\Omega$ is not too sparse, or equivalently $\Omega^c \cap \mathbb{R}^{N+1}_+$ is
not too thin near $\mathcal{O}$. 
The principal purpose of this paper is to prove  
the necessary and sufficient condition for the non-thinness of $\Omega^c \cap \mathbb{R}^{N+1}_+$ near $\mathcal{O}$ 
which is equivalent to the uniqueness of the solution of the
$h$-PDP (or PDP (\eqref{sDproblem}) without specification of the boundary function (or ratio of the boundary function to $h$) at $\mathcal{O}$.

Furthermore, given {\it bounded Borel measurable} function $f=g/h: \partial\Omega\setminus\{\mathcal{O}\} \rightarrow \mathbb{R}$,
we fix an arbitrary finite real number $\bar{f}$, and extend a function
$f$ as $f(\mathcal{O})=\bar{f}$. The extended function
is a bounded Borel measurable on $\partial\Omega$ and
there exists a unique solution $^hH_f^{\Omega}$ of the $h$-PDP, and the unique solution of the PDP \eqref{sDproblem}
is given by \eqref{newPWB}. The major question now becomes:\\

{\it {\bf Problem 1:} How many bounded solutions do we have, or does the constructed
solution depend on} $\bar{f}$ ? \\

Similarly, in the context of the singular PDP \eqref{sDproblem-} the generalized solution is defined by the formula
\begin{equation}\label{tildenewPWB}
H_g^{\tilde{\Omega}} := \tilde{h} \ ^{\tilde{h}}H_f^{\tilde{\Omega}},
\end{equation}
where $\ ^{\tilde{h}}H_f^{\tilde{\Omega}}$ is a Perron solution of the coresponding $\tilde{h}$-PDP with the boundary function $f=g/\tilde{h}$ extended to boundary point $\infty\in\partial\tilde{\Omega}$ as $f(\infty)=\bar{f}$ with arbitrary fixed $\bar{f}\in\mathbb{R}$. 
In particular, boundary function $g:\partial\tilde{\Omega}\to \mathbb{R}$, with $g/\tilde{h}$ being bounded Borel measurable, is resolutive with
\[ \overline{H}_g^{\tilde{\Omega}} \equiv \ \underline{H}_g^{\tilde{\Omega}} \equiv \ H_g^{\tilde{\Omega}}, \]
where
\begin{equation}\label{Eq:W:1:7addtilde}
\overline{H}_g^{\tilde{\Omega}} \equiv \inf \Big \{u\in \mathcal{S}(\tilde{\Omega}): \liminf_{z\to w, z\in \tilde{\Omega}} \frac{u}{\tilde{h}} \ge \limsup_{z\to w, z\in \partial\tilde{\Omega}}\frac{g}{\tilde{h}} \ \quad\text{for all}\ w \in \partial\tilde{\Omega} \Big \}
\end{equation}
\begin{equation}\label{Eq:W:1:8addtilde}
\underline{H}_g^{\tilde{\Omega}} \equiv \sup \Big \{u\in -\mathcal{S}(\tilde{\Omega}): \limsup_{z\to w, z\in \tilde{\Omega}} \frac{u}{\tilde{h}} \le  \liminf_{z\to w, z\in \tilde{\Omega}}\frac{g}{\tilde{h}} \ \quad\text{for all}\ w \in \partial \tilde{\Omega} \Big \}
\end{equation}
The solution $u=H_g^{\tilde{\Omega}}$ of the singular PDP \eqref{sDproblem-} is unique by prescribing the behavior of $u/\tilde{h}$ via relations in \eqref{Eq:W:1:7addtilde},\eqref{Eq:W:1:8addtilde} at every bondary point $w\in\partial\tilde{\Omega}$, including at the boundary singularity point $\infty$. The major open problem is, as it is formulated in {\bf Problem 1}, whether uniqueness holds without prescribing the behavior of $u/\tilde{h}$ at $\infty$.

\subsection{Characterization of the $h$-Parabolic Measure of Singularity Point}\label{E:1:2}
For a given boundary Borel subset $A\subset \partial \Omega$, denote the indicator function
of $A$ as $1_A$. Indicator functions of the Borel measurable subsets of $\partial\Omega$ are 
resolutive (\cite{Doob}). {\it $h$-Parabolic measure} of the boundary Borel subset $A$ is defined as (\cite{Doob}):
\[ \mu_{\Omega}^h(z,A)= \ ^hH_{1_A}^\Omega(z), \]
where $z\in\Omega$ is a reference point. It is said that $A$ is an $h$-parabolic measure null set if $\mu_{\Omega}(\cdot,A)$
vanishes identically in $\Omega$. If this is not the case, $A$ is a set of
positive $h$-parabolic measure.
In particular, the $h$-parabolic measure of $\{\mathcal{O}\}$ is well defined:
\[ \mu_{\Omega}^h(\cdot,\{\mathcal{O}\})= \ ^hH_{1_{\{\mathcal{O}\}}}^\Omega(\cdot). \]
The following formula is true for the solution  $^h H_f^{\Omega}$ of the $h$-PDP \cite{Doob}:
\begin{equation}\label{Eq:W:1:8f}
^h H_f^{\Omega}(z)=\int\limits_{\partial\Omega}f(w)\mu_\Omega^h(z,dw), \ \ z\in\Omega
\end{equation}
Since $f$ is extended to $\{\mathcal{O}\}$ as $f(\mathcal{O})=\bar{f}$, we have
\begin{equation}\label{Eq:W:1:8f-1}
^h H_f^{\Omega}(z)=\int\limits_{\partial\Omega\setminus\{\mathcal{O}\}}f(w)\mu_\Omega^h(z,dw)+\bar{f}\cdot  \ ^hH_{1_{\{\mathcal{O}\}}}^\Omega(z), \ \ z\in\Omega
\end{equation}
This formula implies that the uniqueness of the solution to the $h$-PDP without prescribing the behavior of the solution at the singularity point $\mathcal{O}$, that is to say, the
independence of $^h H_f^{\Omega}$ on $\bar{f}$ is equivalent to whether or not $\mathcal{O}$ is an $h$-parabolic measure null set. Equivalently, according to the formula \eqref{newPWB}, the following formula is true for the unique solution of the PDP \eqref{sDproblem}:
\begin{equation}\label{Eq:W:1:8h}
H_g^{\Omega}(z)=h(z)\int\limits_{\partial\Omega}\frac{g(w)}{h(w)}\mu_\Omega^h(z,dw), \ \ z\in\Omega
\end{equation}
Splitting the integral as in \eqref{Eq:W:1:8f-1} we have
\begin{equation}\label{Eq:W:1:8h-1}
H_g^{\Omega}(z)=h(z)\int\limits_{\partial\Omega\setminus\{\mathcal{O}\}}\frac{g(w)}{h(w)}\mu_\Omega^h(z,dw)+\bar{f}h(z) \ ^hH_{1_{\{\mathcal{O}\}}}^\Omega(z), \ \ z\in\Omega.
\end{equation}
where $\bar{f}$ is a prescribed limit value of $H_g^{\Omega} /h$ at $\mathcal{O}$. Similar to its counterpart \eqref{Eq:W:1:8f-1}, the formula \eqref{Eq:W:1:8h-1} demonstrates that the uniqueness of the solution $u$ of the PDP \eqref{sDproblem} without prescribing $u /h$ at $\mathcal{O}$ is equivalent to whether or not $\mathcal{O}$ is an $h$-parabolic measure null set.

Hence, the following problem is the measure-theoretical counterpart of Problem 1:\\

{\it {\bf Problem 2:} Given $\Omega$, is the $h$-parabolic measure 
of $\{\mathcal{O}\}$ null or positive} ?\\

From the Example~\ref{example1} demonstrated above it follows that in the particular case with $\Omega=\mathbb{R}^{N+1}_+$, we have
\begin{equation}\label{Eq:W:1:8a}
 \mu_{\mathbb{R}^{N+1}_+}^h(\cdot,\{\mathcal{O}\})\equiv 1, \ \ \mu_{\mathbb{R}^{N+1}_+}^h(\cdot,\partial \mathbb{R}^{N+1}_+ -\{\mathcal{O}\})\equiv 0.
\end{equation}
\begin{example}\label{example2} 
For arbitrary $c>0$, consider a domain bounded by the level set of $h$:
\[ \Omega=\{h>c'\}=\{z: |x-\gamma|^2<-2Nt\log\frac{t}{c}, \ 0<t<c\}, \]
where $c'=(4\pi c)^{-\frac{N}{2}}$. It is easy to see that the $h$-parabolic measure of $\{\mathcal{O}\}$ is positive, and we have
\begin{equation}\label{ex2}
\ ^hH_{1_{\{\mathcal{O}\}}}^\Omega(z)=\frac{h(z)-c'}{h(z)}, \ z\in \Omega.
\end{equation}
\end{example}
Both {\bf Problem 1 and 2  are equivalent to the Problem $\mathcal{A}_\gamma$ formulated in Section~\ref{prelude}}. The connection follows from the following formula:
\begin{equation}\label{u*parmeas}
u_*(z)=h(z) \ ^hH_{1_{\{\mathcal{O}\}}}^\Omega(z), \ z\in\Omega
\end{equation}
To establish \eqref{u*parmeas}, first note that the $h$-parabolic function $u_*/h$ is in the lower class of the Perron's solution $^hH_{1_{\{\mathcal{O}\}}}^\Omega$,
which implies that
\begin{equation}\label{u*parmeas1}
\frac{u_*(z)}{h(z)}\leq \  ^hH_{1_{\{\mathcal{O}\}}}^\Omega(z), \ z\in \Omega.
\end{equation}
Moreover, $^hH_{1_{\{\mathcal{O}\}}}^\Omega$ itself is in the lower class of Perron's solution $u_n/h$ of the $h$-PDP in $\Omega_n$ with boundary function
$1_{\partial\Omega_n\cap\{t=n^{-1}\}}$, where $u_n$ is a solution of PDP \eqref{u_n}. Therefore, we have
\begin{equation}\label{u*parmeas2}
 \  ^hH_{1_{\{\mathcal{O}\}}}^\Omega(z)\leq \frac{u_n(z)}{h(z)}, \ z\in \Omega_n.
\end{equation}
passing to the limit as $n\to\infty$, from \eqref{ulim}, \eqref{u*parmeas1} and \eqref{u*parmeas2}, \eqref{u*parmeas} follows.

In light of the measure-theoretical counterpart of the removability of the fundamental singularity, we introduce a concept of $h$-regularity of the 
boundary point $\mathcal{O}$. 
\begin{definition}\label{def2}
$\mathcal{O}$ is said to be $h$-{\it regular} for
$\Omega$ if it is an $h$-parabolic measure null set. Conversely, $\mathcal{O}$ is $h$-{\it irregular} if it has a positive $h$-parabolic measure.
\end{definition}
Hence, Theorem~\ref{wienertest} establishes a criterion for the removability of the fundamental singularity in terms of the necessary and sufficient condition for the $h$-{\it regularity} of $\mathcal{O}$.

Similarly, in the context of the singular PDP \eqref{sDproblem-}, and corresponding $\tilde{h}$-PDP, we have the formulae analogous to \eqref{Eq:W:1:8f-1}, \eqref{Eq:W:1:8h-1} and \eqref{u*parmeas}:
\begin{gather}
^{\tilde{h}} H_f^{\tilde{\Omega}}(z)=\int\limits_{\partial\tilde{\Omega}\setminus\{\infty\}}f(w)\mu_{\tilde{\Omega}}^{\tilde{h}}(z,dw)+\bar{f}\cdot  \ ^{\tilde{h}}H_{1_{\{\infty\}}}^{\tilde{\Omega}}(z), \ \ z\in\tilde{\Omega}\label{Eq:W:1:8f-1tilde}\\
H_g^{\tilde{\Omega}}(z)=\tilde{h}(z)\int\limits_{\partial\tilde{\Omega}\setminus\{\infty\}}\frac{g(w)}{h(w)}\mu_{\tilde{\Omega}}^{\tilde{h}}(z,dw)+\bar{f}\tilde{h}(z) \ ^{\tilde{h}}H_{1_{\{\infty\}}}^{\tilde{\Omega}}(z), \ \ z\in\tilde{\Omega}\label{Eq:W:1:8h-1tilde}\\
\tilde{u}_*(z)=\tilde{h}(z) \ ^{\tilde{h}}H_{1_{\{\infty\}}}^{\tilde{\Omega}}(z), \ z\in\tilde{\Omega}.\label{tildeu*parmeas}
\end{gather}
We introduce a concept of $\tilde{h}$-regularity of the boundary point $\infty$ for $\tilde{\Omega}\subset \mathbb{R}^{N+1}_-$.
\begin{definition}\label{def3}
 $\infty$ is said to be $\tilde{h}$-{\it regular} for
$\tilde{\Omega}$ if it is an $\tilde{h}$-parabolic measure null set. Conversely, $\infty$ is $\tilde{h}$-{\it irregular} if it has a positive $\tilde{h}$-parabolic measure.
\end{definition}
In fact, in the particular case with $\gamma=0, \tilde{h}\equiv 1$, it coincides with the concept of regularity of $\infty$ introduced in \cite{Abdulla2}. Theorem~\ref{wienertest2} presents a criterion for the removability of the fundamental singularity at $\infty$ in terms of the necessary and sufficient condition for the $\tilde{h}$-{\it regularity} of $\infty$. 

\subsection{Boundary Regularity in Singular Dirichlet Problem}\label{E:1:3}
The notion of the $h$-regularity of $\mathcal{O}$ is, in particular, relates to 
the notion of continuity of the solution to the $h$-PDP at $\mathcal{O}$. 

{\it {\bf Problem 3:} Given $\Omega$,
determine whether or not}
\begin{equation}\label{boundreg}
\liminf_{z\to \mathcal{O}, z\in \partial\Omega} f \le \liminf_{z\to \mathcal{O}, z\in \Omega}  \ ^hH_f^{\Omega} \le \limsup_{z\to \mathcal{O}, z\in \Omega}  \ ^hH_f^{\Omega} \le \limsup_{z\to \mathcal{O}, z\in \partial\Omega} f,
\end{equation}
for all bounded $f\in C(\partial \Omega\setminus\{\mathcal{O}\})$.

Note that if $f$ has a limit at $\mathcal{O}$, (\ref{boundreg}) simply means that 
the solution $^hH_f^{\Omega}$ is continuous at $\mathcal{O}$.

The equivalent problem in the context of the PDP \eqref{sDproblem} is the following:

{\it {\bf Problem 3$^\prime$:} Given $\Omega$,
whether or not
\begin{equation}\label{hboundreg}
\liminf_{z\to \mathcal{O}, z\in \partial\Omega} \frac{g}{h} \le \liminf_{z\to \mathcal{O}, z\in \Omega}  \ \frac{H_g^{\Omega}}{h} \le \limsup_{z\to \mathcal{O}, z\in \Omega}  \ \frac{H_g^{\Omega}}{h} \le \limsup_{z\to \mathcal{O}, z\in \partial\Omega} \frac{g}{h},
\end{equation}
for all $g$ such that $\frac{g}{h}\in C(\partial \Omega\setminus\{\mathcal{O}\})$ and bounded.}

In particular, if  $g/h$ has a limit at $\mathcal{O}$, (\ref{hboundreg}) means that 
the limit of the ratio $H_g^{\Omega}/h$ at $\mathcal{O}$ exists and equal to the limit of $g/h$.

\begin{definition}\label{defbr}
Singular PDP \eqref{sDproblem} (and corresponding $h$-PDP) is said to be regular at $\mathcal{O}$ if \eqref{boundreg} and \eqref{hboundreg} are satisfied. It is said to be irregular at $\mathcal{O}$ otherwise.
\end{definition}

Similarly, we introduce the concept of regularity at $\infty$ of the singular PDP \eqref{sDproblem-} (and corresponding $\tilde{h}$-PDP):

\begin{definition}\label{defbrtilde}
Singular PDP \eqref{sDproblem-} (and corresponding $\tilde{h}$-PDP) is said to be regular at $\infty$ if the following conditions are satisfied; 
\begin{gather}
\liminf_{z\to \infty, z\in \partial\tilde{\Omega}} f \le \liminf_{z\to \infty, z\in \Omega}  \ ^{\tilde{h}}H_f^{\tilde{\Omega}} \le \limsup_{z\to \infty, z\in \tilde{\Omega}}  \ ^{\tilde{h}}H_f^{\tilde{\Omega}} \le \limsup_{z\to \infty, z\in \partial\tilde{\Omega}} f,\label{tildehboundreg}\\
\forall \ \text{bounded} \ f\in C(\partial \tilde{\Omega})\nonumber\\
\liminf_{z\to \infty, z\in \partial\tilde{\Omega}} \frac{g}{\tilde{h}} \le \liminf_{z\to \infty, z\in \tilde{\Omega}}  \ \frac{H_g^{\tilde{\Omega}}}{\tilde{h}} \le \limsup_{z\to \infty, z\in \tilde{\Omega}}  \ \frac{H_g^{\tilde{\Omega}}}{\tilde{h}} \le \limsup_{z\to \infty, z\in \partial\tilde{\Omega}} \frac{g}{\tilde{h}},\label{Atildehboundreg}\\
\forall \ \text{bounded} \ \frac{g}{\tilde{h}} \in C(\partial \tilde{\Omega})\nonumber
\end{gather}
It is said to be irregular at $\infty$ otherwise.
\end{definition}

\subsection{Thinness in Parabolic Minimal-Fine Topology}\label{E:1:4}
The notion of $h$-{\it regularity} of $\mathcal{O}$ introduced in Definition~\ref{def2} fits naturally in the framework of parabolic minimal-fine topology. Recall that parabolic fine topology
is the coarsest topology of $\mathbb{R}^{N+1}$ which makes every superparabolic function continuous
\cite{Doob}.  Parabolic fine topology is finer than the Euclidean topology. It is well-known that there is an elegant connection
between the problem of finding the structure of the neighborhood base in
parabolic fine topology and the problem of the regularity of finite boundary points.
Namely, given open set $\Omega \subset \mathbb{R}^{N+1}$, its finite boundary point $x_0$
is irregular if and only if $\Omega$ is a deleted neigborhood of $x_0$
in parabolic-fine topology \cite{Doob}. Equivalently, $\Omega^c$ is called parabolic thin at 
$x_0$. For arbitrary open sets with non-compact boundaries, a similar connection for the point at $\infty$ was lacking. The new concept of the regularity
of the point at $\infty$ introduced in \cite{Abdulla2}, in particular, revealed a similar connection. Wiener test for the regularity of the point at $\infty$ provides 
the criterion for the uniqueness of the bounded solution of the parabolic Dirichlet problem in any open set with non-compact boundary, and equivalently 
provides a criterion for the characterization of the neighborhood base of the point at $\infty$ in parabolic fine topology \cite{Abdulla2}.

Parabolic minimal-fine topology in $\mathbb{R}^{N+1}_+\cup \{\mathcal{O}\}$ is an an extension of the parabolic fine topology of $\mathbb{R}^{N+1}_+$ to boundary point $\mathcal{O}$. Positive parabolic function $h$ is minimal parabolic in $\mathbb{R}^{N+1}_+$, in the sense that it dominates there no positive parabolic function except for its own constant submultiples. $\mathcal{O}$ is a minimal Martin boundary point of $\mathbb{R}^{N+1}_+$ and $h$ is an associated minimal Martin boundary function
with pole at $\mathcal{O}$. 
Parabolic minimal-fine topology has as a relative topology on $\mathbb{R}^{N+1}_+$, the parabolic fine topology \cite{Doob}. \\

\begin{definition}\label{mthinness}
Subset $E\subset \mathbb{R}^{N+1}_+$ is parabolic minimally thin at $\mathcal{O}$ if
\begin{equation}\label{Eq:W:1:10a}
\lim_{\delta\downarrow 0} \ ^h\hat{R}^{E\cap \{t\leq \delta\}}_1(z)=0, \ z\in \mathbb{R}^{N+1}_+
\end{equation}
where $\ ^hR^{E}_1$ be a $h$-{\it reduction} of $1$ on E
\[ ^hR^{E}_1(z)=\inf \{v(z): v \quad\text{is $h$-superparabolic in} \ \mathbb{R}^{N+1}_+, v\ge 0, v\ge 1 \quad\text{on} \ E \} \]
and $\ ^h\hat{R}^{E}_1$ is its lower-semicontinuous regularization, or smoothed $h$-{\it reduction} of $1$ on $E$:
\[ \ ^h\hat{R}^{E}_1(z)=\lim_{\epsilon \to 0} \inf_{w\in B_\epsilon(z)} \ ^hR^{E}_1(w) \]
\end{definition}

The following is the equivalent definition of parabolic minimal thinness at $\mathcal{O}$.\\

\begin{definition}\label{mthinnessalt}
Subset $E\subset \mathbb{R}^{N+1}_+$ is parabolic minimally thin at $\mathcal{O}$ if
\begin{equation}\label{minimalthinnessdefinition}
^h\hat{R}^{E\cap \{t\leq \delta\}}_1(z)\not\equiv 1, \ z\in \mathbb{R}^{N+1}_+
\end{equation}
for some (and equivalently for all) $\delta >0$. 
\end{definition}

A point $\mathcal{O}$ is a parabolic minimal-fine limit point of a set $E\subset \mathbb{R}^{N+1}_+$ if $E$ is not parabolic minimally thin at $\mathcal{O}$. We write this fact as $\mathcal{O} \in E^{pmf}$, where $E^{pmf}$ is the set of parabolic minimal-fine limit points of $E$. If $E$ is parabolic minimally thin at $\mathcal{O}$, then equivalently, $E^c\cap \mathbb{R}^{N+1}_+$ is a deleted parabolic minimal-fine neighborhood of $\mathcal{O}$.  

We can now formulate the topological counterpart of Problems 1-3 \cite{AbdullaI}\\

{\it {\bf Problem 4:} Is the given open set $\Omega\subset \mathbb{R}^{N+1}_+$ a deleted neigborhood of $\mathcal{O}$
in parabolic minimal-fine topology? Equivalently, is $\Omega^c$ parabolic minimally thin at $\mathcal{O}$?
Or conversely, whether or not $\mathcal{O} \in (\Omega^c\cap \mathbb{R}^{N+1}_+)^{pmf}$.} \\

\begin{remark}\label{remarkhred1e1redh}
Parabolic minimal thinness can be equivalently expressed in terms of the {\it reduction} of $h$. Since $h$-{\it reduction} of $1$ and {\it reduction} of $h$ satisfy the relation
\begin{equation}\label{hred1redh}
^h\hat{R}^{E}_1\equiv \frac{\hat{R}^{E}_h}{h}, \ \ ^h R^{E}_1\equiv \frac{R^{E}_h}{h}
\end{equation} 
the requirements \eqref{Eq:W:1:10a} and \eqref{minimalthinnessdefinition} can be replaced with the following relations:
\begin{gather}
\lim_{\delta\downarrow 0} \ \hat{R}^{E\cap \{t\leq \delta\}}_h(z)=0, \ z\in \mathbb{R}^{N+1}_+, \label{Eq:W:1:10a1}\\
\hat{R}^{E\cap \{t\leq \delta\}}_h(z)\not\equiv h, \ z\in \mathbb{R}^{N+1}_+. \label{minimalthinnessdefinition1}
\end{gather}
\end{remark}
%%%%%%%%
%%%%%%%%%
To define the parabolic minimal-fine topology in $\mathbb{R}^{N+1}_-\cup \{\infty\}$, note that the class of minimal parabolic functions in $\mathbb{R}^{N+1}_-$ are positive multiples of 
$\tilde{h}$ for some point $\gamma \in \mathbb{R}^N$ \cite{Doob}. All the minimal parabolic functions in $\mathbb{R}^{N+1}_-$ has a pole at $\infty$. Otherwise speaking. the parabolic minimal-fine topology at $\infty$ in general is dependent on $\gamma$. Therefore, we have to introduce the concept of $\tilde{h}$-parabolic minimal-fine topology at $\infty$.
$\tilde{h}$-parabolic minimal-fine topology in $\mathbb{R}^{N+1}_-\cup \{\infty\}$ is an extension of the parabolic fine topology of $\mathbb{R}^{N+1}_-$ to boundary point $\infty$.  \\

\begin{definition}\label{mthinnessinfinity}
Subset $E\subset \mathbb{R}^{N+1}_-$ is $\tilde{h}$-parabolic minimally thin at $\infty$ if
\begin{equation}\label{Eq:W:1:10a1tilde}
\lim_{\delta\downarrow -\infty} \ ^{\tilde{h}}\hat{R}^{E\cap \{t\leq \delta\}}_1(z)=0, \ z\in \mathbb{R}^{N+1}_-
\end{equation}
where $\ ^{\tilde{h}}R^{E}_1$ be a $\tilde{h}$-{\it reduction} of $1$ on E
\[ ^{\tilde{h}}R^{E}_1(z)=\inf \{v(z): v \quad\text{is $\tilde{h}$-superparabolic in} \ \mathbb{R}^{N+1}_-, v\ge 0, v\ge 1 \quad\text{on} \ E \} \]
and $\ ^{\tilde{h}}\hat{R}^{E}_1$ is its lower-semicontinuous regularization, or smoothed $\tilde{h}$-{\it reduction} of $1$ on $E$:
\[ \ ^{\tilde{h}}\hat{R}^{E}_1(z)=\lim_{\epsilon \to 0} \inf_{w\in B_\epsilon(z)} \ ^{\tilde{h}}R^{E}_1(w) \]
\end{definition}

The following is the equivalent definition of $\tilde{h}$-parabolic minimal thinness at $\infty$.\\

\begin{definition}\label{mthinnessinfinityalt}
Subset $E\subset \mathbb{R}^{N+1}_-$ is $\tilde{h}$-parabolic minimally thin at $\infty$ if
\begin{equation}\label{tildeminimalthinnessdefinition}
^{\tilde{h}}\hat{R}^{E\cap \{t\leq \delta\}}_1(z)\not\equiv 1, \ z\in \mathbb{R}^{N+1}_-
\end{equation}
for some $\delta<0$.
\end{definition}

A point $\infty$ is $\tilde{h}$-parabolic minimal-fine limit point of a set $E\subset \mathbb{R}^{N+1}_-$ if $E$ is not $\tilde{h}$-parabolic minimally thin at $\infty$. We write this fact as $\infty \in E^{\tilde{h}pmf}$, where $E^{\tilde{h}pmf}$ is the set of $\tilde{h}$-parabolic minimal-fine limit points of $E$. If $E$ is $\tilde{h}$-parabolic minimally thin at $\infty$, then equivalently, $E^c\cap \mathbb{R}^{N+1}_-$ is a deleted $\tilde{h}$-parabolic minimal-fine neighborhood of $\infty$.  

Topological counterpart of Problems 1-3 formulated in \cite{AbdullaI} reads:\\

{\it {\bf Problem 4:} Is the given open set $\Omega\subset \mathbb{R}^{N+1}_-$ a deleted neigborhood of $\infty$
in $\tilde{h}$-parabolic minimal-fine topology? Equivalently, is $\Omega^c$ $\tilde{h}$-parabolic minimally thin at $\infty$?
Or conversely, whether or not $\infty \in (\Omega^c\cap \mathbb{R}^{N+1}_-)^{\tilde{h}pmf}$.} \\

\begin{remark}\label{hred1e1redh12}
$\tilde{h}$-parabolic minimal thinness can be equivalently expressed in terms of {\it reduction} of $\tilde{h}$. Since $\tilde{h}$-{\it reduction} of $1$ and {\it reduction} of $\tilde{h}$ satisfy the relation
\begin{equation}\label{hred1redh11}
^{\tilde{h}}\hat{R}^{E}_1\equiv \frac{\hat{R}^{E}_{\tilde{h}}}{\tilde{h}}, \ \ ^{\tilde{h}} R^{E}_1\equiv \frac{R^{E}_{\tilde{h}}}{\tilde{h}}
\end{equation} 
the requirements \eqref{Eq:W:1:10a1} and \eqref{tildeminimalthinnessdefinition} can be replaced with the following relations:
\begin{gather}
\lim_{\delta\downarrow -\infty} \ \hat{R}^{E\cap \{t\leq \delta\}}_{\tilde{h}}(z)=0, \ z\in \mathbb{R}^{N+1}_-, \label{Eq:W:1:10a2}\\
\hat{R}^{E\cap \{t\leq \delta\}}_{\tilde{h}}(z)\not\equiv \tilde{h}, \ z\in \mathbb{R}^{N+1}_-. %\label{minimalthinnessdefinition1~}
\end{gather}
\end{remark}
\subsection{Asymptotic Laws for the $h$-Brownian Motion}\label{hBrownian}
In the probabilistic context, Theorems~\ref{wienertest} and \ref{wienertest2} establishes an asymptotic laws for the conditional Brownian motion \cite{Doob}. 

Let $\{x(t)=(x_1(t),...,x_N(t)): t\geq 0, P_\bullet\}$ be an $N$-dimensional $h$-Brownian process, and $P_\bullet(B)$ is a probability of the event
$B$ as a function of the starting point $x(\tau)$ with $\tau>0$ \cite{Doob}. $h$-Brownian motion $x(t)$ from a point $x(\tau)$ is an almost surely continuous process whose sample functions never leave $\mathbb{R}^{N+1}_+$ and proceed downward, that is, in the direction of decreasing $t$. In fact, almost every path starting at $x(\tau)$ has a finite lifetime $\tau$ and tends to the boundary point $\mathcal{O}$ as $t\downarrow 0$ \cite{Doob}. Consider the event 
\[ B= [ \text{the set} \ \{0<t\leq \tau: x(t)\in \Omega^c\} \ \text{clusters to} \ 0] \] 
The probabilistic analog of Theorem~\ref{wienertest} states that
\[ P(B) =0 \ \text{or} \ 1  \ \text{according as} \ \sum\limits_n2^{-\frac{nN}{2}}C_h(E_n) < \ \text{or} \ =+\infty. \]

Similarly, Theorem~\ref{wienertest2} establishes an asymptotic law for the $\tilde{h}$-Brownian processes. Let $\{x(t)=(x_1(t),...,x_N(t)): t<0, P_\bullet\}$ be an $N$-dimensional $\tilde{h}$-Brownian process, and $P_\bullet(B)$ is a probability of the event
$B$ as a function of the starting point $x(\tau)$ with $\tau<0$ \cite{Doob}. $\tilde{h}$-Brownian motion $x(t)$ from a point $x(\tau)$ is an almost surely continuous process whose sample functions never leave $\mathbb{R}^{N+1}_-$ and proceed downward, that is, in the direction of decreasing $t$. Almost every path starting at $x(\tau)$ tends to the boundary point $\infty$ as $t\downarrow -\infty$ \cite{Doob}. Given open set $\Omega\subset\mathbb{R}^{N+1}_-$, consider the event 
\[ B= [ \ \text{the set} \ \{-\infty<t\leq \tau : x(t)\in \Omega^c\} \  \text{clusters to} \  -\infty \ ].  \]
The probabilistic analog of Theorem~\ref{wienertest2} states that
\[ P(B) =0 \ \text{or} \ 1  \ \text{according as} \  \sum\limits_n 2^{-\frac{nN}{2}} C_{\tilde{h}}(\tilde{E}_n) < \ \text{or} \ =+\infty. \]
%%%%%%%%%%%
%%%%%%%%%%%%%
Theorem~\ref{wienertest} (or Theorem~\ref{wienertest2}) expresses the solutions to equivalent
Problems 1-4 in terms of the Wiener-type criteria for the $h$-{\it regularity} of $\mathcal{O}$ (or $\tilde{h}$-{\it regularity} of $\infty$) respectively.

%%%%%%%%%%%%%
%%%%%%%%%%%%%%
\section{The Main Results}\label{mainresult}

We now reformulate the main results of Theorems~\ref{wienertest}, \ref{wienertest2} in a broader context as a solution of the equivalent Problems 1-4. 

\begin{theorem}\label{main theorem}
For arbitrary open set $\Omega\subset\mathbb{R}^{N+1}_+$ the following conditions are equivalent:
\begin{enumerate}
\item $\mathcal{O}$ is $h$-regular (or $h$-irregular).
\item Singular Parabolic Dirichlet Problem \eqref{sDproblem}, and equivalently $h$-PDP has a unique (or infinitely many) solution(s).
\item Singular Parabolic Dirichlet Problem \eqref{sDproblem}, and equivalently $h$-PDP is regular (or irregular) at $\mathcal{O}$.
\item $\Omega^c\cap \mathbb{R}^{N+1}_+$ is not parabolic minimally thin at $\mathcal{O}$ (or it is, and equivalently $\Omega$ is a deleted neighborhood of $\mathcal{O}$ in parabolic minimal fine topology)
\item The series \eqref{wtest} diverges (or converges).
\end{enumerate}
\end{theorem}

\begin{theorem}\label{main theorem2}
For arbitrary open set $\tilde{\Omega}\subset\mathbb{R}^{N+1}_-$ the following conditions are equivalent:
\begin{enumerate}
\item $\infty$ is $\tilde{h}$-regular (or $\tilde{h}$-irregular).
\item Singular Parabolic Dirichlet Problem \eqref{sDproblem-}, and equivalently $\tilde{h}$-PDP has a unique (or infinitely many) solution(s).
\item Singular Parabolic Dirichlet Problem \eqref{sDproblem-}, and equivalently $\tilde{h}$-PDP is regular (or irregular) at $\infty$.
\item $\tilde{\Omega}^c\cap \mathbb{R}^{N+1}_-$ is not $\tilde{h}$-parabolic minimally thin at $\infty$ (or it is, and equivalently $\tilde{\Omega}$ is a deleted neighborhood of $\infty$ in $\tilde{h}$-parabolic minimal fine topology)
\item The series \eqref{wtest2} diverges (or converges).
\end{enumerate}
\end{theorem}

\section{Geometric Test for the Removability of the Fundamental Singularity}\label{geometrictest}
Theorems~\ref{main theorem} and \ref{main theorem2} present a full characterization of the removability of the fundamental singularity for arbitrary open sets expressed in terms of the fine-topological thinness of the exterior set near the singularity point. In the special case when the boundary of the open set is locally represented by a graph near the singularity point, thinness criteria for the removability of the singularity can be expressed in terms of the regularity on the boundary manifold near the singularity point. In a recent paper \cite{AbdullaI}, we proved the Kolmogorov-Petrovsky-type test for a special case of domains with a boundary formed by the surface of revolution around the time axis near the singularity point. Here, we present sharp geometric criteria for the class of open sets whose boundary locally near the singularity point is a surface of the form
\[ x_1=\phi(\bar{x},t), \ \bar{x}=(x_2,...,x_n),\]
and the major problem is {\it to express the problem of the removability of the singularity in terms of the minimal regularity of $\phi$ at the singularity point.}
\subsection{Geometric Iterated Logarithm Test}\label{iteratedlogtest} For a fixed point $\gamma=(\gamma_1,\bar{\gamma})\in\mathbb{R}^{N}$ consider upper parabolic neighborhood of $(\bar{\gamma},0)$
in $\mathbb{R}^{N-1}\times \mathbb{R}_+$
\[ G=\{(\bar{x},t)\in\mathbb{R}^{N-1}\times\mathbb{R}_+: 0<\xi<\delta\}, \]
where $\xi:=\alpha t +\beta |\bar{x}-\bar{\gamma}|^2$, and $\alpha,\beta,\delta >0$. Assume that after suitable rotation of $x$-coordinates an open set $\Omega\in\mathbb{R}^{N+1}_+$ is of the form 
\begin{equation}\label{omegagtest}
 \Omega=\{z=(x_1,\bar{x},t)\in\mathbb{R}^{N+1}_+: \phi(\bar{x},t)<x_1<\delta, (\bar{x},t)\in G\},
 \end{equation}
where $\phi\in C(G\cup \{(\bar{\gamma},0)\})$ with $\phi(\bar{\gamma},0)=\gamma_1$. We prove that the $h$-regularity of the boundary point $\mathcal{O}$ crucially depends on the upper bound of the following function 
\begin{equation}\label{ratiophi}
p(\bar{x},t)\equiv \frac{\phi(\bar{\gamma},0)-\phi(\bar{x},t)}{\xi^{\frac{1}{2}}}, \ (\bar{x},t)\in G.
\end{equation}
\begin{definition}\label{pdef}
Function $p$ defined in \eqref{ratiophi} iis called a parabolic $Lip(1,\frac{1}{2}; \alpha,\beta)$ quotient of $\phi$ at $(\bar{\gamma},0)$
\end{definition}
If $p$ is bounded above, then $\phi$ satisfies one-side $Lip(1,\frac{1}{2})$ condition with respect to $(\bar{x},t)$ in $\xi^{\frac{1}{2}}$-norm at the point $(\bar{\gamma},0)$. This condition is satisfactory for the $h$-regularity of $\mathcal{O}$, and it implies the geometric "exterior hyperbolic paraboloid condition" which is analogous to the exterior cone condition for the elliptic PDEs (see Section~\ref{hpctest}). This $h$-regularity result is independent of positive constants $\alpha$ and $\beta$ defining the parabolic norm. Indeed, if $p$ is bounded above in $\xi^{\frac{1}{2}}$-norm for some fixed $\alpha>0,\beta>0$, then it is so for arbitrary fixed values of $\alpha>0,\beta>0$. 
In terms of the one-sided H\"{o}lder condition $Lip(1,\frac{1}{2})$ is an optimal space for the $h$-regularity of $\mathcal{O}$. Optimal condition for the $h$-regularity of $\mathcal{O}$ can be expressd in terms of the growth rate of the parabolic $Lip(1,\frac{1}{2}; \alpha,\beta)$ quotient of $\phi$, which still provides satisfactory thinness of the exterior set near $\mathcal{O}$ for the removability of the singularity. To identify an optimal growth rate of $Lip(1,\frac{1}{2}; \alpha,\beta)$ in the following we define a key class $\mathcal{R}$.
\begin{definition}\label{classR}
Let $\mathcal{R}$ be a class of functions $\rho\in C^2(0,1)$ with $\rho(\eta)>1$, and
\begin{gather}
 \rho(\eta)\uparrow +\infty, \ \eta\rho^{-1}(\eta)\rho'(\eta)=o(1), \ \eta^2\rho^{-1}(\eta)\rho''(\eta)=o(1) \ \quad\text{as} \ \eta \downarrow 0,\label{rhodiff}\\
\int\limits_{0^+}\frac{d\eta}{\eta \rho(\eta)}=+\infty.\label{rhoint}
\end{gather}
\end{definition}
Typical examples of functions $\rho\in\mathcal{R}$ are the following:
\begin{equation}\label{rhoexamples}
\rho(\eta)=|\log\eta|,\ \rho(\eta)=|\log\eta| \prod_{k=2}^n\log_k\eta,  n=2,3,...,
\end{equation}
where we use the following notation:
\[ \log_2\eta=\log|\log \eta|, \log_n\eta=\log\log_{n-1}\eta, \ n\geq 3. \]

Our main geometric test for the $h$-regularity of the singularity point $\mathcal{O}$ reads:
\begin{theorem}\label{geometricmaintheorem}
Let $\gamma=(\gamma_1,\bar{\gamma})\in\mathbb{R}\times\mathbb{R}^{N-1}$ is fixed, and $\Omega\in\mathbb{R}^{N+1}_+$ be an open set such that for some $\delta>0$, $\Omega\cap\{t<\delta\}$ is of the form \eqref{omegagtest}. Assume that for some $\delta>0$, a parabolic $Lip(1,\frac{1}{2}; \alpha,\beta)$ quotient of $\phi$ at $(\bar{\gamma},0)$ satisfies the upper bound
\begin{equation}\label{pbound}
p(\bar{x},t)\leq 2 \Big [ \log \rho\Big (\frac{4t^2}{\xi}\Big )\Big ]^{\frac{1}{2}}, \ \quad\text{in} \ G,
\end{equation}
where $\xi=\alpha t + \beta |\bar{x}-\bar{\gamma}|^2$, $\rho\in \mathcal{R}$, and
\[   0<\alpha+2(N-1)\beta \leq 1. \]
Then $\mathcal{O}$ is $h$-regular for $\Omega$, that is to say the fundamental singularity is removable at $\mathcal{O}$, and the singular parabolic DP has a unique solution in $\Omega$. 
\end{theorem}  

Consider the following model domains in $\mathbb{R}^{N+1}_+$ in the critical case when the conditions of Theorem~\ref{geometricmaintheorem} are satisfied with inequality replaced with equality in \eqref{pbound}:
\begin{gather}
\mathcal{G}^1_\rho=\{z: |x_1-\gamma_1|^2<4\xi\log\rho\Big (\frac{4t^2}{\xi}\Big ), 0<\xi < \delta\}\label{G1rho}\\
\mathcal{G}^2_\rho=\{z: -2\Big (\xi\log\rho\Big (\frac{4t^2}{\xi}\Big )\Big )^{\frac{1}{2}}<x_1-\gamma_1<\delta, 0<\xi < \delta\}\label{G2rho}
\end{gather}
The singularity point $\mathcal{O}$ is $h$-regular for $\mathcal{G}^1_\rho$ and $\mathcal{G}^2_\rho$ for examples of $\rho$ given in \eqref{rhoexamples}. Consider the case with $\rho(\eta)=|\log \eta|$. The following counterexample demonstrates the sharpness of the result of Theorem~\ref{geometricmaintheorem}.
\begin{corollary}\label{hregcountex}
Let 
\begin{gather}
\mathcal{G}^1_{|\log|}=\Big \{z: |x_1-\gamma_1|^2<4\xi\log_2\Big (\frac{4t^2}{\xi}\Big ), 0<\xi < \delta\Big \}\label{rholog-11}\\
\mathcal{G}^2_{|\log|}=\Big \{z: -2\Big (\xi\log_2\Big (\frac{4t^2}{\xi}\Big )\Big )^{\frac{1}{2}}<x_1-\gamma_1<\delta, 0<\xi < \delta\Big\}. \label{rholog-1}
\end{gather}
Then $\mathcal{O}$ is $h$-regular (or $h$-irregular), that is to say the singular DP has a unique (or infinitely many) solutions in $\mathcal{G}^1_{|\log|}$ and $\mathcal{G}^2_{|\log|}$ if $\alpha+2(N-1)\beta\leq 1$ (or $\alpha>1,\beta>0$).
\end{corollary} 
Clearly, the same result is true if we take any other examples of $\rho$ from \eqref{rhoexamples}. 

Another conclusion of Theorem~\ref{geometricmaintheorem} is that the criteria for the removability of a fundamental singularity are different for the heat equation
\begin{equation}\label{heat1}
u_t=\Delta u,
\end{equation}
and 
\begin{equation}\label{heata}
u_t=a\Delta u, \ 0<a<1.
\end{equation}
\begin{corollary}\label{hregcountex2}
Let $0<a<1$ and $0<\epsilon \leq 1-a$ are given numbers. Consider the domains $\mathcal{G}^1_{|\log|}$ and $\mathcal{G}^2_{|\log|}$ with $\alpha=1-\frac{\epsilon}{2}$, $\beta=\frac{\epsilon}{4(N-1)}$. Then the boundary singularity point $\mathcal{O}$ is $h$-regular for the domains $\mathcal{G}^1_{|\log|}$ and $\mathcal{G}^2_{|\log|}$ regarding the heat equation \eqref{heat1}, and the singular PDP has a unique solution. On the other hand, $\mathcal{O}$ is $h$-irregular for the domains $\mathcal{G}^1_{|\log|}$ and $\mathcal{G}^2_{|\log|}$ regarding the heat equation \eqref{heata}, and the corresponding singular PDP has infinitely many solutions.
\end{corollary}
Theorem~\ref{geometricmaintheorem} provides a general geometric condition for the $h$-regularity of boundary singularity point $\mathcal{O}$ for an arbitrary oprn set $\Omega\subset\mathbb{R}^{N+1}_+$. Consider a set
\[
\mathcal{A}_\rho=\Big \{z: -\delta<x_1<-2\Big ((\alpha t+ \beta|\bar{x}|^2)\log\rho\Big (\frac{4t^2}{\alpha t+ \beta|\bar{x}|^2}\Big )\Big )^{\frac{1}{2}}, 0<\alpha t+ \beta|\bar{x}|^2 < \delta\Big\}
\]
and assume that $\alpha,\beta,\rho$ satisfy conditions of Theorem~\ref{geometricmaintheorem}. Consider the rigid body displacements of $\mathcal{A}_\rho$ composed of $\gamma$-translation in $x$-space, and rotation around the $(\gamma,t)$-axis.
\begin{definition}
We shall say that $\Omega\subset\mathbb{R}^{N+1}_+$ satisfies the exterior $\mathcal{A}_\rho$-condition at the boundary point $\mathcal{O}=(\gamma,0)$, if for some $\delta>0$ after the above mentioned displacement $\mathcal{A}_\rho$ lies in the exterior of $\Omega\cap\{z: t\leq \delta\}$. 
\end{definition}
Theorem~\ref{geometricmaintheorem} implies the following result:
\begin{corollary}
$\mathcal{O}$ is $h$-regular and the singular parabolic DP has a unique solution in $\Omega$, if $\Omega$ satisfies the exterior $\mathcal{A}_\rho$-condition at $\mathcal{O}$.
\end{corollary}

Next, we formulate a similar geometric criterion for the $\tilde{h}$-regularity of $\infty$, and for the unique solvability of the singular parabolic Dirichlet problem in an arbitrary open set $\tilde{\Omega}\subset\mathbb{R}^{N+1}_-$. Let 
\[ \tilde{G}=\{(\bar{y},t)\in\mathbb{R}^{N-1}\times\mathbb{R}_-: -\infty<\tau<\delta\}, \]
where $\delta<0$. Assume that after suitable rotation of $y$-coordinates an open set $\tilde{\Omega}\in\mathbb{R}^{N+1}_-$ is of the form 
\begin{equation}\label{tildeomegagtest}
 \tilde{\Omega}=\{z=(y_1,\bar{y},\tau)\in\mathbb{R}^{N+1}_-: \tilde{\phi}(\bar{y},\tau)<y_1<+\infty, (\bar{y},\tau)\in\tilde{G}\},
 \end{equation}
where $\tilde{\phi}\in C(\tilde{G})$. We prove that the $\tilde{h}$-regularity of $\infty$ for $\tilde{\Omega}$ crucially depends on the upper bound of the following function 
\begin{equation}\label{ratiophitilde}
\tilde{p}(\bar{y},\tau)\equiv \frac{-2\tau\gamma_1-\tilde{\phi}(\bar{y},\tau)}{[-\alpha\tau+\beta|\bar{y}+2\tau\bar{\gamma}|^2]^{\frac{1}{2}}}, \ (\bar{y},\tau)\in \tilde{G},
\end{equation}
where $\alpha,\beta>0$ are fixed numbers. Upper bound on $\tilde{p}$ restricts the decay rate of $\tilde{\phi}$ to $-\infty$. 
If $\gamma=0$ ($\tilde{h}\equiv 1$), per Definition~\ref{pdef}, it is natural to call a function $\tilde{p}$, a parabolic $Lip(1,\frac{1}{2}; \alpha,\beta)$ quotient of $\tilde{\phi}$ at $\infty$. 
In general, if $\gamma\neq 0$, a function $\tilde{p}$ is derived by applying the following {\it shear transformation} of $\mathbb{R}^{N+1}_-$ to the parabolic $Lip(1,\frac{1}{2}; \alpha,\beta)$ quotient of $\tilde{\phi}$ at $\infty$:
\begin{equation}\label{shear}
\mathcal{S}_\gamma: (y_1,\bar{y},\tau)\in\mathbb{R}^{N+1}_- \mapsto \mathcal{S}(y_1,\bar{y},\tau)=(y_1+2\gamma_1\tau,\bar{y}+2\bar{\gamma}\tau,\tau)\in \mathbb{R}^{N+1}_-
\end{equation}
$\mathcal{S}_\gamma$ performs horizontal displacement parallel to $x$-space in the direction of the vector $-2\gamma$ with the factor $\tau$, and $\{\tau=0\}$ is a fixed hyperplane. Precisely, the transformation $\mathcal{S}_\gamma$ performs a {\it laminar flow} in $\mathbb{R}^{N+1}_-$ in the direction of the vector $-2\gamma$ between parallel planes in a relative motion with the speed propotional to the signed distance from the fixed hyperplane $\{\tau=0\}$.
\begin{definition}\label{tildepdef}
Function $\tilde{p}$ defined in \eqref{ratiophitilde} is called an {\it $\mathcal{S}_\gamma$-sheared parabolic $Lip(1,\frac{1}{2}; \alpha,\beta)$ quotient of $\tilde{\phi}$ at $\infty$.}
\end{definition}

If $\tilde{p}$ is bounded above, then $\tilde{\phi}$ satisfies one-sided decay rate restriction as follows:
\begin{equation}\label{exthyppar}
\tilde{\phi}(\bar{y},\tau)+2\tau\gamma_1\geq -C[-\alpha\tau+\beta|\bar{y}+2\tau\bar{\gamma}|^2]^{\frac{1}{2}}, \ \quad\text{in} \ \tilde{G}.
\end{equation}
This condition is satisfactory for the $\tilde{h}$-regularity of $\infty$. In the particular case, when $\gamma=0$ ($\tilde{h}\equiv 1$), the condition \eqref{exthyppar} is equivalent to the geometric "exterior hyperbolic paraboloid condition" for the uniqueness of the bounded solution of the parabolic DP in arbitrary open sets revealed in \cite{Abdulla3} (see Section~\ref{hpctest}).  Hence, the condition \eqref{exthyppar} provides a geometric {\it exterior $\mathcal{S}_\gamma$-sheared hyperbolic paraboloid} condition for the $\tilde{h}$-regularity of $\infty$ (see Section~\ref{hpctest}). 

In terms of power-like functions, a lower bound in \eqref{exthyppar} is an optimal bound for the $\tilde{h}$-regularity of $\infty$. Optimal condition for the $\tilde{h}$-regularity of $\infty$ can be expressed in terms of the growth rate of the $\mathcal{S}_\gamma$-sheared quotient $\tilde{p}$ to $+\infty$, which still provides satisfactory thinness of the exterior set near $\infty$ for the removability of the singularity. To identify an optimal growth rate of $\mathcal{S}_\gamma$-sheared quotient $\tilde{p}$ to $+\infty$, we define a key class $\mathcal{R}^{'}$.
\begin{definition}\label{classR1}
Let $\mathcal{R}^{'}$ be a class of functions $\rho\in C^2(-\infty,-1), \rho(\eta)>1$, satisfying \eqref{rhodiff} as $\eta\downarrow -\infty$; and
\begin{equation}\label{tilderhoint}
\int\limits_{-\infty}\frac{d\eta}{\eta \rho(\eta)}=-\infty.
\end{equation}
\end{definition}
Typical examples of functions $\rho\in\mathcal{R}^{'}$ are the following:
\begin{equation}\label{tilderhoexamples}
\rho(\eta)=\log|\eta|,\ \rho(\eta)=\log|\eta| \prod_{k=2}^n\log_k|\eta|, \   n=2,3,...,
\end{equation}

Our main geometric test for the $\tilde{h}$-regularity of the singularity point $\infty$ reads:
\begin{theorem}\label{tildegeometricmaintheorem}
Let $\gamma=(\gamma_1,\bar{\gamma})\in\mathbb{R}\times\mathbb{R}^{N-1}$ is fixed, and $\tilde{\Omega}\in\mathbb{R}^{N+1}_-$ be an open set such that for some $\delta<0$, $\tilde{\Omega}\cap\{t<\delta\}$ is of the form \eqref{tildeomegagtest}. Assume that for some $\delta<0$, $\mathcal{S}_\gamma$-sheared quotient of $\tilde{\phi}$ at $\infty$ satisfies the upper bound
\begin{equation}\label{tildepbound}
\tilde{p}(\bar{y},\tau)\leq 2 [ \log \rho(\tilde{\xi})]^{\frac{1}{2}}, \ \quad\text{in} \ \tilde{G}.
\end{equation}
where $\tilde{\xi}=\alpha\tau-\beta |\bar{y}+2\bar{\gamma}\tau|^2$, $\rho\in\mathcal{R}^{'}$ and
\[0<\alpha+2(N-1)\beta\leq 1. \]
Then $\infty$ is $\tilde{h}$-regular for $\tilde{\Omega}$, that is to say the fundamental singularity at $\infty$ is removable, and the singular parabolic DP has a unique solution in $\tilde{\Omega}$. 
\end{theorem}  
%%%%%%%%
%%%%%%%%
%%%%%%%

Consider the domains in $\mathbb{R}^{N+1}_-$ in the critical case when the conditions of Theorem~\ref{tildegeometricmaintheorem} are satisfied with inequality replaced with equality in \eqref{tildepbound}:
\begin{gather}
\mathcal{L}^1_\rho=\{z: |y_1+2\tau\gamma_1|^2<-4\tilde{\xi}\log\rho(\tilde{\xi}), (\bar{y},\tau)\in\tilde{G}\}\label{L1rho}\\
\mathcal{L}^2_\rho=\{z\in\mathbb{R}^{N+1}_-: -2\tau\gamma_1-2\Big (-\tilde{\xi}\log\rho(\tilde{\xi})\Big )^{\frac{1}{2}}<y_1<+\infty, (\bar{y},\tau)\in\tilde{G}\}\label{L2rho}
\end{gather}
The singularity point at $\infty$ is $\tilde{h}$-regular for $\mathcal{L}^1_\rho$ and $\mathcal{L}^2_\rho$ with $\rho$ being selected as in \eqref{rhoexamples}. Consider $\mathcal{L}^1_\rho$ and $\mathcal{L}^2_\rho$ with $\rho(\eta)=\log|\eta|$. The following counterexample demonstrates the sharpness of the result of Theorem~\ref{tildegeometricmaintheorem}.
\begin{corollary}\label{tildehregcountex}
Let 
\begin{gather}
\mathcal{L}^1_{\log|\cdot|}=\{z: |y_1+2\tau\gamma_1|^2<-4\tilde{\xi}\log_2|\tilde{\xi}|,  (\bar{y},\tau)\in\tilde{G}\}\label{tilderholog-11}\\
\mathcal{L}^2_{\log|\cdot|}=\{z: -2\tau\gamma_1-2\Big (-\tilde{\xi}\log_2|\tilde{\xi}|\Big )^{\frac{1}{2}}<y_1<+\infty,  (\bar{y},\tau)\in\tilde{G}\}. \label{tilderholog-1}
\end{gather}
Then $\infty$ is $\tilde{h}$-regular (or $\tilde{h}$-irregular), that is to say the singular PDP has a unique (or infinitely many) solutions in $\mathcal{L}^1_{\log|\cdot|}$ and $\mathcal{L}^2_{\log|\cdot|}$ if $\alpha+2(N-1)\beta\leq 1$ (or $\alpha>1,\beta>0$).
\end{corollary} 
The same result is true if we take any other examples of $\rho$ from \eqref{tilderhoexamples}. 

The next corollary demonstrates that the criteria for the $\tilde{h}$-regularity of $\infty$ are different for the heat equations \eqref{heat1} and \eqref{heata}.
\begin{corollary}\label{tildehregcountex2}
Let $0<a<1$ and $0<\epsilon \leq 1-a$ are given numbers. Consider the domains $\mathcal{L}^1_{\log|\cdot|}$ and $\mathcal{L}^2_{\log|\cdot|}$ with $\alpha=1-\frac{\epsilon}{2}$, $\beta=\frac{\epsilon}{4(N-1)}$. Then $\infty$ is $\tilde{h}$-regular for the domains $\mathcal{L}^1_{\log|\cdot|}$ and $\mathcal{L}^2_{\log|\cdot|}$ regarding the heat equation \eqref{heat1}, and the singular PDP has a unique solution. On the other hand, $\infty$ is $\tilde{h}$-irregular for the domains $\mathcal{L}^1_{\log|\cdot|}$ and $\mathcal{L}^2_{\log|\cdot|}$ regarding the heat equation \eqref{heata}, and the corresponding singular PDP has infinitely many solutions.
\end{corollary}
Theorem~\ref{tildegeometricmaintheorem} provides a general geometric condition for the $\tilde{h}$-regularity of $\infty$ for an arbitrary open set $\tilde{\Omega}\subset\mathbb{R}^{N+1}_-$. Consider a set
\[
\tilde{\mathcal{A}}_\rho=\Big \{z: y_1<-2\Big ((-\alpha \tau+ \beta|\bar{y}|^2)\log\rho(\alpha\tau-\beta|\bar{y}|^2 )\Big )^{\frac{1}{2}}, \bar{y}\in \mathbb{R}^{N-1}, \tau<\delta\Big\}
\]
and assume that $\alpha,\beta,\rho$ satisfy conditions of Theorem~\ref{tildegeometricmaintheorem}. 

\begin{definition}
Let $\gamma\in\mathbb{R}^{N}$ is fixed. We shall say that an open set $\Omega\subset\mathbb{R}^{N+1}_-$ satisfies the {\it $\mathcal{S}_\gamma$-sheared exterior $\tilde{\mathcal{A}}_\rho$-condition at $\infty$}, if for some $\delta<0$ the set $\mathcal{S}_\gamma\tilde{\mathcal{A}}_\rho$ lies in the exterior of $\Omega\cap\{z: t\leq \delta\}$. 
\end{definition}
Theorem~\ref{tildegeometricmaintheorem} implies the following result:
\begin{corollary}\label{extArho}
$\infty$ is $\tilde{h}$-regular and the singular parabolic DP has a unique solution in $\Omega\subset\mathbb{R}^{N+1}_-$, if $\Omega$ satisfies the {\it $\mathcal{S}_\gamma$-sheared exterior $\tilde{\mathcal{A}}_\rho$-condition at $\infty$}
\end{corollary}
Note that if $\gamma=0$ ($\tilde{h}\equiv 1$, then $\mathcal{S}_\gamma$ is an identity map, and the result of Corollary~\ref{extArho} coincides with the exterior $\mathcal{A}_\rho$ condition for the regularity of $\infty$ and uniqueness of bounded solutions of the parabolic DP revealed in \cite{Abdulla2} (Corollary 1.1, page 9).  
%%%%%%%%%%
%%%%%%%%%
%%%%%%%%%
\subsection{Exterior Hyperbolc Paraboloid Condition}\label{hpctest}
In this section, we describe an exterior hyperbolic paraboloid condition for the removability of the fundamental singularity for the heat equation. In terms of the one-sided H\"{o}lder condition on the boundary manifold for the $h$-regularity of the finite boundary singularity point, $Lip(1,\frac{1}{2})$ is an optimal space. The subsurface of the hyperbolic paraboloid is a typical example of this space. Similarly, in terms of the one-sided power-like decay rate of the boundary manifold for the $\tilde{h}$-regularity of $\infty$, $\mathcal{S}_\gamma$-shearing of the hyperbolic paraboloid provides an optimal rate. Exterior hyperbolic paraboloid condition is a natural analogy of the exterior cone condition for elliptic PDEs. The parabolic scaling dictates its relevance. It was first introduced in \cite{AbdullaSiam} as a geometric condition for the boundary regularity. In \cite{Abdulla2} it is proved that the exterior hyperbolic paraboloid condition at infinity is sufficient for the regularity of $\infty$, and for the existence of the unique bounded solution of the parabolic DP in an arbitrary open set.

Let $\gamma=(\gamma_1,\bar{\gamma})\in\mathbb{R}^{N}$ be fixed. Consider model domains
\begin{gather}
\mathcal{D}^1=\Big \{z: |x_1-\gamma_1|^2<M^2(t+|\bar{x}-\bar{\gamma}|^2), 0<t+|\bar{x}-\bar{\gamma}|^2 < \delta\Big \}\label{EHP1}\\
\mathcal{D}^2=\Big \{z: -M(t+|\bar{x}-\bar{\gamma}|^2)^{\frac{1}{2}}<x_1-\gamma_1<M\delta^{\frac{1}{2}}, 0<t+|\bar{x}-\bar{\gamma}|^2 < \delta\Big\},\label{EHP2}
\end{gather}
where $M>0, \delta>0$. The domains are formed by the intersection of the hyperbolic paraboloid (Figure 3)
\begin{figure}[ht]\label{hyperpar}
    \centering
    \includegraphics[width=0.9\textwidth, trim=0.5cm 0.5cm 0.5cm 0.5cm, clip]{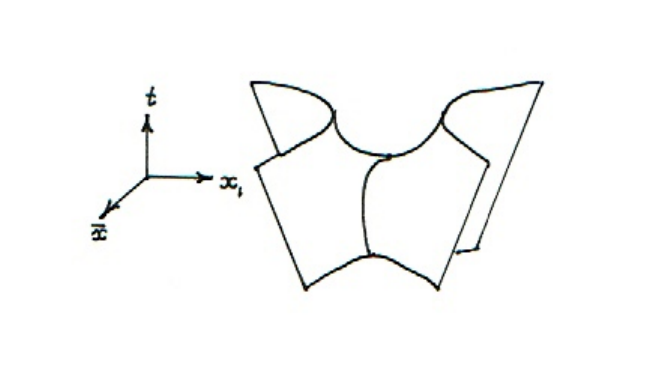}
    \caption{Hyperbolic paraboloid with saddle point at $\gamma$}
\end{figure}
\[ (x_1-\gamma_1)^2=M^2(t+|\bar{x}-\bar{\gamma}|^2) \]
with the parabolic cylinder 
\[ \mathcal{C}=\{z: 0<t+|\bar{x}-\bar{\gamma}|^2 < \delta, |x_1-\gamma_1|<M\delta^{\frac{1}{2}}\} \]
The domain $\mathcal{D}_1$ is a part of the cylinder $\mathcal{C}$ between two symmetric subsurfaces $x_1=\phi_{\mp}(\bar{x},t)\equiv \gamma_1\mp M(t+|\bar{x}-\bar{\gamma}|^2)^{\frac{1}{2}}$ (Figure 4). 
\begin{figure}[ht]\label{D1}
    \centering
    \includegraphics[width=0.9\textwidth, trim=0.3cm 0.3cm 0.3cm 0.3cm, clip]{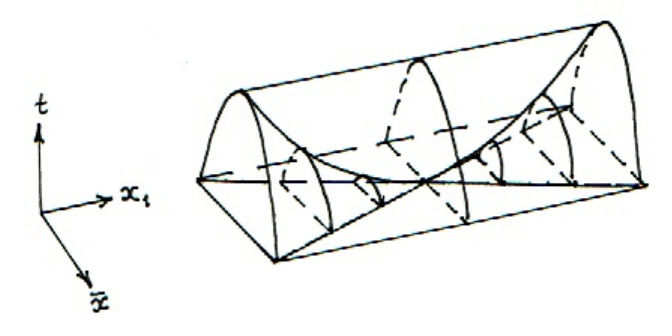}
    \caption{The domain $\mathcal{D}^1$}
\end{figure}

The domain $\mathcal{D}_2$ is a part of the cylinder $\mathcal{C}$ lying on one side of the surface $x_1=\phi_-(\bar{x},t)$ (Figure 5).
\begin{figure}[ht]\label{D2}
    \centering
    \includegraphics[width=0.9\textwidth, trim=0.3cm 0.3cm 0.3cm 0.3cm, clip]{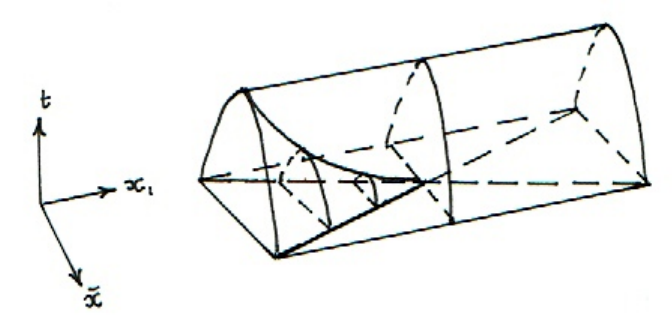}
    \caption{The domain $\mathcal{D}^2$}
\end{figure}

 It is easy to see that the parabolic $Lip(1,\frac{1}{2}; 1,1)$ quotient of $\phi_-$ at $(\bar{\gamma},0)$ is identically equal to $M$. Accordingly, for any fixed $\alpha,\beta>0$, $Lip(1,\frac{1}{2}; \alpha,\beta)$ quotient of $\phi_-$ at $(\bar{\gamma},0)$ is bounded by $M/(min(\alpha,\beta))^{\frac{1}{2}}$. 

\begin{theorem}\label{OregD1D2}
Boundary singularity point $\mathcal{O}=(\gamma,0)$ is $h$-regular for $\mathcal{D}_1$ (or $\mathcal{D}_2$), that is to say the fundamental singularity at $\mathcal{O}$ is removable, and the singular PDP has a unique solution in $\mathcal{D}_1$ (or $\mathcal{D}_2$).
\end{theorem}

Theorem~\ref{OregD1D2} provides a simple geometric sufficient condition for the $h$-regularity of the boundary singularity point $\mathcal{O}$ for arbitrary open set $\Omega\subset\mathbb{R}^{N+1}_+$. Let $\mathcal{P}$ be a region bounded by the subsurface of the hyperbolic paraboloid and hyperplanes $\{x_1=-M\delta^{\frac{1}{2}}\}$ and $\{t=0\}$  (Figure 6)
\begin{equation}\label{exteriorhpcondition}
\mathcal{P}= \{z: -M\delta^{\frac{1}{2}}<x_1< -M(t+|\bar{x}|^2)^{\frac{1}{2}}, 0< t+ |\bar{x}|^2 < \delta\}
\end{equation}
\begin{figure}[ht]\label{tilde}
    \centering
    \includegraphics[width=0.9\textwidth, trim=0.1cm 0.1cm 0.1cm 0.1cm, clip]{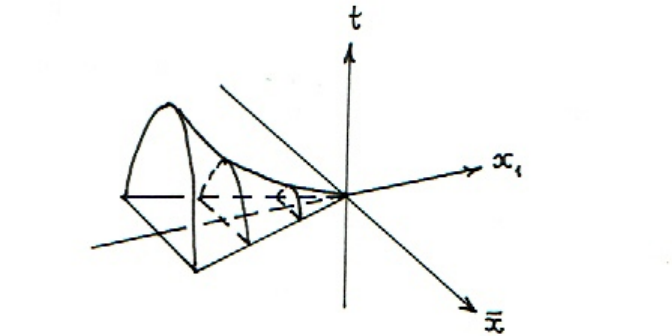}
    \caption{The region $\mathcal{P}$}
\end{figure}
Consider the rigid body displacements of $\mathcal{P}$ composed of $\gamma$-translation in $x$-space, and rotation around the $(\gamma,t)$-axis.
\begin{definition}
We shall say that $\Omega\subset\mathbb{R}^{N+1}_+$ satisfies the exterior hyperbolic paraboloid condition at the boundary point $\mathcal{O}=(\gamma,0)$, if for some $M>0$ and $\delta>0$ after the above mentioned displacement $\mathcal{P}$ lies in the exterior of $\Omega\cap\{z: t\leq \delta\}$. 
\end{definition}
Theorem~\ref{OregD1D2} implies the following result:
\begin{corollary}\label{extHPcondition}
$\mathcal{O}$ is $h$-regular and the singular parabolic DP has a unique solution in $\Omega$, if $\Omega$ satisfies the exterior hyperbolic paraboloid condition at $\mathcal{O}$.
\end{corollary}
Next, we formulate a similar geometric criterion for the $\tilde{h}$-regularity of $\infty$, and for the unique solvability of the singular parabolic Dirichlet problem in an arbitrary open set $\tilde{\Omega}\subset\mathbb{R}^{N+1}_-$. Consider model domains
\begin{gather}
\mathcal{E}_1= \{z: |y_1+2\tau\gamma_1|^2<M^2(-\tau+|\bar{y}+2\tau\bar{\gamma}|^2), \bar{y}\in \mathbb{R}^{N-1}, \tau<\delta  \}\label{tildeEHP1}\\
\mathcal{E}_2= \{z: -2\tau\gamma_1-M(-\tau+|\bar{y}+2\tau\bar{\gamma}|^2)^{\frac{1}{2}}<y_1<+\infty, \bar{y}\in \mathbb{R}^{N-1}, \tau<\delta\},\label{tildeEHP2}
\end{gather}
where $M>0, \delta<0$. The domains are formed by the $\mathcal{S}_\gamma$-sheared hyperbolic paraboloid
\begin{equation}\label{SgammaHP}
 (y_1+2\tau\gamma_1)^2=M^2(-\tau+|\bar{y}+2\tau\bar{\gamma}|^2)
 \end{equation}
If $\gamma=0$ ($\mathcal{S}_\gamma$ is an identity map), \eqref{SgammaHP} is an hyperbolic paraboloid
\begin{equation}\label{S0HP}
y_1^2=M^2(-\tau+|\bar{y}|^2).
\end{equation}
In general, \eqref{SgammaHP} is an image of the hyperbolic paraboloid \eqref{S0HP} under the shear transformation $\mathcal{S}_\gamma$. 
The domain $\mathcal{E}_1$ is a part of the half-space $\{t<\delta\}$ between two subsurfaces of \eqref{SgammaHP}:  $y_1=\tilde{\phi}_{\mp}(\bar{y},\tau)\equiv -2\gamma_1\tau\mp M(-\tau+|\bar{y}+2\bar{\gamma}\tau|^2)^{\frac{1}{2}}$ (Fig.4). The domain $\mathcal{E}_2$ is a part of the half-space $\{t<\delta\}$ lying on one side of the surface $y_1=\tilde{\phi}_-(\bar{y},\tau)$. It is easy to see that the $\mathcal{S}_\gamma$-sheared parabolic $Lip(1,\frac{1}{2}; 1,1)$ quotient of $\tilde{\phi}_-$ at $\infty$ is identically equal to $M$. Accordingly, for any fixed $\alpha,\beta>0$, $Lip(1,\frac{1}{2}; \alpha,\beta)$ quotient of $\phi_-$ at $(\bar{\gamma},0)$ is bounded by $M/min(\alpha^{\frac{1}{2}},\beta^{\frac{1}{2}})$. 

\begin{theorem}\label{tildeOregD1D2}
Boundary point $\infty$ is $\tilde{h}$-regular for $\mathcal{E}_1$ (or $\mathcal{E}_2$), that is to say the fundamental singularity at $\infty$ is removable, and the singular PDP has a unique solution in $\mathcal{E}_1$ (or $\mathcal{E}_2$).
\end{theorem}

Theorem~\ref{tildeOregD1D2} provides a simple geometric sufficient condition for the $\tilde{h}$-regularity of the boundary point $\infty$ for arbitrary open set $\Omega\subset\mathbb{R}^{N+1}_-$. Consider a region (Figure 7)
\begin{equation}\label{tildeexteriorhpcondition}
\mathcal{\tilde{P}}= \{z: -\infty<y_1< -M(-\tau+|\bar{y}|^2)^{\frac{1}{2}}, \bar{y}\in\mathbb{R}^{N-1}, \tau<\delta  \}
\end{equation}
\begin{figure}[ht]\label{tildeP}
    \centering
    \includegraphics[width=0.8\textwidth, trim=0.3cm 0.3cm 0.3cm 0.3cm, clip]{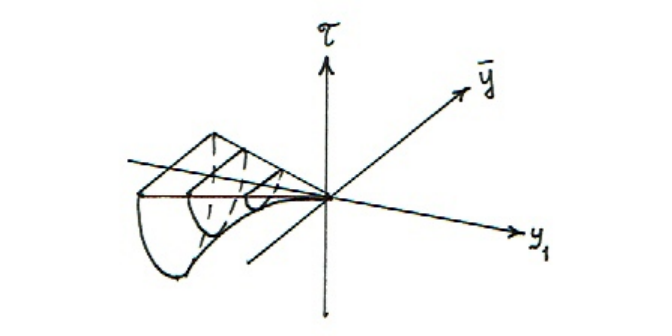}
    \caption{The region $\tilde{\mathcal{P}}$}
\end{figure}

\begin{definition}
Let $\gamma\in\mathbb{R}^{N}$ is fixed. We shall say that an open set $\Omega\subset\mathbb{R}^{N+1}_-$ satisfies the {\it $\mathcal{S}_\gamma$-sheared exterior hyperbolic paraboloid condition at $\infty$}, if for some $\delta<0$ the set $\mathcal{S}_\gamma \mathcal{\tilde{P}}$ lies in the exterior of $\Omega\cap\{z: t\leq \delta\}$. 
\end{definition}
Theorem~\ref{tildegeometricmaintheorem} implies the following result:
\begin{corollary}\label{extSgammahyperpar}
$\infty$ is $\tilde{h}$-regular and the singular parabolic DP has a unique solution in $\Omega\subset\mathbb{R}^{N+1}_-$, if $\Omega$ satisfies the {\it $\mathcal{S}_\gamma$-sheared exterior hyperbolic paraboloid condition at $\infty$}
\end{corollary}
Note that if $\gamma=0$ ($\tilde{h}\equiv 1$, then $\mathcal{S}_\gamma$ is an identity map, and the result of Corollary~\ref{extSgammahyperpar} coincides with the exterior hyperbolic paraboloid condition for the regularity of $\infty$ and uniqueness of bounded solutions of the parabolic DP introduced in \cite{Abdulla2} (Corollary 1.2, page 10).

\section{$h$-Parabolic Potential Theory}\label{E:1:7}
%%%%%%%%%%%
The equivalence of two problems formulated in $\mathbb{R}^{N+1}_+$ and $\mathbb{R}^{N+1}_-$ is a consequence of the
{\bf Appell transformation}.
Consider a homeomorphism $A: \mathbb{R}^{N+1}_+\cup\{\mathcal{O}\} \mapsto \mathbb{R}^{N+1}_-\cup\{\infty\}$ with
\begin{equation}\label{Appellmap}
\left\{
\begin{array}{l}
(x,t)\in\mathbb{R}^{N+1}_+ \mapsto A(x,t)=\Big (\frac{x}{2t}, -\frac{1}{4t}\Big )\in \mathbb{R}^{N+1}_-; \ \ A(\mathcal{O})=\infty\\
(x,t)\in\mathbb{R}^{N+1}_- \mapsto A^{-1}(x,t)=\Big (-\frac{x}{2t}, -\frac{1}{4t}\Big )\in \mathbb{R}^{N+1}_+; \ \ A^{-1}(\infty)=\mathcal{O}
\end{array}\right.
\end{equation}
Let $\mathcal{P}(\Omega)$ be a class of parabolic functions
in an open set $\Omega$. Given open set $\Omega\subset\mathbb{R}^{N+1}_+$, the Appell transformation is a homeomorphism ${\bf A}: \mathcal{P}(\Omega) \mapsto \mathcal{P}(A\Omega)$ defined as
\begin{equation}\label{Atransform}
\left\{
\begin{array}{l}
 \Omega\subset\mathbb{R}^{N+1}_+: u\in\mathcal{P}(\Omega) \mapsto {\bf A}u(z)=(-\frac{\pi}{t})^{\frac{N}{2}} e^{-\frac{|x|^2}{4t}}u(A^{-1}(z))\in \mathcal{P}(A\Omega) \\ \\
 \Omega\subset\mathbb{R}^{N+1}_-: v\in\mathcal{P}(\Omega) \mapsto {\bf A}^{-1}v(z)=F(z) v(A(z))\in \mathcal{P}(A^{-1}\Omega)
  \end{array}\right.
\end{equation}
The claim follows from the following formula:
\begin{equation}\label{appelltransform}
\left\{
\begin{array}{l}
\mathcal{H}[{\bf A}u(z)]=\frac{\pi^{N/2}}{4} (-t)^{-\frac{N}{2}-2}e^{-\frac{|x|^2}{4t}} \mathcal{H}[u(A^{-1}(z)], \ z\in A\Omega\subset \mathbb{R}^{N+1}_-\\ \\
\mathcal{H}[{\bf A}^{-1}v(z)]=\frac{1}{4t^2} F(z) \mathcal{H}[v(A(z)], \ z\in A^{-1}\Omega\subset \mathbb{R}^{N+1}_+.
\end{array}\right.
\end{equation}
In particular, the Appell transform of $h$ is given by
\begin{equation}\label{afs1}
\tilde{h}(x,t)={\bf A}h(x,t)= e^{\langle x,\gamma\rangle + |\gamma|^2 t},
\end{equation}
as it is defined in \eqref{afs}. In fact, the Appell transformation generates an homomorphism ${\bf A}: \mathcal{S}(\Omega) \mapsto \mathcal{S}(A\Omega)$, where $\mathcal{S}(\Omega)$ and  $\mathcal{S}(A\Omega)$ denote the class of superparabolic functions in $\Omega$ and $A\Omega$ respectively. The Appell transformation
establishes equivalence between the singular PDPs \eqref{sDproblem}, \eqref{sDproblem-} and the corresponding Problems 1-3 (see Lemmas IV.1, IV.2 of \cite{AbdullaI}). 

In particular, the function $u$ is $h$-parabolic (or $h$-super- or $h$-subparabolic) in $\Omega\subset\mathbb{R}^{N+1}_+$ if and only if $u(A^{-1}(w))$ is $\tilde{h}$-parabolic (or $\tilde{h}$-super- or $\tilde{h}$-subparabolic) in $A\Omega\subset \mathbb{R}^{N+1}_-$ ({Lemma IV.1, \cite{AbdullaI}). The claim follows from the following formulae : if $u:\mathbb{R}^{N+1}_+\to \mathbb{R}$ is smooth, then
\begin{equation}\label{HuA-1}
\frac{\mathcal{H}[h(z)u(z)]}{h(z)}\Big |_{z=A^{-1}(w)}=4\tau^2\frac{\mathcal{H}[\tilde{h}(w)u(A^{-1}(w))]}{\tilde{h}(w)}.
\end{equation}
Equivalently,  if $u:\mathbb{R}^{N+1}_-\to \mathbb{R}$ is smooth, then
\begin{equation}\label{HuA}
\frac{\mathcal{H}[\tilde{h}(w)u(w)]}{\tilde{h}(w)}\Big |_{w=A(z)}=4t^2\frac{\mathcal{H}[h(z)u(A(z))]}{h(z)}.
\end{equation}
{\it Proof of \eqref{HuA-1} and \eqref{HuA}}: Given a smooth function  $u:\mathbb{R}^{N+1}_+\to \mathbb{R}$, by \eqref{Atransform} we have
\[ {\bf A}(uh)(w)=\Big(-\frac{\pi}{\tau}\Big )^{\frac{N}{2}}e^{-\frac{|y|^2}{4\tau}}u(A^{-1}(w))h(A^{-1}(w))=u(A^{-1}(w))\tilde{h}(w). \]
From \eqref{appelltransform} it follows that
\[ \mathcal{H}[{\bf A}(uh)](y,\tau)=\frac{\pi^{N/2}}{4}(-\tau)^{-\frac{N}{2}-2}e^{-\frac{|y|^2}{4\tau}} \mathcal{H}[u(x,t)h(x,t)]|_{(x,t)=A^{-1}(y,\tau)}, \]
and hence
\[ \mathcal{H}[u(x,t)h(x,t)]|_{(x,t)=A^{-1}(y,\tau)}=\frac{\mathcal{H}[\tilde{h}(y,\tau)u(A^{-1}(y,\tau))]}{\frac{\pi^{N/2}}{4}(-\tau)^{-\frac{N}{2}-2}e^{-\frac{|y|^2}{4\tau}}}, \]
which implies \eqref{HuA-1}. Simirlarly, given a smooth function  $u:\mathbb{R}^{N+1}_-\to \mathbb{R}$, by \eqref{Atransform} we have
\[ {\bf A}^{-1}(\tilde{h}u)(z)=F(z)\tilde{h}(A(z))u(A(z))=h(z)u(A(z)). \]
From \eqref{appelltransform} it follows that
\[ \mathcal{H}[{\bf A}^{-1}(\tilde{h}u)](z)=\frac{1}{4t^2}F(z)\mathcal{H}[\tilde{h}(w)u(w)]|_{w=A(z)}, \]
and hence
\[ \mathcal{H}[\tilde{h}(w)u(w)]|_{w=A(z)}=4t^2\frac{\mathcal{H}[h(z)u(A(z))]}{F(z)}, \]
which implies \eqref{HuA}.   \hfill{$\square$}

Next lemma demonstrates that the Appell transformation
generates a one-to-one mapping between $h$- and $\tilde{h}$-reductions of $1$, and establishes an equivalence between Problems 4 in the context of singular PDPs \eqref{sDproblem}, \eqref{sDproblem-}. 
\begin{lemma}\label{reductionof1}
Given an open set $\Omega\in\mathbb{R}^{N+1}_+$ and $\delta>0$, $\ ^h\hat{R}^{\Omega^c\cap \{t\leq \delta\}}_1$ is a smoothed reduction of $1$ in $\Omega^c\cap\{t\leq \delta\}$ if and only if $\ ^h\hat{R}^{\Omega^c\cap \{t\leq \delta\}}_1(A^{-1}w)$ is a smoothed reduction of $1$ in $(A\Omega)^c \cap \{\tau\leq -1/4\delta\}$, i.e.
\begin{equation}\label{reductionof1mapping}
\left\{
\begin{array}{l}
 \Omega\subset\mathbb{R}^{N+1}_+, \delta>0 : \ ^h\hat{R}^{\Omega^c\cap \{t\leq \delta\}}_1(A^{-1}w)=\ ^{\tilde{h}}\hat{R}^{(A\Omega)^c\cap \{\tau\leq -1/4\delta\}}_1(w), \ w\in\mathbb{R}^{N+1}_- \\
\tilde{\Omega}\subset\mathbb{R}^{N+1}_-, \delta<0 : \ ^{\tilde{h}}\hat{R}^{\tilde{\Omega}^c\cap \{\tau\leq \delta\}}_1(Az)=\ ^h\hat{R}^{(A^{-1}\Omega)^c\cap \{t\leq -1/4\delta\}}_1(z), \ z\in\mathbb{R}^{N+1}_+ \\
  \end{array}\right.
\end{equation}
\end{lemma}
Proof of Lemma~\ref{reductionof1} follows directly from the definition of reduction and Lemma IV.1 of \cite{AbdullaI}. 
 
The next lemma demonstrates a one-to-one mapping between $h$- and $\tilde{h}$-potentials, and their corresponding measures.
 \begin{lemma}\label{hpotentialmeasure}
 (i) \ Given a compact set $K\subset\mathbb{R}^{N+1}_+$, $\mathbb{P}^h_\mu$ is an $h$-potential of measure $\mu\in\mathcal{M}_K$ if and only if 
$\mathbb{P}^h_\mu(A^{-1}z)$ is an $\tilde{h}$-potential of measure $\mu(A^{-1})\in \mathcal{M}_{AK}$, i.e.
\begin{equation}\label{Atransformpotential}
\left\{
\begin{array}{l}
 K\subset\mathbb{R}^{N+1}_+ \ \mu\in\mathcal{M}_K: \ \mathbb{P}^h_\mu(A^{-1}z)=\mathbb{P}^{\tilde{h}}_{\mu(A^{-1})}(z), \ z\in\mathbb{R}^{N+1}_- \\
 K\subset\mathbb{R}^{N+1}_- \ \mu\in\mathcal{M}_K: \ \mathbb{P}^{\tilde{h}}_\mu(Az)=\mathbb{P}^{h}_{\mu(A)}(z), \ z\in\mathbb{R}^{N+1}_+
  \end{array}\right. 
\end{equation}
(ii) \ $\lambda\in\mathcal{M}_K$ is $h$-capacitary measure of compact $K\subset\mathbb{R}^{N+1}_+$ if and only if $\lambda(A^{-1})$ is $\tilde{h}$-capacitary measure of $AK$. \\
(iii) \ For a compact $K\subset\mathbb{R}^{N+1}_+$ with $h$-capacitary measure $\lambda\in \mathcal{M}_K$ we have
\begin{equation}\label{capAK}
C_{\tilde{h}}(AK)=\lambda(K)=C_h(K).
\end{equation}         
Conversely, for a compact $K\subset\mathbb{R}^{N+1}_-$ with $\tilde{h}$-capacitary measure $\tilde{\lambda}\in \mathcal{M}_K$ we have
\begin{equation}\label{capA-1K}
C_{h}(A^{-1}K)=\tilde{\lambda}(K)=C_{\tilde{h}}(K).
\end{equation}  
 \end{lemma}
 {\it Proof.} (i) \ Given compact $K\subset\mathbb{R}^{N+1}_+$ and $\mu\in\mathcal{M}_K$, we have
 \begin{equation}\label{Atransformhpotential0}
 \mathbb{P}^h_\mu(A^{-1}z)=\int\limits_{K}\frac{F(A^{-1}z-w)}{h(A^{-1}z)h_*(w)}d\mu(w)=\int\limits_{AK}\frac{F(A^{-1}z-A^{-1}w)}{h(A^{-1}z)h_*(A^{-1}w)}d\mu(A^{-1}w).
 \end{equation}
 Since
 \begin{gather}
 h(A^{-1}z)=\Big (-\frac{\pi}{t}\Big )^{-\frac{N}{2}}e^{\frac{|x+2t\gamma|^2}{4t}}=\Big (-\frac{\pi}{t}\Big )^{-\frac{N}{2}}e^{\frac{|x|^2}{4t}}\tilde{h}(z),\nonumber\\
 h_*(A^{-1}w)=(-4\pi\tau)^{\frac{N}{2}}e^{-\frac{|y+2\tau\gamma|^2}{4\tau}},\label{h*A-1w}\\
 F\Big (A^{-1}(z)-A^{-1}(w)\Big )
=  \Big (\frac{\pi(t-\tau)}{t\tau}\Big )^{-\frac{N}{2}} e^{-\frac{|\tau x-ty|^2}{4t\tau(t-\tau)}} 1_{\{t>\tau\}}(t)\nonumber
 \end{gather}
 by using the identity
\begin{equation}\label{ident}
\frac{|\tau x-ty|^2}{4t\tau(t-\tau)}-\frac{|y|^2}{4\tau}=\frac{|x-y|^2}{4(t-\tau)}-\frac{|x|^2}{4t},
\end{equation}
from \eqref{Atransformhpotential0} it follows that for all $z\in\mathbb{R}^{N+1}_-$
\begin{equation}
 \mathbb{P}^h_\mu(A^{-1}z)=\int\limits_{AK}\frac{F(z-w)}{\tilde{h}(z)\tilde{h}_*(w)}d\mu(A^{-1}w)=\mathbb{P}^{\tilde{h}}_{\mu(A^{-1})}(z)\nonumber
\end{equation}
Similarly, given compact $K\subset\mathbb{R}^{N+1}_-$ and $\mu\in\mathcal{M}_K$, we have
 \begin{equation}\label{Atransformhpotential3}
 \mathbb{P}^{\tilde{h}}_\mu(Az)=\int\limits_{K}\frac{F(Az-w)}{\tilde{h}(Az)\tilde{h}_*(w)}d\mu(w)=\int\limits_{A^{-1}K}\frac{F(Az-Aw)}{\tilde{h}(Az)\tilde{h}_*(Aw)}d\mu(Aw).
 \end{equation}
 Since
 \begin{gather}
 \tilde{h}(Az)=e^{\frac{2\langle \gamma, x\rangle -|\gamma|^2}{4t}}=(4\pi t)^{\frac{N}{2}}e^{\frac{|x|^2}{4t}}h(z),\nonumber\\
 \tilde{h}_*(Aw)=\Big (\frac{\pi}{\tau}\Big )^{-\frac{N}{2}}h_*(w)e^{-\frac{|y|^2}{4\tau}}\label{h*Aw}\\
 F(Az-Aw)
=  \Big (\frac{\pi(t-\tau)}{t\tau}\Big )^{-\frac{N}{2}} e^{-\frac{|\tau x-ty|^2}{4t\tau(t-\tau)}} 1_{\{t>\tau\}}(t),\nonumber
 \end{gather}
 by using the identity \eqref{ident}
from \eqref{Atransformhpotential3} it follows that for all $z\in\mathbb{R}^{N+1}_+$
\begin{equation}
 \mathbb{P}^{\tilde{h}}_\mu(Az)=\int\limits_{A^{-1}K}\frac{F(z-w)}{h(z)h_*(w)}d\mu(Aw)=\mathbb{P}^{h}_{\mu(A)}(z)\nonumber
\end{equation}
 (ii) \  This claim is a direct consequence of \eqref{Atransformpotential} and the definition of $h$- and $\tilde{h}$-capacities. \\
(iii) \  To establish \eqref{capAK}, note that from (ii) it follows 
\begin{equation}\label{capAK1}
C_{\tilde{h}}(AK)=\lambda(A^{-1})(AK)=\lambda(K)=C_h(K),
\end{equation}  
and similar proof applies to \eqref{capA-1K}.  \hfill{$\square$}

In the next lemma, we establish a one-to-one mapping between the heat potentials
\[ P^h_\mu(z)\equiv h\mathbb{P}^h_\mu(z)=\int\limits_{\mathbb{R}^{N+1}_+}\frac{F(z-w)}{h_*(w)}\,d\mu(w); \ P^{\tilde{h}}_\mu(z)\equiv \tilde{h}\mathbb{P}^{\tilde{h}}_\mu(z)=\int\limits_{\mathbb{R}^{N+1}_-}\frac{F(z-w)}{\tilde{h}_*(w)}\,d\mu(w) \]
\begin{lemma}\label{hpotentialmeasureold}
(i) \ For any fixed $w=(y,\tau)\in\mathbb{R}^{N+1}_+$, the Appell transform of $F(z-w), z=(x,t)\in \mathbb{R}^{N+1}_+$ is
\begin{equation}\label{AF(z-w)}
{\bf A}F(z-w)=(-4\pi\tilde{\tau})^{\frac{N}{2}}e^{-\frac{|\tilde{y}|^2}{4\tilde{\tau}}}F(z-\tilde{w}), \ z\in\mathbb{R}^{N+1}_-,
\end{equation}
where $\tilde{w}=Aw\in \mathbb{R}^{N+1}_-$. \\
(ii) \ For any fixed $w=(y,\tau)\in\mathbb{R}^{N+1}_-$, the inverse Appell transform of $F(z-w), z=(x,t)\in \mathbb{R}^{N+1}_-$ is
\begin{equation}\label{A-1F(z-w)}
{\bf A}^{-1}F(z-w)=\Big (\frac{\tilde{\tau}}{\pi}\Big )^{\frac{N}{2}}e^{-\frac{|\tilde{y}|^2}{4\tilde{\tau}}}F(z-\tilde{w}), \ z\in\mathbb{R}^{N+1}_+,
\end{equation}
where $\tilde{w}=A^{-1}w\in \mathbb{R}^{N+1}_+$. \\
(iii) \ Given a compact set $K\subset\mathbb{R}^{N+1}_+$, $P^h_\mu$ is a potential of measure $\frac{\mu}{h_*}\in\mathcal{M}_K$ if and only if 
${\bf A}P^h_\mu$ is a potential of measure $\frac{\mu(A^{-1})}{\tilde{h}_*}\in \mathcal{M}_{AK}$, i.e.
\begin{equation}\label{Atransformpotential1}
\left\{
\begin{array}{l}
 K\subset\mathbb{R}^{N+1}_+ \ \mu\in\mathcal{M}_K: \ {\bf A}P^h_\mu=P^{\tilde{h}}_{\mu(A^{-1})} \\
 K\subset\mathbb{R}^{N+1}_- \ \mu\in\mathcal{M}_K: \ {\bf A}^{-1}P^{\tilde{h}}_\mu=P^{h}_{\mu(A)} 
  \end{array}\right.
\end{equation}     
\end{lemma}
{\it Proof.} (i) \ We have
\begin{gather}
{\bf A}F(z-w)=\Big (-\frac{\pi}{t}\Big )^{\frac{N}{2}} e^{-\frac{|x|^2}{4t}}F\Big (A^{-1}(z)-A^{-1}(\tilde{w})\Big )\nonumber\\
= \Big (-\frac{\pi}{t}\Big )^{\frac{N}{2}} e^{-\frac{|x|^2}{4t}} \Big (\frac{\pi(t-\tilde{\tau})}{t\tilde{\tau}}\Big )^{-\frac{N}{2}} e^{-\frac{|\tilde{\tau}x-t\tilde{y}|^2}{4t\tilde{\tau}(t-\tilde{\tau})}} 1_{\{t>\tilde{\tau}\}}(t).\label{AF1}
\end{gather}
Using the identity \eqref{ident} 
from \eqref{AF1}, \eqref{AF(z-w)} follows. \\
(ii) \ We have
\begin{gather}
{\bf A}^{-1}F(z-w)=F(x,t) F\Big (A(z)-A(\tilde{w})\Big )=\nonumber\\
F(x,t) \Big (\frac{\pi(t-\tilde{\tau})}{t\tilde{\tau}}\Big )^{-\frac{N}{2}} e^{-\frac{|\tilde{\tau}x-t\tilde{y}|^2}{4t\tilde{\tau}(t-\tilde{\tau})}} 1_{\{t>\tilde{\tau}\}}(t).\label{A-11}
\end{gather}
Using the identity \eqref{ident}, from \eqref{A-11}, \eqref{A-1F(z-w)} follows.\\
(iii) \  Given a compact $K\subset\mathbb{R}^{N+1}_+$ and a measure $\mu\in\mathcal{M}_K$, from \eqref{AF(z-w)} it follows
\begin{equation}\label{APh1}
{\bf A}P^h_\mu(z)=\int\limits_K\frac{{\bf A}F(z-w)}{h_*(w)}d\mu(w)=\int\limits_{K}\frac{(-4\pi\tilde{\tau})^{\frac{N}{2}}e^{-\frac{|\tilde{y}|^2}{4\tilde{\tau}}}F(z-\tilde{w})}{h_*(w)}\,d\mu(w), 
\end{equation}
where $\tilde{w}=A(w)$. Changing the variable $A(w) \to w$ from \eqref{APh1} we deduce
\begin{equation}\label{APh2}
{\bf A}P^h_\mu(z)=\int\limits_{AK}\frac{(-4\pi\tau)^{\frac{N}{2}}e^{-\frac{|y|^2}{4\tau}}F(z-w)}{h_*(A^{-1}w)}\,d\mu(A^{-1}w).
\end{equation}
By \eqref{h*A-1w} we have 
\[ {\bf A}P^h_\mu(z)= \int\limits_{AK} \frac{e^{-\frac{|y|^2}{4\tau}}F(z-w)}{e^{-\frac{|y+2\tau\gamma|^2}{4t}}}\,d\mu(A^{-1}w)= \int\limits_{AK}\frac{F(z-w)}{\tilde{h}_*(w)}\,d\mu(A^{-1}w)=P^{\tilde{h}}_{\mu(A^{-1})}(z). \]
On the contrary, given a compact $K\subset\mathbb{R}^{N+1}_-$ and a measure $\mu\in\mathcal{M}_K$, from \eqref{A-1F(z-w)} it follows
\begin{equation}\label{A-1Ph1}
 {\bf A}^{-1}P^{\tilde{h}}_{\mu}(z)=\int\limits_K\frac{{\bf A}^{-1}F(z-w)}{\tilde{h}_*(w)}d\mu(w)=\int\limits_{K}\Big (\frac{\tilde{\tau}}{\pi}\Big)^{\frac{N}{2}}e^{-\frac{|\tilde{y}|^2}{4\tilde{\tau}}}\frac{F(z-\tilde{w})}{\tilde{h}(w)}\,d\mu(w), 
\end{equation}
where $\tilde{w}=A^{-1}(w)$. Changing the variable $A^{-1}(w)\to w$ from \eqref{A-1Ph1} we deduce
\begin{equation}\label{A-1Ph2}
 {\bf A}^{-1}P^{\tilde{h}}_\mu(z)=\int\limits_{A^{-1}K}\Big (\frac{\tau}{\pi}\Big)^{\frac{N}{2}}e^{-\frac{|y|^2}{4\tau}}\frac{F(z-w)}{\tilde{h}_*(Aw)}\,d\mu(Aw).
\end{equation}
By \eqref{h*Aw} we have 
\[ {\bf A}^{-1}P^{\tilde{h}}_\mu(z)=\int\limits_{A^{-1}K}\frac{F(z-w)}{h_*(w)}\,d\mu(Aw)=P^{h}_{\mu(A)}(z), \]
 which completes the proof of \eqref{Atransformpotential1}. 
\hfill{$\square$}

%%%%%%%

Next lemmas express the properties of $h$-capacity and $h$-capacitary potentials.
\begin{lemma}\label{hpotentialvsreduction1}
(i) \ For any compact $K\subset \mathbb{R}^{N+1}_+$ there is a unique $h$-capacitary measure $\lambda\in\mathcal{M}_{\partial K}$ such that
\begin{equation}\label{hpotentialvsreduction}
\ ^h\hat{R}^{K}_1 \equiv \mathbb{P}^h_\lambda, \ \  \hat{R}^{K}_h \equiv P^h_\lambda \ \ on \ \ \mathbb{R}^{N+1}_+ 
\end{equation}
where $P^h_\lambda$ be a heat potential with associated Riesz measure $\frac{\lambda}{h_*}$. \\
(ii) \ $^h\hat{R}^{K}_1$ is $h$-superparabolic in $\mathbb{R}^{N+1}_+$, and $h$-parabolic in $\mathbb{R}^{N+1}_+\setminus \partial K$; \ $ \hat{R}^{K}_h \equiv P^h_\lambda$ is 
superparabolic in $\mathbb{R}^{N+1}_+$, and parabolic in $\mathbb{R}^{N+1}_+\setminus \partial K$.\\  
(iii) \ Furthermore,
\begin{equation}\label{hpotentialdistribution}
\mathcal{H}[P^h_\lambda]\equiv\mathcal{H}[ h  \mathbb{P}^h_\lambda]=\frac{\lambda}{h_*} \ \text{in the sense of distributions on} \ \mathbb{R}^{N+1}_+ 
\end{equation}
(iv) \ As a set function defined in the class of compact subsets of $\mathbb{R}^{N+1}_+$, $C_h(\cdot)$ is a topological precapacity. i.e.
\begin{enumerate}
\item $C_h(\cdot)$ is strongly subadditive:
\begin{enumerate}
\item $C_h(K_1) \leq C_h(K_2)$ if $K_1 \subset K_2$;
\item $C_h(K_1\cup K_2)+C_h(K_1\cap K_2) \leq C_h(K_2)+C_h(K_2)$;
\end{enumerate}
\item If $K_n$ is a monotone sequence of compact sets with compact limit $K$, then
\[  \lim_{n\to +\infty} C_h(K_n)=C_h(K) \]
\end{enumerate}
\end{lemma}   
Similar result in the $(\mathbb{R}^{N+1}_-,\tilde{h})$ setting reads:
\begin{lemma}\label{tildehpotentialvsreduction1}
(i) \ For any compact $K\subset \mathbb{R}^{N+1}_-$ there is a unique $\tilde{h}$-capacitary measure $\lambda\in\mathcal{M}_{\partial K}$ such that
\begin{equation}\label{tildehpotentialvsreduction}
\ ^{\tilde{h}}\hat{R}^{K}_1 \equiv \mathbb{P}^{\tilde{h}}_\lambda, \ \ \hat{R}^{K}_{\tilde{h}} \equiv P^{\tilde{h}}_\lambda \ \ on \ \ \mathbb{R}^{N+1}_- 
\end{equation}
where $P^{\tilde{h}}_\lambda$ be a heat potential with associated Riesz measure $\frac{\lambda}{\tilde{h}_*}$. \\
(ii) \ $^{\tilde{h}}\hat{R}^{K}_1$ is $\tilde{h}$-superparabolic in $\mathbb{R}^{N+1}_-$, and $\tilde{h}$-parabolic in $\mathbb{R}^{N+1}_-\setminus \partial K$; \ $ \hat{R}^{K}_{\tilde{h}} \equiv P^{\tilde{h}}_\lambda$ is 
superparabolic in $\mathbb{R}^{N+1}_-$, and parabolic in $\mathbb{R}^{N+1}_-\setminus \partial K$.\\
(iii) \ Furthermore,
\begin{equation}\label{tildehpotentialdistribution}
\mathcal{H}[P^{\tilde{h}}_\lambda]\equiv\mathcal{H}[ \tilde{h} \mathbb{P}^{\tilde{h}}_\lambda]=\frac{\lambda}{\tilde{h}_*} \ \text{in the sense of distributions on} \ \mathbb{R}^{N+1}_-. 
\end{equation}
(iv) \ As a set function defined in the class of compact subsets of $\mathbb{R}^{N+1}_-$, $C_{\tilde{h}}(\cdot)$ is a topological precapacity. i.e.
\begin{enumerate}
\item $C_{\tilde{h}}(\cdot)$ is strongly subadditive:
\begin{enumerate}
\item $C_{\tilde{h}}(K_1) \leq C_{\tilde{h}}(K_2)$ if $K_1 \subset K_2$;
\item $C_{\tilde{h}}(K_1\cup K_2)+C_{\tilde{h}}(K_1\cap K_2) \leq C_{\tilde{h}}(K_2)+C_{\tilde{h}}(K_2)$;
\end{enumerate}
\item If $K_n$ is a monotone sequence of compact sets with compact limit $K$, then
\[  \lim_{n\to +\infty} C_{\tilde{h}}(K_n)=C_{\tilde{h}}(K) \]
\end{enumerate}
\end{lemma}   
%%%%%%
{\it Proof}. Proofs of Lemmas~\ref{hpotentialvsreduction1} and ~\ref{tildehpotentialvsreduction1} is based on the following scheme. 
\begin{equation}\label{lemmasproofscheme}
Lem.~\ref{tildehpotentialvsreduction1}|_{\gamma=0} \stackrel{A^{-1}}{\Longrightarrow} Lem.~\ref{hpotentialvsreduction1}|_{\gamma=0} \stackrel{x\mapsto x+\gamma}{\Longrightarrow} Lem.~\ref{hpotentialvsreduction1}|_{\gamma\neq 0}\stackrel{A}{ \Longrightarrow} Lem.~\ref{tildehpotentialvsreduction1}|_{\gamma\neq 0}   
\end{equation}
The result of Lemma~\ref{tildehpotentialvsreduction1} with $\gamma=0$ is well known\cite{Doob,watson}. Note that in this case $\tilde{h}\equiv 1$, $\mathbb{P}^{\tilde{h}}_\lambda \equiv P^{\tilde{h}}_\lambda$ be heat potential, and $C_{\tilde{h}}(K)\equiv cap(K)$ be a classical thermal capacity of a compact $K$. As it is indicated in the first chain of \eqref{lemmasproofscheme}, the claims of the Lemma~\ref{tildehpotentialvsreduction1} with $\gamma=0$ imply the corresponding claims of the Lemma~\ref{hpotentialvsreduction1} with $\gamma=0$ through application of the map $A^{-1}$ of \eqref{Appellmap} and using the results of Lemma~\ref{hpotentialmeasure}, Lemma~\ref{hpotentialmeasureold} (iii), and Remarks~\ref{hred1e1redh} and ~\ref{hred1e1redh12}. The statements of Lemma~\ref{hpotentialvsreduction1} with $\gamma\neq 0$ directly follow from the corresponding statements of Lemma~\ref{hpotentialvsreduction1} with $\gamma=0$ by using the translation map $x\mapsto x+\gamma$. Finally, the statements of the Lemma~\ref{hpotentialvsreduction1} with $\gamma\neq 0$ are translated to the corresponding statements of the 
Lemma~\ref{tildehpotentialvsreduction1} with $\gamma\neq 0$ through application of the map $A$ of \eqref{Appellmap} and using again the results of Lemma~\ref{hpotentialmeasure}, Lemma~\ref{hpotentialmeasureold} (iii), and Remarks~\ref{hred1e1redh} and ~\ref{hred1e1redh12}. Note that the properties \eqref{hpotentialdistribution} and \eqref{tildehpotentialdistribution} follow directly from the properties of heat potentials \cite{Doob,watson}, since $h  \mathbb{P}^h_\lambda= P^h_\lambda$ and
$ \tilde{h} \mathbb{P}^{\tilde{h}}_\lambda =P^{\tilde{h}}_\lambda$ are heat potentials with corresponding Riesz measures $\frac{\lambda}{h}$ and $\frac{\lambda}{\tilde{h}}$ respectively.

{\it Alternative proof of  \eqref{hpotentialdistribution} and \eqref{tildehpotentialdistribution}:} By choosing $u=\mathbb{P}^h_\lambda$ in \eqref{HuA-1}, and taking into account \eqref{Atransformpotential} we have
\begin{equation}\label{HuA-1altproof1}
\frac{\mathcal{H}[h(z)\mathbb{P}^h_\lambda(z)]}{h(z)}\Big |_{z=A^{-1}(w)}=4\tau^2\frac{\mathcal{H}[\tilde{h}(w)\mathbb{P}^{\tilde{h}}_{\lambda(A^{-1})}(w)]}{\tilde{h}(w)}.
\end{equation}
If $\gamma=0$ ($\tilde{h}\equiv 1$), then being a heat potential $\mathbb{P}^{\tilde{h}}_{\lambda(A^{-1})}$ satisfies
\begin{equation}\label{altproof2}
\mathcal{H}[\mathbb{P}^{\tilde{h}}_{\lambda(A^{-1})}(w)]=\lambda(A^{-1}), \ \text{in the sense of distributions on} \ \mathbb{R}^{N+1}_-,
\end{equation}
that is to say,
\begin{equation}\label{altproof3}
\int\limits_{\mathbb{R}^{N+1}_-}\mathbb{P}^{\tilde{h}}_{\lambda(A^{-1})}(w) \mathcal{H}^*f(w)\,dw=\int\limits_{\mathbb{R}^{N+1}_-}f(w)\,d\lambda(A^{-1}w), \ \forall f\in C_0^\infty(\mathbb{R}^{N+1}_-).
\end{equation}
Assumung that $\mathbb{P}^{h}_{\lambda}$ is smooth and using \eqref{HuA-1altproof1},\eqref{altproof3} we derive
\begin{gather}
\int\limits_{\mathbb{R}^{N+1}_+}h(z)\mathbb{P}^{h}_{\lambda}(z)\mathcal{H}^*[f(z)]\,dz= \int\limits_{\mathbb{R}^{N+1}_+}\mathcal{H}[h(z)\mathbb{P}^{h}_{\lambda}(z)] f(z)\,dz\nonumber\\ =\int\limits_{\mathbb{R}^{N+1}_-}\mathcal{H}[h(z)\mathbb{P}^{h}_{\lambda}(z)]|_{z=A^{-1}(w)}f(A^{-1}w) (-2\tau)^{-N-2}\,dw\nonumber\\
=\int\limits_{\mathbb{R}^{N+1}_-}\mathcal{H}[\mathbb{P}^{\tilde{h}}_{\lambda(A^{-1})}](w)h(A^{-1}w)f(A^{-1}w)(-2\tau)^{-N}\,dw\nonumber\\ =\int\limits_{\mathbb{R}^{N+1}_-}\mathbb{P}^{\tilde{h}}_{\lambda(A^{-1})}(w) \mathcal{H}^*[h(A^{-1}w)f(A^{-1}w)(-2\tau)^{-N}]\,dw\nonumber\\ =\int\limits_{\mathbb{R}^{N+1}_-}h(A^{-1}w)f(A^{-1}w)(-2\tau)^{-N}\,d\lambda(A^{-1}w)=\int\limits_{\mathbb{R}^{N+1}_+}(2t)^Nh(z)f(z)\,d\lambda(z)\nonumber\\
=\int\limits_{\mathbb{R}^{N+1}_+}f(z)\frac{\,d\lambda(z)}{h_*(z)}\,d\l, \ \forall f\in C_0^\infty(\mathbb{R}^{N+1}_+).\label{altproof4}
\end{gather}
Using approximation argument, \eqref{altproof4} is valid for $h\mathbb{P}^{h}_{\lambda}\in L^1_{loc}(\mathbb{R}^{N+1}_+)$ which implies \eqref{hpotentialdistribution}. 
The claim \eqref{hpotentialdistribution} in the case $\gamma\neq 0$ easily follows using a translation map $x\mapsto x+\gamma$. 

To prove \eqref{tildehpotentialdistribution}, choose $u=\mathbb{P}^{\tilde{h}}_\lambda$ in \eqref{HuA}, and take into account \eqref{Atransformpotential} to derive
\begin{equation}\label{HuAaltproof1}
\frac{\mathcal{H}[\tilde{h}(w)\mathbb{P}^{\tilde{h}}_\lambda(w)]}{\tilde{h}(w)}\Big |_{w=A(z)}=4t^2\frac{\mathcal{H}[h(z)\mathbb{P}^{h}_{\lambda(A)}(z)]}{h(z)}.
\end{equation}
From \eqref{hpotentialdistribution} it follows that
\begin{equation}\label{haltproof2}
\mathcal{H}[h(z)\mathbb{P}^{h}_{\lambda(A)}(z)]=\frac{\lambda(A)}{h_*} \ \text{in the sense of distributions on} \ \mathbb{R}^{N+1}_+,
\end{equation}
that is to say,
\begin{equation}\label{haltproof3}
\int\limits_{\mathbb{R}^{N+1}_+}h(z)\mathbb{P}^{h}_{\lambda(A)}(z) \mathcal{H}^*f(z)\,dz=\int\limits_{\mathbb{R}^{N+1}_+}\frac{f(z)}{h(z)}\,d\lambda(Az), \ \forall f\in C_0^\infty(\mathbb{R}^{N+1}_+).
\end{equation}
Assumung that $\mathbb{P}^{\tilde{h}}_{\lambda}$ is smooth and using \eqref{HuAaltproof1},\eqref{haltproof3} we derive
\begin{gather}
\int\limits_{\mathbb{R}^{N+1}_-}\tilde{h}(w)\mathbb{P}^{\tilde{h}}_{\lambda}(w)\mathcal{H}^*[f(w)]\,dw= \int\limits_{\mathbb{R}^{N+1}_-}\mathcal{H}[\tilde{h}(w)\mathbb{P}^{\tilde{h}}_{\lambda}(w)] f(w)\,dw\nonumber\\ =\int\limits_{\mathbb{R}^{N+1}_+}\mathcal{H}[\tilde{h}(w)\mathbb{P}^{\tilde{h}}_{\lambda}(w)]|_{w=A(z)}f(Az) (-2t)^{-N-2}\,dz\nonumber\\
=\int\limits_{\mathbb{R}^{N+1}_+}\tilde{h}(Az)\frac{\mathcal{H}[h(z)\mathbb{P}^{h}_{\lambda(A)}(z)]}{h(z)}(2t)^{-N}f(Az)\,dz\nonumber\\ =\int\limits_{\mathbb{R}^{N+1}_+}h(z)\mathbb{P}^{h}_{\lambda(A)}(z) \mathcal{H}^*[\tilde{h}(Az)h^{-1}(z)f(Az)(2t)^{-N}]\,dz\nonumber\\ =\int\limits_{\mathbb{R}^{N+1}_+}\tilde{h}(Az)f(Az)(2t)^{-N}\,\frac{\,d\lambda(Az)}{h(z)h^*(z)}\nonumber\\
=\int\limits_{\mathbb{R}^{N+1}_-}\frac{(-2\tau)^{N}\tilde{h}(w)f(w)}{h(A^{-1}w)h_*(A^{-1}w)}\,d\lambda(w)
=\int\limits_{\mathbb{R}^{N+1}_-}f(w)\frac{\,d\lambda(w)}{\tilde{h}_*(w)}, \ \forall f\in C_0^\infty(\mathbb{R}^{N+1}_-).\label{haltproof4}
\end{gather}
Using approximation argument, \eqref{haltproof4} is valid for $\tilde{h}\mathbb{P}^{\tilde{h}}_{\lambda}\in L^1_{loc}(\mathbb{R}^{N+1}_-)$ which implies \eqref{tildehpotentialdistribution}. \hfill{$\square$}

As a topological precapacity, $C_h$ has an extension to countably strongly subadditive Choquet capacity $C_h(\cdot)$ relative to the class of compact subsets of $\mathbb{R}^{N+1}_+$. Recall that the set function $C: 2^{\mathbb{R}^{N+1}_+} \to [0,+\infty]$ is called a Choquet capacity if it satisfies the conditions $(1a),(2)$ in a statement (iii) of the Lemma~ \ref{hpotentialvsreduction1}. Extend $C_h(\cdot)$ to $2^{\mathbb{R}^{N+1}_+}$ as follows:
\begin{gather}
C_h(A)=\sup \{ C_h(K): K \quad\text{is compact}, \ K\subset A\}, \ A\subset \mathbb{R}^{N+1}_+ \quad\text{is open} \label{extopen}\\
C_h(A)=\inf \{ C_h(B): B \quad\text{is open}, \ A\subset B\}, \ A\subset \mathbb{R}^{N+1}_+ \quad\text{is arbitrary} \label{extarb}
\end{gather}
Likewise, $C_{\tilde{h}}$ is extended to $2^{\mathbb{R}^{N+1}_-}$.
The celebrated Choquet capacitability theorem states:
\begin{lemma}\label{Choquet}\cite{Choquet} The set function $C_h(\cdot): 2^{\mathbb{R}^{N+1}_+} \to [0,+\infty]$ (or \ $C_{\tilde{h}}(\cdot): 2^{\mathbb{R}^{N+1}_-} \to [0,+\infty] $) is a countably strongly subadditive Choquet capacity
defined on all subsets of $\mathbb{R}^{N+1}_+$ (or $\mathbb{R}^{N+1}_-$). All the Borel (and even analytic) subsets $E \subset \mathbb{R}^{N+1}_+$ (or $E \subset \mathbb{R}^{N+1}_-$) are capacitable, i.e.
\[ C_h(E)=\sup \{ C_h(K): K\quad\text{is compact}, \ K\subset E\}\]
\[ ( \text{or} \ C_{\tilde{h}}(E)=\sup \{ C_{\tilde{h}}(K): K\quad\text{is compact}, \ K\subset E\} ) \]
\end{lemma}
%%%%%%%%%

\section{Averaiging Property of $h$-Parabolic Functions, and Harnack Estimates}\label{averaging}

The goal of this section is to prove some essential geometric estimates of the $h$-heat balls, averaging property of $h$-parabolic functions and Harnack estimates on these sets.

The next lemma expresses some key geometric estimates of the $h$-heat balls, and the averaging property of $h$-parabolic functions on these sets.
Consider an $h$-heat ball of center $z_0=(\gamma,t_0)$, $t_0>0$ and radius $4t_0^2c(1+4t_0c)^{-1}$ with $c>0$ (Figure 8):
\begin{gather}
B(z_0,c)=\{z\in\mathbb{R}^{N+1}_+: \frac{F(z_0-z)}{h(z_0)h_*(z)}>(4\pi c)^{-\frac{N}{2}}\} \nonumber\\
\equiv \Big \{z\in\mathbb{R}^{N+1}_+: \ |x-\gamma|<R_c(t), \ \frac{t_0}{1+4t_0c}<t<t_0\Big \}, \nonumber
\end{gather}
where 
\[ R_c(t):=\Big [\frac{2N}{t_0}t(t-t_0)\log \frac{t_0-t}{4tt_0c}\Big ]^{\frac{1}{2}}\]
is the radius of the spherical cros section of $B(z_0,c)$ at $t$. Note that 
\[ R_c(t_0)=R_c\Big ( \frac{t_0}{1+4t_0c} \Big )=0, \]
and it attains its maximum at interior of the interval $((1+4t_0c)^{-1}t_0, t_0)$. 
\begin{figure}[ht]\label{fig3}
    \centering
    \includegraphics[width=0.9\textwidth, trim=1cm 8cm 1cm 7cm, clip]{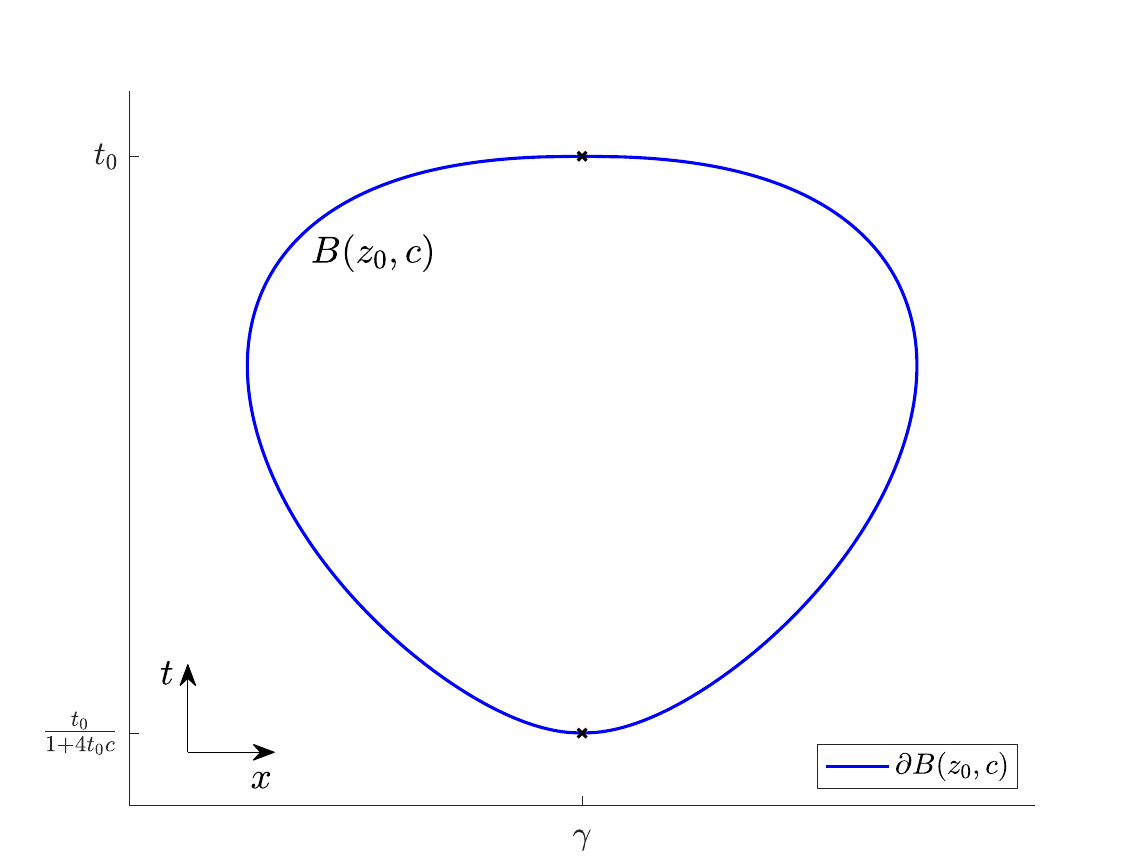}
    \caption{Sketch of $h$-heat ball $B(z_0,c)$ for $N=1$.}
\end{figure}
Hence, $h$-heat ball $B(z_0,c)$ is a bounded open set surrounded by the level surface of the $h_*$-parabolic function 
\[ \frac{F(z_0-z)}{h(z_0)h_*(z)}.\] 
\begin{lemma}\label{avprophparabolic}
If $u:\mathbb{R}^{N+1}_+\to \mathbb{R}$ is smooth, then\\
(i) \ the function
\begin{equation}\label{phih}
 \phi(c)=\frac{t_0^2}{(4c)^{\frac{N}{2}}}\int\limits_{B(z_0,c)}u(x,t) \frac{|x-\gamma|^2}{t^{N+2}|t-t_0|^2}\,dx\,dt
 \end{equation}
is differentiable for $c>0$, and
\begin{equation}\label{phiprimeh}
 \phi'(c)=\frac{Nt_0}{(4c)^{\frac{N+1}{2}}}\int\limits_{B(z_0,c)}\frac{\mathcal{H}[h(x,t)u(x,t)]}{h(x,t)}\frac{R_c^2(t)-|x-\gamma|^2}{t^{N+1}(t-t_0)}\,dx\,dt.
 \end{equation}
(ii) \ There exists a positive constant $C$ such that if $\mathcal{H}[h(x,t)u(x,t)]\leq 0$ in $B(z_0,2c)$, $c>0$, then
\begin{equation}\label{phi2c-phich} 
 \phi(2c)-\phi(c)\geq \frac{C}{c^{\frac{N}{2}}}\int\limits_{B(z_0,c/2)}\frac{-\mathcal{H}[h(x,t)u(x,t)]}{(2t)^Nh(x,t)}\,dx\,dt. 
 \end{equation}
(iii) If $u$ is $h$-parabolic in an open subset of $\mathbb{R}^{N+1}_+$ containing $\overline{B(z_0,c)}$, then it satisfies the averaging formula
\begin{equation}\label{hparaverage}
 u(\gamma,t_0)=\frac{4t_0^2}{(\pi c)^{\frac{N}{2}}}\int\limits_{B(z_0,c)}u(x,t) \frac{|x-\gamma|^2}{(4t)^{N+2}|t-t_0|^2}\,dx\,dt.
 \end{equation}
\end{lemma}

To formulate similar result in the $(\mathbb{R}^{N+1}_-,\tilde{h})$ setting, consider
$\tilde{h}$-{\it heat ball} of center $\overline{z}=(-2\gamma\tau_0, \tau_0)$, $\tau_0<0$ and radius $c>0$ (Figure 9):
\begin{gather}
\tilde{B}(\overline{z},c)=\{w\in\mathbb{R}^{N+1}_-: \frac{F(\bar{z}-w)}{\tilde{h}(\bar{z})\tilde{h}_*(w)}>(4\pi c)^{-\frac{N}{2}}\} \nonumber\\
\equiv \{w=(y,\tau)\in\mathbb{R}^{N+1}_-: \ |y+2\tau\gamma|^2<\tilde{R}_c^2(\tau), \ \tau_0-c<\tau<\tau_0\Big \},
\end{gather}
where
\[ \tilde{R}_c(\tau):=\Big [2N(\tau-\tau_0)\log \frac{\tau_0-\tau}{c}\Big ]^{\frac{1}{2}}\]
is the radius of the spherical cross section of $\tilde{B}(\overline{z},c)$ at $\tau$. Note that 
\[ \tilde{R}_c(\tau_0)=\tilde{R}_c(\tau_0-c)=0, \]
and it attains its maximum value $\Big(\frac{2Nc}{e}\Big)^{\frac{1}{2}}$ at the unique point $\tau=\tau_0-\frac{c}{e}$. 
\begin{figure}[ht]\label{fig4}
    \centering
    \includegraphics[width=0.9\textwidth, trim=1cm 8cm 1cm 7cm, clip]{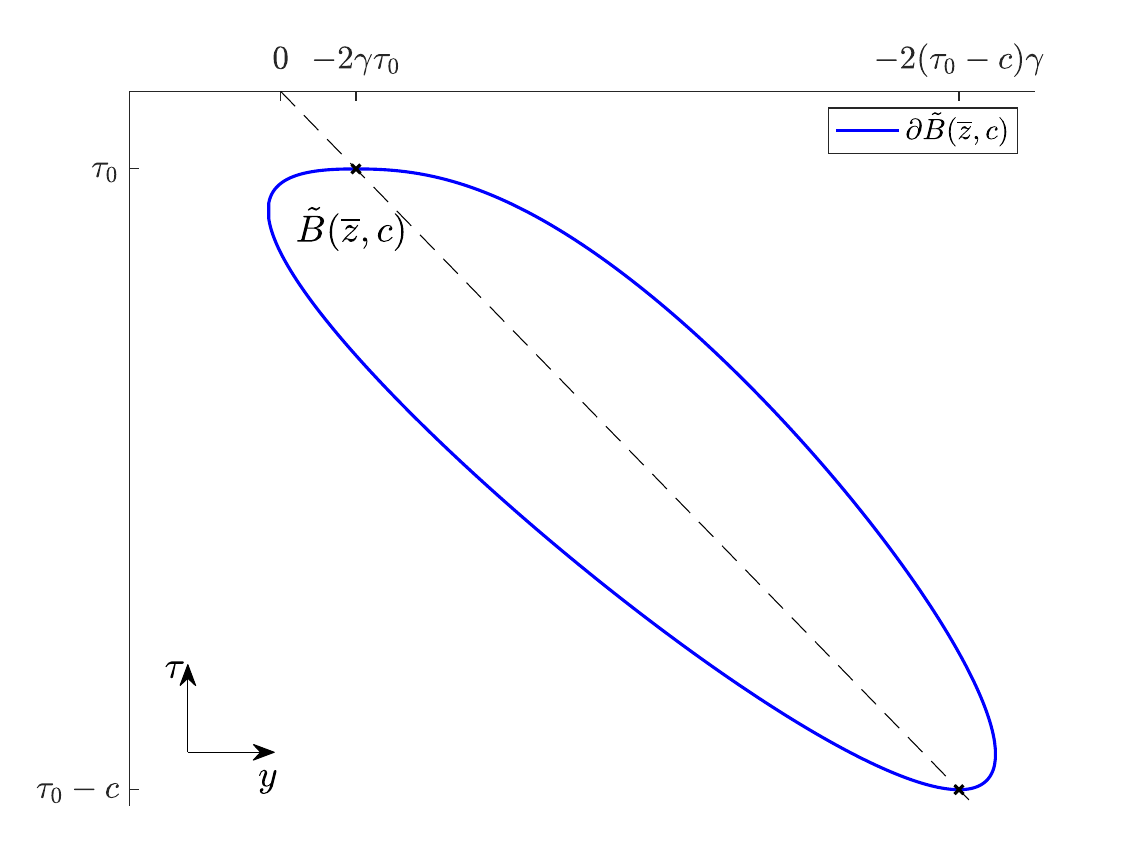}
    \caption{Sketch of $\tilde{h}$-heat ball $\tilde{B}(\bar{z}_0,c)$ for $N=1$.}
\end{figure}
Hence, $\tilde{h}$-heat ball $\tilde{B}(\bar{z},c)$ is a bounded open set surrounded by the level surface of the $\tilde{h}_*$-parabolic function 
\[ \frac{F(\bar{z}-w)}{\tilde{h}(\bar{z})\tilde{h}_*(w)}.\] 

\begin{lemma}\label{avprophtildeparabolic}
If $u:\mathbb{R}^{N+1}_-\to \mathbb{R}$ is smooth, then\\
(i) \ the function
\begin{equation}\label{phitildeh}
 \tilde{\phi}(c)=\frac{1}{c^{\frac{N}{2}}}\int\limits_{\tilde{B}(\overline{z},c)}u(y,\tau) \frac{|y+2\tau\gamma|^2}{|\tau-\tau_0|^2}\,dy\,d\tau 
 \end{equation}
is differentiable for $c>0$, and
\begin{equation}\label{phiprimetildeh}
 \tilde{\phi}'(c)=\frac{N}{2c^{\frac{N+1}{2}}}\int\limits_{\tilde{B}(\overline{z}_0,c)}\frac{\mathcal{H}[\tilde{h}(y,\tau)u(y,\tau)]}{\tilde{h}(y,\tau)}\frac{\tilde{R}_c^2(\tau)-|y+2\gamma \tau|^2}{\tau-\tau_0}\,dy\,d\tau.
 \end{equation}
(ii) \ There exists a positive constant $C$ such that if $\mathcal{H}[\tilde{h}(x,t)u(x,t)]\leq 0$ in $\tilde{B}(\overline{z},2c)$, $c>0$, then
\begin{equation}\label{phi2c-phictildeh}
 \tilde{\phi}(2c)-\tilde{\phi}(c)\geq \frac{C}{c^{\frac{N}{2}}}\int\limits_{\tilde{B}(\overline{z},c/2)}\frac{-\mathcal{H}[\tilde{h}(y,\tau)u(y,\tau)]}{\tilde{h}(y,\tau)}\,dy\,d\tau. 
 \end{equation}
(iii) If $u$ is $\tilde{h}$-parabolic in an open subset of $\mathbb{R}^{N+1}_-$ containing $\overline{\tilde{B}(\overline{z},c)}$, then it satisfies the averaging formula
\begin{equation}\label{tildehparaverage}
 u(-2\gamma \tau_0,\tau_0)=\frac{1}{2^{N+1}(\pi c)^{\frac{N}{2}}}\int\limits_{\tilde{B}(\overline{z},c)}u(y,\tau) \frac{|y+2\gamma \tau|^2}{|\tau-\tau_0|^2}\,dy\,d\tau.
 \end{equation}
\end{lemma}
{\it Proof}. Proofs of Lemmas~\ref{avprophparabolic} and ~\ref{avprophtildeparabolic} are based on the following scheme. 
\begin{equation}\label{lemmasproofscheme2}
Lem.~\ref{avprophtildeparabolic}|_{\gamma=0} \stackrel{A^{-1}}{\Longrightarrow} Lem.~\ref{avprophparabolic}|_{\gamma=0} \stackrel{x\mapsto x+\gamma}{\Longrightarrow} Lem.~\ref{avprophparabolic}|_{\gamma\neq 0}\stackrel{A}{ \Longrightarrow} Lem.~\ref{avprophtildeparabolic}|_{\gamma\neq 0}   
\end{equation}
Lemma~\ref{avprophtildeparabolic} with $\gamma=0$ is well known (see Lemma 3.1 in \cite{evansgariepy} and \cite{watson}). Note that in this case $\tilde{h}\equiv 1$, and $\tilde{B}(\overline{z},c)$ is an heat ball of center $\overline{z}$ and radius $c$. As it is indicated in the first chain of \eqref{lemmasproofscheme2}, the claims of the Lemma~\ref{avprophtildeparabolic} with $\gamma=0$ imply the corresponding claims of the Lemma~\ref{avprophparabolic} with $\gamma=0$ through application of the map $A^{-1}$ of \eqref{Appellmap}. Indeed, given a smooth function $u:\mathbb{R}^{N+1}_+\to \mathbb{R}$, it is not difficult to verify that \eqref{phitildeh} with $u(y,\tau)$ and $\tilde{B}(\overline{z},c)$ replaced with $u(A^{-1}(y,\tau))$ and $\tilde{B}(Az_0,c)$ respectively, is transformed to \eqref{phih} after the change of variable $w=Az$ in the integral. Therefore to prove \eqref{phiprimeh}, \eqref{phi2c-phich}, \eqref{hparaverage} we only need to pursue change of variable $w=Az$ in the corresponding formulae \eqref{phiprimetildeh}, \eqref{phi2c-phictildeh}, \eqref{tildehparaverage}. To change the variable in the term with heat operator under the integrals \eqref{phiprimetildeh}, \eqref{phi2c-phictildeh} we use the formula \eqref{HuA} with $\tilde{h}\equiv 1$ and $h\equiv F$. 
Using \eqref{HuA}, the formulae \eqref{phiprimeh}, \eqref{phi2c-phich}, \eqref{hparaverage} follow from the corresponding formulae \eqref{phiprimetildeh}, \eqref{phi2c-phictildeh}, \eqref{tildehparaverage} via the change of variable $w=Az$. Hence, Lemma~\ref{avprophparabolic} is proved if $\gamma=0$.

The claims of Lemma~\ref{avprophparabolic} with $\gamma\neq 0$ easily follow from the case $\gamma=0$ by applying the translation $x \mapsto x+\gamma$, and noting that the function $u:\mathbb{R}^{N+1}_+\to \mathbb{R}$ is $F$-parabolic (or $F$-super- or subparabolic) if and only if the the function $u(x-\gamma,t)$ is $h$-parabolic (or $h$-super- or subparabolic). 

Finally, we prove that Lemma~\ref{avprophtildeparabolic} with $\gamma\neq 0$ follows from the Lemma~\ref{avprophparabolic}. First, we easily verify that given a smooth function $u:\mathbb{R}^{N+1}_-\to \mathbb{R}$, \eqref{phih} with $u(x,t)$ and $B(z_0,c)$ replaced with $u(A(x,t))$ and $B(A^{-1}\overline{z},c)$ respectively, is transformed to \eqref{phitildeh} after the change of variable $z=A^{-1}w$ in the integral. Therefore to prove \eqref{phiprimetildeh}, \eqref{phi2c-phictildeh}, \eqref{tildehparaverage} we only need to pursue change of variable $z=A^{-1}w$ in the corresponding formulae \eqref{phiprimeh}, \eqref{phi2c-phich}, \eqref{hparaverage}. To change the variable in the term with the heat operator under the integrals \eqref{phiprimeh}, \eqref{phi2c-phich} we use the formula \eqref{HuA-1}. 
Using \eqref{HuA-1}, the formulae \eqref{phiprimetildeh}, \eqref{phi2c-phictildeh}, \eqref{tildehparaverage} follow from the corresponding formulae \eqref{phiprimeh}, \eqref{phi2c-phich}, \eqref{hparaverage} via the change of variable $z=A^{-1}w$. Hence, Lemma~\ref{avprophtildeparabolic} is proved if $\gamma\neq 0$. \hfill{$\square$}

In the next two lemmas, we formulate important Harnack-type estimates for the $h$- and $\tilde{h}$-parabolic functions. Let 
\[Q(z_0,c)=\Big \{z\in\mathbb{R}^{N+1}_+: \ |x-\gamma|<R_c(t), \ \frac{t_0}{1+3t_0c}<t<t_0\Big \}, \]
be an $h$-heat ball $B(z_0,c)$ with the region below $t=t_0(1+3t_0c)^{-1}$ removed. 
\begin{lemma}\label{harnackhparabolic}
Let $u\geq 0$ be an $h$-parabolic function in $Q(z_0,2c), c>0$, and suppose that $u$ is continuous at each point 
of $\partial Q(z_0,2c)$, except possibly at $z_0$. Then there exists a positive constant $C$, depending only on $N$ and not on c, such that (Figure 10)
\begin{equation}\label{harnackhpar}
\dashint_{|x|\leq \frac{(3Nc)^{\frac{1}{2}}t_0}{1+6ct_0}}u\Big ( x,\frac{t_0}{1+6ct_0}\Big )\,dx \leq C \inf\limits_{B(z_0,\frac{3c}{4})} u
\end{equation}
\begin{figure}[ht]\label{fig5}
    \centering
    \includegraphics[width=0.9\textwidth, trim=1cm 8cm 1cm 7cm, clip]{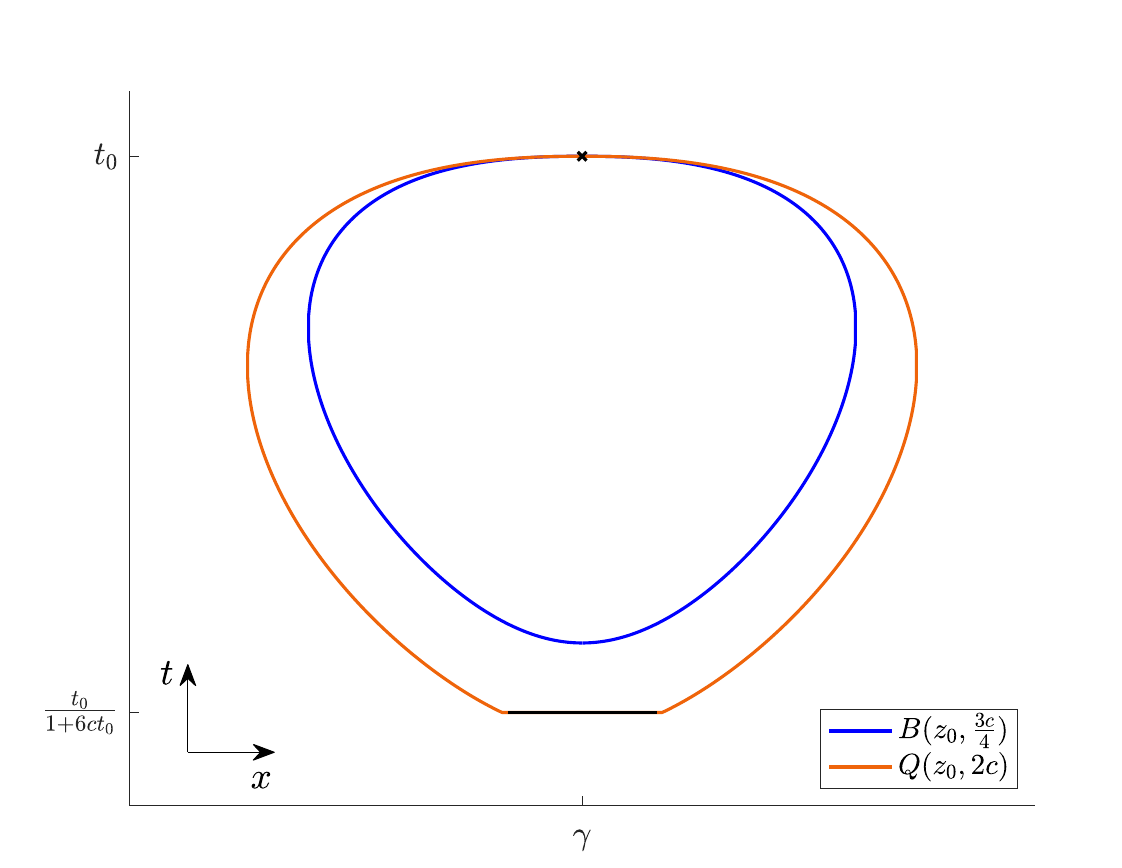}
    \caption{Sketch of $h$-heat ball $B(z_0,\frac{3c}{4})$ and $Q(z_0,2c)$ for $N=1$.}
\end{figure}
\end{lemma}
To formulate similar result for the $\tilde{h}$-parabolic functions, let 
\[\tilde{Q}(\bar{z},c)=\Big \{w\in\mathbb{R}^{N+1}_-: \ |y+2\gamma\tau|<\tilde{R}_c(\tau), \ \tau_0-\frac{3c}{4}<\tau<\tau_0\Big \}, \]
be an $\tilde{h}$-ball $\tilde{B}(\bar{z},c)$ with the region below $\tau=\tau_0-\frac{3c}{4}$ removed. 
\begin{lemma}\label{harnacktildehparabolic}
Let $u\geq 0$ be an $\tilde{h}$-parabolic function in $\tilde{Q}(\bar{z},2c), c>0$, and suppose that $u$ is continuous at each point 
of $\partial \tilde{Q}(\bar{z},2c)$, except possibly at $\bar{z}$. Then there exists a positive constant $C$, depending only on $N$ and not on $c$, such that (Figure 11)
\begin{equation}\label{harnacktildehpar}
\dashint_{|y+2\gamma (\tau_0-\frac{3c}{2})|\leq \frac{(3Nc)^{\frac{1}{2}}}{2}}u\Big ( y,\tau_0-\frac{3c}{2}\Big )\,dx \leq C \inf\limits_{\tilde{B}(\bar{z},\frac{3c}{4})} u
\end{equation}
\end{lemma}
\begin{figure}[ht]\label{fig6}
    \centering
    \includegraphics[width=0.9\textwidth, trim=1cm 8cm 1cm 7cm, clip]{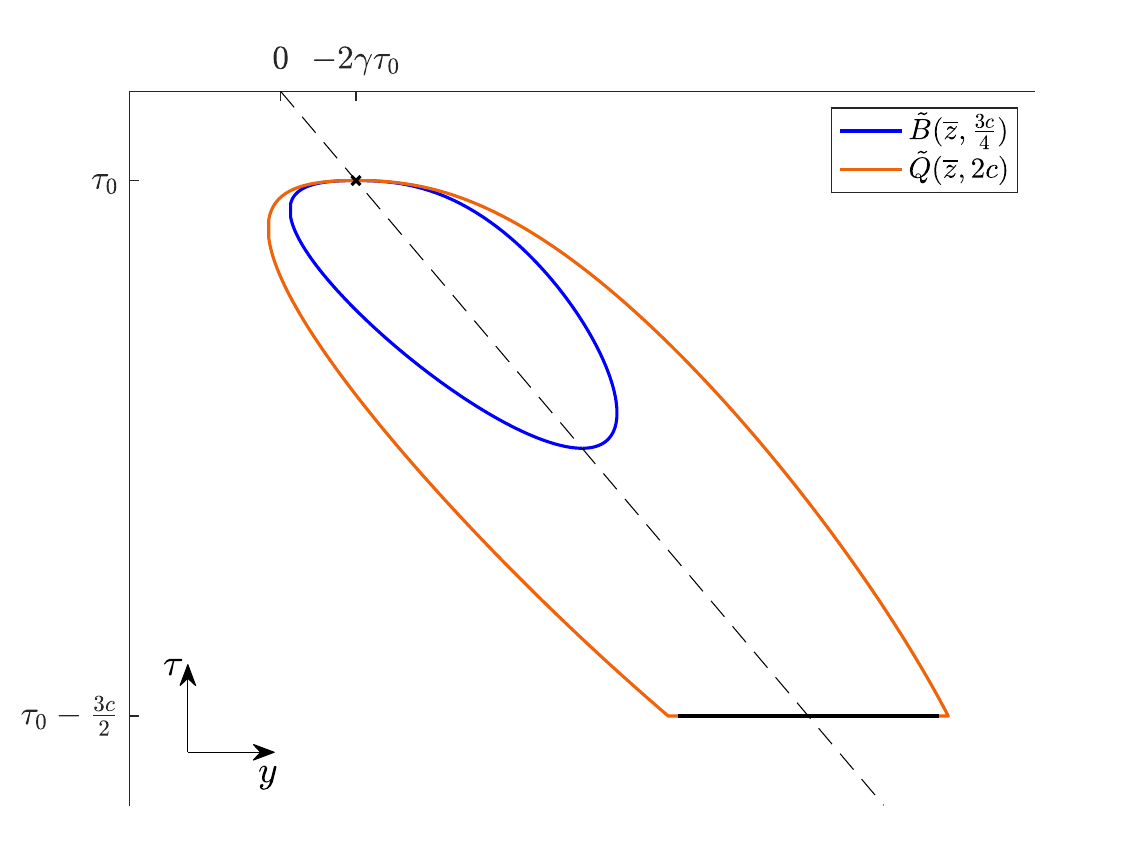}
    \caption{Sketch of $\tilde{h}$-heat ball $\tilde{B}(\bar{z},\frac{3c}{4})$ and $\tilde{Q}(\bar{z},2c)$ for $N=1$.}
\end{figure}
Note that in the particular case  $\gamma=0$ ($\tilde{h}\equiv 1$), Lemma \ref{harnacktildehparabolic} is known (\cite{evansgariepy}, Lemma 3.2). Using this fact, proof of Lemmas \ref{harnackhparabolic} and \ref{harnacktildehparabolic} can be pursued by applying the scheme \eqref{lemmasproofscheme2} and Lemma IV.1 of \cite{AbdullaI}. %%%%%%%%%%%%%%%%%%%

%%%%%%%%%%%%%%%%%%%%%%%%
\section{Dual $h_*$-Parabolic Potential Theory}
All the concepts and results of the potential theory generated by the heat operator $\mathcal{H}$ have a dual counterpart in the framework of the potential theory for the backward heat operator.
Smooth solutions of the equation $\mathcal{H}_*u=0$ will be called $*$-{\it parabolic} functions. A function $u$ is called $*$-{\it superparabolic} (or $*$-{\it subparabolic}) in open set $\Omega$ if $u(x,-t)$ is superparabolic (or subparabolic) in 
\[\hat{\Omega}:=\{(x,t): (x,-t)\in \Omega\}. \]

Let $\mathcal{P}_*(\Omega)$ be a class of $*$-{\it parabolic} functions
in an open set $\Omega$. Given open set $\Omega\subset\mathbb{R}^{N+1}_-$, the $*$-{\bf Appell transformation} is a homeomorphism ${\bf A}_*: \mathcal{P}_*(\Omega) \mapsto \mathcal{P}_*(A^{-1}\Omega)$ defined as
\begin{equation}\label{*Atransform}
\left\{
\begin{array}{l}
 \Omega\subset\mathbb{R}^{N+1}_-: u\in\mathcal{P}_*(\Omega) \mapsto {\bf A}_*u(z)=(\frac{\pi}{t})^{\frac{N}{2}} e^{\frac{|x|^2}{4t}}u(Az)\in \mathcal{P}_*(A^{-1}\Omega) \\ \\
 \Omega\subset\mathbb{R}^{N+1}_+: v\in\mathcal{P}_*(\Omega) \mapsto {\bf A}_*^{-1}v(z)=F(x,-t) v(A^{-1}(z))\in \mathcal{P}_*(A\Omega)
  \end{array}\right.
\end{equation}
The claim follows from the following formula:
\begin{equation}\label{*appelltransform}
\left\{
\begin{array}{l}
\mathcal{H}_*[{\bf A}_*u(z)]=\frac{\pi^{N/2}}{4} t^{-\frac{N}{2}-2}e^{\frac{|x|^2}{4t}} \mathcal{H}_*[u(A(z)], \ z\in A^{-1}\Omega\subset \mathbb{R}^{N+1}_+\\ \\
\mathcal{H}_*[{\bf A}_*^{-1}v(z)]=\frac{1}{4t^2} F(x,-t) \mathcal{H}_*[v(A^{-1}(z)], \ z\in A\Omega\subset \mathbb{R}^{N+1}_-.
\end{array}\right.
\end{equation}
We have
\begin{equation}\label{*afs1}
h_*={\bf A}_*\tilde{h}_*.
\end{equation}
As its counterpart,  the $*$-Appell transformation generates an homeomorphism ${\bf A}_*: \mathcal{S}_*(\Omega) \mapsto \mathcal{S}_*(A^{-1}\Omega)$, where $\mathcal{S}_*(\Omega)$ and  $\mathcal{S}_*(A^{-1}\Omega)$ denote the class of $*$-superparabolic functions in $\Omega\subset\mathbb{R}^{N+1}_-$ and $A^{-1}\Omega$ respectively. 

Next, we introduce the dual concept of $h$-parabolic function. 
A function $u=v/h_*$ is called an $h_*$-{\it parabolic}, $h_*$-{\it superparabolic}, or $h_*$-{\it subparabolic} in open subset $\Omega\subset \mathbb{R}^{N+1}_+$ if $v$ is $*$-{\it parabolic}, $*$-{\it superparabolic}, or $*$-{\it subparabolic} in $\Omega$. Note that for the sake of simplicity, we do a slight abuse of notation by writing $h_*$-{\it parabolic} instead of 
$h_*$-$*$-{\it parabolic}. Similarly, a function $u=v/\tilde{h}_*$ is called an $\tilde{h}_*$-{\it parabolic}, $\tilde{h}_*$-{\it superparabolic}, or $\tilde{h}_*$-{\it subparabolic} in open subset $\Omega\subset \mathbb{R}^{N+1}_-$ if $v$ is $*$-{\it parabolic}, $*$-{\it superparabolic}, or $*$-{\it subparabolic} in $\Omega$.

Next, we establish the dual of the formulae \eqref{HuA-1}, \eqref{HuA}. If $u:\mathbb{R}^{N+1}_-\to \mathbb{R}$ is smooth, then
\begin{equation}\label{*HuA-1}
\frac{\mathcal{H}_*[\tilde{h}_*(w)u(w)]}{\tilde{h}_*(z)}\Big |_{w=A(z)}=4t^2\frac{\mathcal{H}_*[h_*(z)u(A(z))]}{h_*(z)}.
\end{equation}
Equivalently,  if $u:\mathbb{R}^{N+1}_+\to \mathbb{R}$ is smooth, then
\begin{equation}\label{*HuA}
\frac{\mathcal{H}_*[h_*(z)u(z)]}{h_*(z)}\Big |_{z=A^{-1}(w)}=4\tau^2\frac{\mathcal{H}_*[\tilde{h}_*(w)u(A^{-1}(w))]}{\tilde{h}_*(w)}.
\end{equation}
{\it Proof of \eqref{*HuA-1} and \eqref{*HuA}}: Given a smooth function  $u:\mathbb{R}^{N+1}_-\to \mathbb{R}$, by \eqref{*Atransform} we have
\[ {\bf A}_*(u\tilde{h}_*)(z)=\Big(\frac{\pi}{t}\Big )^{\frac{N}{2}}e^{\frac{|x|^2}{4t}}u(A(z))\tilde{h}_*(A(z))=u(A(z))h_*(z). \]
From \eqref{*appelltransform} it follows that
\[ \mathcal{H}_*[{\bf A}_*(u\tilde{h}_*)](x,t)=\frac{\pi^{N/2}}{4}t^{-\frac{N}{2}-2}e^{\frac{|x|^2}{4t}} \mathcal{H}_*[u(y,\tau)\tilde{h}_*(y,\tau)]|_{(y,\tau)=A(x,t)}, \]
and hence
\[  \mathcal{H}_*[u(y,\tau)\tilde{h}_*(y,\tau)]|_{(y,\tau)=A(x,t)}=\frac{\mathcal{H}_*[h_*(x,t)u(A(x,t))]}{\frac{\pi^{N/2}}{4}t^{-\frac{N}{2}-2}e^{\frac{|x|^2}{4t}}}, \]
which implies \eqref{*HuA-1}. Similarly, given a smooth function  $u:\mathbb{R}^{N+1}_+\to \mathbb{R}$, by \eqref{*Atransform} we have
\[ {\bf A}_*^{-1}(h_*u)(w)=F(y,-\tau)h_*(A^{-1}(w))u(A^{-1}(w))=\tilde{h}_*(w)u(A^{-1}(w)). \]
From \eqref{*appelltransform} it follows that
\[ \mathcal{H}_*[{\bf A}_*^{-1}(h_*u)](w)=\frac{1}{4\tau^2}F(y,-\tau)\mathcal{H}_*[h_*(z)u(z)]|_{z=A^{-1}(w)}, \]
and hence
\[ \mathcal{H}_*[h_*(z)u(z)]|_{z=A^{-1}(w)}=4\tau^2\frac{\mathcal{H}_*[\tilde{h}_*(w)u(A^{-1}(w))]}{F(y,-\tau)}, \]
which implies \eqref{*HuA}.   \hfill{$\square$}

In particular, the function $u$ is $h_*$-parabolic (or $h_*$-super- or $h_*$-subparabolic) in $\Omega\subset\mathbb{R}^{N+1}_+$ if and only if $u(A^{-1}(w))$ is $\tilde{h}_*$-parabolic (or $\tilde{h}_*$-super- or $\tilde{h}_*$-subparabolic) in $A\Omega\subset \mathbb{R}^{N+1}_-$. 

The concepts $h_*$- and $\tilde{h}_*$-capacities are defined similarly as in Definitions 1.1 and 1.2:

\begin{definition}
For a compact set  $K\subset \mathbb{R}^{N+1}_+$, the $h_*$-capacity of $K$ is
\begin{equation}\label{*ChK}
 C_{h_*}(K)\equiv \max \{\mu(K): \ \mu\in \mathcal{M}_K, \ \hat{\mathbb{P}}^{h_*}_\mu\leq 1 \ \text{in} \ \mathbb{R}^{N+1}_+\},
 \end{equation}
 where
 \begin{equation}\label{*hpotential}
\hat{\mathbb{P}}^{h_*}_\mu(w)=\int\limits_{\mathbb{R}^{N+1}_+}\frac{F(z-w)}{h(z) h_*(w)}d\mu(z), \ w\in\mathbb{R}^{N+1}_+
\end{equation}
be an $h_*$-potential of measure $\mu\in\mathcal{M}_K$.
\end{definition}
For a compact set  $K\subset \mathbb{R}^{N+1}_-$, the $\tilde{h}_*$-capacity of $K$ is
\begin{equation}\label{*CtildehK}
 C_{\tilde{h}_*}(K)\equiv \max \{\mu(K): \ \mu\in \mathcal{M}_K, \ \hat{\mathbb{P}}^{\tilde{h}_*}_\mu\leq 1 \ \text{in} \ \mathbb{R}^{N+1}_-\},
 \end{equation}
 where
 \begin{equation}\label{*tildehpotential}
\hat{\mathbb{P}}^{\tilde{h}_*}_\mu(w)=\int\limits_{\mathbb{R}^{N+1}_-}\frac{F(z-w)}{\tilde{h}(z)\tilde{h}_*(w)}d\mu(z), \ w\in\mathbb{R}^{N+1}_-.
\end{equation}
be an $\tilde{h}_*$-potential of measure $\mu\in\mathcal{M}_K$.

Next, we formulate the dual of Lemma~\ref{hpotentialmeasure}, which demonstrates that the one-to-one mapping between $h_*$- and $\tilde{h}_*$-potentials, and their corresponding measures.
 \begin{lemma}\label{*hpotentialmeasure}
 (i) \ Given a compact set $K\subset\mathbb{R}^{N+1}_-$, $\hat{\mathbb{P}}^{\tilde{h}_*}_\mu$ is an $\tilde{h}_*$-potential of measure $\mu\in\mathcal{M}_K$ if and only if 
$\hat{\mathbb{P}}^{\tilde{h}_*}_\mu(Az)$ is an $h_*$-potential of measure $\mu(A)\in \mathcal{M}_{A^{-1}K}$, i.e.
\begin{equation}\label{*Atransformpotential}
\left\{
\begin{array}{l}
 K\subset\mathbb{R}^{N+1}_- \ \mu\in\mathcal{M}_K: \ \hat{\mathbb{P}}^{\tilde{h}_*}_\mu(Az)=\hat{\mathbb{P}}^{h_*}_{\mu(A)}(z), \ z\in\mathbb{R}^{N+1}_+ \\
 K\subset\mathbb{R}^{N+1}_+ \ \mu\in\mathcal{M}_K: \ \hat{\mathbb{P}}^{h_*}_\mu(A^{-1}w)=\hat{\mathbb{P}}^{\tilde{h}_*}_{\mu(A^{-1})}(w), \ w\in\mathbb{R}^{N+1}_-
  \end{array}\right.
\end{equation}
(ii) \ $\lambda_*\in\mathcal{M}_K$ is $\tilde{h}_*$-capacitary measure of compact $K\subset\mathbb{R}^{N+1}_-$ if and only if $\lambda_*(A)$ is $h_*$-capacitary measure of $A^{-1}K$. \\
(iii)  \ For a compact $K\subset\mathbb{R}^{N+1}_-$ with $\tilde{h}_*$-capacitary measure $\tilde{\lambda}_*\in \mathcal{M}_K$ we have
\begin{equation}\label{*capAK}
C_{h_*}(A^{-1}K)=\tilde{\lambda}_*(K)=C_{\tilde{h}_*}(K).
\end{equation}         
Conversely, for a compact $K\subset\mathbb{R}^{N+1}_+$ with $h_*$-capacitary measure $\lambda_*\in \mathcal{M}_K$ we have
\begin{equation}\label{*capA-1K}
C_{\tilde{h}_*}(AK)=\lambda_*(K)=C_{h_*}(K).
\end{equation}  
 \end{lemma}
Proof of Lemma~\ref{*hpotentialmeasure} is identical to the proof of Lemma~\ref{hpotentialmeasure}.

 \begin{lemma}\label{Ch=Ch*Ctildeh=Ctildeh*}
 (i) For arbitrary compact $K\in\mathbb{R}^{N+1}_+$ we have
 \begin{equation}\label{Ctildeh=Ctildeh*}
 C_h(K)=C_{h_*}(K),
 \end{equation}
 and for arbitrary compact $K\in\mathbb{R}^{N+1}_-$
  \begin{equation}\label{Ch=Ch*}
 C_{\tilde{h}}(K)=C_{\tilde{h}_*}(K).
 \end{equation}
 \end{lemma}
 {\it Proof.} Let $K\subset \mathbb{R}^{N+1}_+$ be a compact, and $\mu\in\mathcal{M}(K)$ and $\nu\in\mathcal{M}(K)$ are respectively it's $h$-capacitary, and $h_*$-capacitary measures, i.e.
 \[ C_h(K)=\mu(K), \ C_{h_*}(K)=\nu(K). \]
 Let $\{K_i\}$ be a monotone sequence of compact subsets with limit $K$, and such that $K$ is contained in the interior of each $K_i$. Let $\mu_i\in\mathcal{M}(K_i)$ and $\nu_i\in\mathcal{M}(K_i)$ are respectively $h$-capacitary, and $h_*$-capacitary measures of $K_i$, i.e.
 \[ C_h(K_i)=\mu_i(K), \ C_{h_*}(K_i)=\nu_i(K). \]
 Then we have
 \begin{eqnarray}
 C_h(K)=\mu(K)=\int_K1\,d\mu(w)=\int_{K_i} \hat{\mathbb{P}}^{h_*}_{\nu_i}(w)\,d\mu(w) =\nonumber\\ \int_{K_i} \mathbb{P}^{h}_{\mu}(z)\,d\nu_i(z)
 \leq \int_{K_i}1\,d\nu_i(z)=\nu_i(K_i)=C_{h_*}(K_i), \label{cap=cap*1}
 \end{eqnarray}
 and 
  \begin{eqnarray}
 C_{h_*}(K)=\nu(K)=\int_K1\,d\nu(z)=\int_{K_i} \mathbb{P}^{h}_{\mu_i}(z)\,d\nu(z) =\nonumber\\ \int_{K_i} \hat{\mathbb{P}}^{h_*}_{\nu}(w)\,d\mu_i(w)
 \leq \int_{K_i}1\,d\mu_i(w)=\mu_i(K_i)=C_{h}(K_i), \label{cap=cap*2}
 \end{eqnarray}
 By monotone convergenve property of capacities (Lemma~\ref{hpotentialvsreduction1} (iv) and Lemma~\ref{tildehpotentialvsreduction} (iv)) from \eqref{cap=cap*1} and \eqref{cap=cap*2} we deduce that
 \[ C_h(K)\leq C_{h_*}(K), \ C_{h_*}(K)\leq C_h(K)\]
 which proves \eqref{Ctildeh=Ctildeh*}. Similar proof applies to \eqref{Ch=Ch*}.  \hfill{$\square$}
 
 In the next lemma, we establishe a one-to-one mapping between the $*$-heat potentials
\[ \hat{P}^{h_*}_\mu(w)\equiv h_*\hat{\mathbb{P}}^{h_*}_\mu(w)=\int\limits_{\mathbb{R}^{N+1}_+}\frac{F(z-w)}{h(z)}\,d\mu(z); \ \hat{P}^{\tilde{h}_*}_\mu(w)\equiv \tilde{h}_*\hat{\mathbb{P}}^{\tilde{h}_*}_\mu(w)=\int\limits_{\mathbb{R}^{N+1}_-}\frac{F(z-w)}{\tilde{h}(z)}\,d\mu(z) \]
 \begin{lemma}\label{hpotentialmeasureolddual}
(i) \ For any fixed $z=(x,t)\in\mathbb{R}^{N+1}_-$, the $*$-Appell transform of $F(z-w), w=(y,\tau)\in \mathbb{R}^{N+1}_-$ is
\begin{equation}\label{A*F(z-w)}
{\bf A}_*F(z-w)=\frac{F(\bar{z}-w)}{F(\bar{z})}, \ w\in\mathbb{R}^{N+1}_+,
\end{equation}
where $\bar{z}=A^{-1}z\in \mathbb{R}^{N+1}_+$. \\
(ii) \ For any fixed $z=(x,t)\in\mathbb{R}^{N+1}_+$, the inverse $*$-Appell transform of $F(z-w), w=(y,\tau)\in \mathbb{R}^{N+1}_+$ is
\begin{equation}\label{A*-1F(z-w)}
{\bf A}_*^{-1}F(z-w)=\Big (-\frac{\tilde{\tau}}{\pi}\Big )^{\frac{N}{2}}e^{\frac{|\tilde{y}|^2}{4\tilde{\tau}}}F(\tilde{w}-w), \ w\in\mathbb{R}^{N+1}_-,
\end{equation}
where $\tilde{w}=Az\in \mathbb{R}^{N+1}_-$. \\
(iii) \ Given a compact set $K\subset\mathbb{R}^{N+1}_-$, $\hat{P}^{\tilde{h}_*}_\mu$ is an $*$-potential of measure $\frac{\mu}{\tilde{h}}\in\mathcal{M}_K$ if and only if 
${\bf A}_*\hat{P}^{\tilde{h}_*}_\mu$ is an $*$-potential of measure $\frac{\mu(A)}{h}\in \mathcal{M}_{A^{-1}K}$, i.e.
\begin{equation}\label{A*transformpotential1}
\left\{
\begin{array}{l}
 K\subset\mathbb{R}^{N+1}_- \ \mu\in\mathcal{M}_K: \ {\bf A}_*\hat{P}^{\tilde{h}_*}_\mu=\hat{P}^{h_*}_{\mu(A)} \\
 K\subset\mathbb{R}^{N+1}_+ \ \mu\in\mathcal{M}_K: \ {\bf A}_*^{-1}\hat{P}^{h_*}_\mu=\hat{P}^{\tilde{h}_*}_{\mu(A^{-1})} 
  \end{array}\right.
\end{equation}       
\end{lemma}
{\it Proof.} (i) \ We have
\begin{gather}
{\bf A}_*F(z-w)=\Big (\frac{\pi}{\tau}\Big )^{\frac{N}{2}} e^{\frac{|y|^2}{4\tau}}F\Big (A(\bar{z})-A(w)\Big )\nonumber\\
= \Big (\frac{\pi}{\tau}\Big )^{\frac{N}{2}} e^{\frac{|y|^2}{4\tau}} \Big (\frac{\pi(\bar{t}-\tau)}{\bar{t}\tau}\Big )^{-\frac{N}{2}} e^{-\frac{|\tau \bar{x}-\bar{t}y|^2}{4\bar{t}\tau(\bar{t}-\tau)}} 1_{\{\bar{t}>\tau\}}(\tau).\label{A*F1}
\end{gather}
Using the identity \eqref{ident},
from \eqref{A*F1}, \eqref{A*F(z-w)} follows. \\
(ii) \ We have
\begin{gather}
{\bf A}_*^{-1}F(z-w)=F(y,-\tau) F\Big (A^{-1}\tilde{w}-A^{-1}w)\Big )\nonumber\\
= (-4\pi \tau)^{-\frac{N}{2}}e^{\frac{|y|^2}{4\tau}} \Big (\frac{\pi(\tilde{\tau}-\tau)}{\tilde{\tau}\tau}\Big )^{-\frac{N}{2}} e^{-\frac{|\tau\tilde{y}-\tilde{\tau}y|^2}{4\tau\tilde{\tau}(\tilde{\tau}-\tau)}} 1_{\{\tau<\tilde{\tau}\}}(\tau).\label{A*-11}
\end{gather}
Using the identity \eqref{ident}, from \eqref{A*-11}, \eqref{A*-1F(z-w)} follows.\\
(iii) \  Given a compact $K\subset\mathbb{R}^{N+1}_-$ and a measure $\mu\in\mathcal{M}_K$, from \eqref{A*F(z-w)} it follows
\begin{equation}\label{A*Ph1}
{\bf A}_*\hat{P}^{\tilde{h}_*}_\mu(z)=\int\limits_K\frac{{\bf A}_*F(w-z)}{\tilde{h}(w)}d\mu(w)=\int\limits_{K}\frac{F(A^{-1}w-z)}{F(A^{-1}w)\tilde{h}(w)}\,d\mu(w), 
\end{equation}
Changing the variable $w \to A(w)$ from \eqref{A*Ph1} we deduce
\begin{equation}\label{A*Ph2}
{\bf A}_*\hat{P}^{\tilde{h}_*}_\mu(z)=\int\limits_{A^{-1}K}\frac{F(w-z)}{F(w)\tilde{h}(Aw)}\,d\mu(Aw)=\int\limits_{A^{-1}K}\frac{F(w-z)}{h(w)}\,d\mu(Aw)
\end{equation}
which implies the first statement of  \eqref{A*transformpotential1}.
On the contrary, given a compact $K\subset\mathbb{R}^{N+1}_+$ and a measure $\mu\in\mathcal{M}_K$, from \eqref{A*-1F(z-w)} it follows
\begin{equation}\label{A*-1Ph1}
 {\bf A}_*^{-1}\hat{P}^{h_*}_{\mu}(z)=\int\limits_K\frac{{\bf A}_*^{-1}F(w-z)}{h(w)}d\mu(w)=\int\limits_{K}\Big (-\frac{\tilde{\tau}}{\pi}\Big)^{\frac{N}{2}}e^{\frac{|\tilde{y}|^2}{4\tilde{\tau}}}\frac{F(\tilde{w}-z)}{h(w)}\,d\mu(w), 
\end{equation}
where $\tilde{w}=Aw$. Changing the variable $w \to A^{-1}(w)$ from \eqref{A*-1Ph1} we deduce
\begin{equation}\label{A*-1Ph2}
  {\bf A}_*^{-1}\hat{P}^{h_*}_{\mu}(z)=\int\limits_{AK} \Big (-\frac{\tau}{\pi}\Big)^{\frac{N}{2}}e^{\frac{|y|^2}{4\tau}} \frac{F(w-z)}{h(A^{-1}w)}\,d\mu(A^{-1}w)=\int\limits_{AK} \frac{F(w-z)}{\tilde{h}(w)}\,d\mu(A^{-1}w),
\end{equation}
which implies the second statement of  \eqref{A*transformpotential1}.  \hfill{$\square$}
 %%%%%%%%%%%%%%
 %%%%%%%%%%%%%%fix it to *-setting%%%%%%%
 
 Next lemmas are counterparts of Lemmas~\ref{hpotentialvsreduction1} and \ref{tildehpotentialvsreduction1} expressing the properties of $h_*$-capacity and $h_*$-capacitary potentials.
\begin{lemma}\label{h*potentialvsreduction1}
(i) \ For any compact $K\subset \mathbb{R}^{N+1}_+$ there is a unique $h_*$-capacitary measure $\lambda\in\mathcal{M}_{\partial K}$ such that
\begin{equation}\label{h*potentialvsreduction}
\ ^{h_*}\hat{R}^{K}_1 \equiv \hat{\mathbb{P}}^{h_*}_\lambda, \ \  \hat{R}^{K}_{h_*} \equiv \hat{P}^{h_*}_\lambda \ \ on \ \ \mathbb{R}^{N+1}_+ 
\end{equation}
where $\ ^{h_*}\hat{R}^{K}_1$ be a smoothed $h_*$-reduction of $1$ on $K$, $\hat{R}^{K}_{h_*}$ be smoothed $*$-reduction of $h_*$ on $K$, and $\hat{P}^{h_*}_\lambda$ be a $*$-heat potential with associated Riesz measure $\frac{\lambda}{h}$. \\
(ii) \ $\ ^{h_*}\hat{R}^{K}_1$ is $h_*$-superparabolic in $\mathbb{R}^{N+1}_+$, and $h_*$-parabolic in $\mathbb{R}^{N+1}_+\setminus \partial K$; \ $ \hat{R}^{K}_{h_*} \equiv P^{h_*}_\lambda$ is 
$*$-superparabolic in $\mathbb{R}^{N+1}_+$, and $*$-parabolic in $\mathbb{R}^{N+1}_+\setminus \partial K$.\\  
(iii) \ Furthermore,
\begin{equation}\label{h*potentialdistribution}
\mathcal{H}_*[\hat{P}^{h_*}_\lambda]\equiv\mathcal{H}_*[ h_* \hat{\mathbb{P}}^{h_*}_\lambda]=\frac{\lambda}{h} \ \text{in the sense of distributions on} \ \mathbb{R}^{N+1}_+ 
\end{equation}
(iv) \ As a set function defined in the class of compact subsets of $\mathbb{R}^{N+1}_+$, $C_{h_*}(\cdot)$ is a topological precapacity (see Lemma~\ref{hpotentialvsreduction1}, (iv)). 
As a topological precapacity, $C_{h_*}$ has an extension to countably strongly subadditive Choquet capacity $C_{h_*}(\cdot)$ relative to the class of compact subsets of $\mathbb{R}^{N+1}_+$ (see \eqref{extopen},\eqref{extarb}). The set function $C_{h_*}(\cdot): 2^{\mathbb{R}^{N+1}_+} \to [0,+\infty]$ is a countably strongly subadditive Choquet capacity
defined on all subsets of $\mathbb{R}^{N+1}_+$. All the Borel (and even analytic) subsets $E \subset \mathbb{R}^{N+1}_+$ (or $E \subset \mathbb{R}^{N+1}_-$) are $C_{h_*}$-capacitable (see Lemma~\ref{Choquet}).\\
(v) \ The set functions $C_h(\cdot)$ and $C_{h_*}(\cdot)$ are identical.
\end{lemma}   
Similar result in the $(\mathbb{R}^{N+1}_-,\tilde{h}_*)$ setting reads:
\begin{lemma}\label{tildeh*potentialvsreduction1}
(i) \ For any compact $K\subset \mathbb{R}^{N+1}_-$ there is a unique $\tilde{h}_*$-capacitary measure $\lambda\in\mathcal{M}_{\partial K}$ such that
\begin{equation}\label{tildeh*potentialvsreduction}
\ ^{\tilde{h}_*}\hat{R}^{K}_1 \equiv \hat{\mathbb{P}}^{\tilde{h}_*}_\lambda, \ \  \hat{R}^{K}_{\tilde{h}_*} \equiv \hat{P}^{\tilde{h}_*}_\lambda \ \ on \ \ \mathbb{R}^{N+1}_- 
\end{equation}
where $\ ^{\tilde{h}_*}\hat{R}^{K}_1$ be a smoothed $\tilde{h}_*$-reduction of $1$ on $K$, $\hat{R}^{K}_{\tilde{h}_*}$ be smoothed $*$-reduction of $\tilde{h}_*$ on $K$, and $\hat{P}^{\tilde{h}_*}_\lambda$ be a $*$-heat potential with associated Riesz measure $\frac{\lambda}{\tilde{h}}$. \\
(ii) \ $\ ^{\tilde{h}_*}\hat{R}^{K}_1$ is $\tilde{h}_*$-superparabolic in $\mathbb{R}^{N+1}_-$, and $\tilde{h}_*$-parabolic in $\mathbb{R}^{N+1}_-\setminus \partial K$; \ $ \hat{R}^{K}_{\tilde{h}_*} \equiv P^{\tilde{h}_*}_\lambda$ is 
$*$-superparabolic in $\mathbb{R}^{N+1}_-$, and $*$-parabolic in $\mathbb{R}^{N+1}_-\setminus \partial K$.\\  
(iii) \ Furthermore,
\begin{equation}\label{tildeh*potentialdistribution}
\mathcal{H}_*[\hat{P}^{\tilde{h}_*}_\lambda]\equiv\mathcal{H}_*[ \tilde{h}_* \hat{\mathbb{P}}^{\tilde{h}_*}_\lambda]=\frac{\lambda}{\tilde{h}} \ \text{in the sense of distributions on} \ \mathbb{R}^{N+1}_- 
\end{equation}
(iv) \ As a set function defined in the class of compact subsets of $\mathbb{R}^{N+1}_-$, $C_{\tilde{h}_*}(\cdot)$ is a topological precapacity (see Lemma~\ref{hpotentialvsreduction1}, (iv)). 
As a topological precapacity, $C_{\tilde{h}_*}$ has an extension to countably strongly subadditive Choquet capacity $C_{\tilde{h}_*}(\cdot)$ relative to the class of compact subsets of $\mathbb{R}^{N+1}_-$ (see \eqref{extopen},\eqref{extarb}). The set function $C_{\tilde{h}_*}(\cdot): 2^{\mathbb{R}^{N+1}_-} \to [0,+\infty]$ is a countably strongly subadditive Choquet capacity
defined on all subsets of $\mathbb{R}^{N+1}_-$. All the Borel (and even analytic) subsets $E \subset \mathbb{R}^{N+1}_-$ (or $E \subset \mathbb{R}^{N+1}_-$) are $C_{\tilde{h}_*}$-capacitable (see Lemma~\ref{Choquet}).\\
(v) \ The set functions $C_{\tilde{h}}(\cdot)$ and $C_{\tilde{h}_*}(\cdot)$ are identical. 
\end{lemma}   

Proof of items (i)-(iv) of Lemmas~\ref{h*potentialvsreduction1} and \ref{tildeh*potentialvsreduction1} is identical to the proofs of Lemmas~\ref{hpotentialvsreduction1} and \ref{tildehpotentialvsreduction1}. Claims of items (v) follow from Lemma~\ref{Ch=Ch*Ctildeh=Ctildeh*}.
 %%%%%%%%%%%%%%%%%%%%%%%%%%%%
 %%%%%%%%%%%%%%%%%%%%%%%%%%%%
 \section{Proofs of Theorems~\ref{wienertest} and \ref{wienertest2}}\label{prooftheorem11}
 Due to equivalence of problems $\mathcal{A}_\gamma$ we assume $\gamma=0$.
 
 \textbf{Proof of $h$-irregularity of $\mathcal{O}$}. Note that due to assertion (2) of Lemma IV.2 of \cite{AbdullaI}, without loss of generality, we can assume $\Omega$ is contained in $\{z: 0<t<1\}$. Recall also that by assumption \eqref{assonomega} we have $\mathcal{O}\in \Lambda(z_0,\Omega)$ for some $z_0\in\Omega$. 
 
 Assume that the series \eqref{wtest} converges, and prove that $u_*\not\equiv 0$. Recall that (formula (34) of \cite{AbdullaI}) we have
 \begin{equation}\label{u*H1infty}
 u_*(z)=h(z) \ ^hH_{1_{\{\mathcal{O}\}}}^\Omega(z), \ z\in\Omega.
 \end{equation}
 Therefore, it is enough to prove that $^hH_{1_{\{\mathcal{O}\}}}^\Omega\not\equiv 0$. In fact, we are going to prove that
 \begin{equation}\label{Oirregular}
 \limsup\limits_{z\to \mathcal{O}} \ ^hH_{1_{\{\mathcal{O}\}}}^\Omega(z)=1.
 \end{equation}
 For arbitrary small $0<\epsilon < 1/2$, choose $m$ so large that
 \begin{equation}\label{wser<infty}
 \sum\limits_{n=m+1}^\infty 2^{-\frac{nN}{2}}C_h(E_n) < (2\pi)^{\frac{N}{2}} \epsilon,
 \end{equation}
and $\Lambda(z_0,\Omega)\cap B(\bar{z},2^m)\neq \emptyset$.  Let $\Omega_1$ be the interior of $\Omega \cup \overline{B(\bar{z},2^m)}$. Since outside  $B(\bar{z},2^m)$, $\Omega$ and $\Omega_1$ coincide, clearly the condition \eqref{wser<infty} is satisfied for $\Omega_1$. Also note that $\Lambda(z_0,\Omega)\subset \Lambda(\bar{z},\Omega_1)$, and therefore, $\mathcal{O}\in \partial\Lambda(\bar{z},\Omega_1)$. 

Consider a sequence of monotonically increasing $h$-parabolic functions
\[ G_M(z)=\sum\limits_{n=m+1}^M \mathbb{P}_{\mu_n}^h(z), \ M=m+1,m+2,... \ x\in \Big ( \cup_{n=m+1}^{+\infty}E_n\Big )^c, \]
where $\mathbb{P}_{\mu_n}^h$ be an $h$-capacitary potential of $E_n$, and $\mu_n$ be a corresponding $h$-capacitary measure. Clearly, we have
\[ G(z)=\sum\limits_{n=m+1}^{+\infty} \mathbb{P}_{\mu_n}^h(z)\equiv \sup\limits_M G_M(z), \ z\in \Omega_1. \]
According to Harnack convergence theorem \cite{Doob} in each connected component $\Lambda(z,\Omega_1)$, either $G\equiv +\infty$ or $G$ is an $h$-parabolic function.  
By using assertion (i) of Lemma~\ref{hpotentialvsreduction1}, we deduce
\begin{equation}\label{Gzbar}
G(\bar{z})=\sum\limits_{n=m+1}^{+\infty}\int_{E_n}\frac{F(\bar{z}-w)}{h(\bar{z})h_*(w)}d\mu_n(w) \leq (2\pi)^{-N/2} \sum\limits_{n=m+1}^\infty 2^{-\frac{nN}{2}}C_h(E_n)<\epsilon.
\end{equation}
Hence, $G$ is $h$-parabolic in $\Lambda(\bar{z},\Omega_1)$, and it is a generalized solution of the $h$-parabolic Dirichlet problem in $\Lambda(\bar{z},\Omega_1)$ under the boundary condition
\begin{equation}\label{Gboundary}
\sum\limits_{n=m+1}^{+\infty} \mathbb{P}_{\mu_n}^h(z)\geq 1 \ on \ \partial\Lambda(\bar{z},\Omega_1)\cap \{t<1\}.
\end{equation}
Note that the boundary function vanishes at the boundary point $\mathcal{O}$. Since the generalized solution is order preserving it follows that
\begin{equation}\label{orderpres}
0\leq \ ^hH_{1-1_{\{\mathcal{O}\}}}^{\Lambda(\bar{z},\Omega_1)}\leq G \ \quad\text{on} \ \Lambda(\bar{z},\Omega_1).
\end{equation}
From \eqref{Gzbar} and \eqref{Gboundary} it follows that
\begin{equation}\label{liminfG>}
0\leq \liminf_{z\to\mathcal{O}} G \leq \epsilon.
\end{equation}
Accordingly, by \eqref{orderpres} we have
\begin{equation}\label{liminfhPDP}
0\leq \liminf_{z\to\mathcal{O}}  \ ^hH_{1-1_{\{\mathcal{O}\}}}^{\Lambda(\bar{z},\Omega_1)} \leq \epsilon.
\end{equation}
Since
\begin{equation}
 \ ^hH_{1-1_{\{\mathcal{O}\}}}^{\Lambda(\bar{z},\Omega_1)} =1- \ ^hH_{1_{\{\mathcal{O}\}}}^{\Lambda(\bar{z},\Omega_1)},
\end{equation}
from \eqref{liminfhPDP} it follows that
\begin{equation}\label{limsuphPDP}
1-\epsilon\leq \limsup_{z\to\mathcal{O}}  \ ^hH_{1_{\{\mathcal{O}\}}}^{\Lambda(\bar{z},\Omega_1)} \leq 1.
\end{equation}
Hence, we have $ \ ^hH_{1_{\{\mathcal{O}\}}}^{\Lambda(\bar{z},\Omega_1)}\not\equiv 0$, and by the monotonicity of $h$-parabolic measure with respect to expanding open sets (Lemma IV.2 statement 1(b) of \cite{AbdullaI}) it follows that $ \ ^hH_{1_{\{\mathcal{O}\}}}^{\Omega_1}\not\equiv 0$. Since outside  $B(\bar{z},2^m)$, $\Omega$ and $\Omega_1$ coincide, we deduce that $ \ ^hH_{1_{\{\mathcal{O}\}}}^{\Omega}\not\equiv 0$ (Lemma IV.2 statement 2 of \cite{AbdullaI}). Hence, $\mathcal{O}$ is $h$-irregular boundary point for $\Omega$, and singularity at $\mathcal{O}$ is not removable. We are still going to prove the relation \eqref{Oirregular}. From the construction of the open set $\Omega_1$, and from the definition of the $h$-parabolic measure of $\mathcal{O}$ it follows that
\begin{eqnarray}
 \ ^hH_{1_{\{\mathcal{O}\}}}^{\Lambda(\bar{z},\Omega_1)} \Big|_{\{t<\frac{1}{1+4\cdot 2^m} \}}\equiv \ ^hH_{1_{\{\mathcal{O}\}}}^{\Lambda(\bar{z},\Omega_1)\cap\{t<\frac{1}{1+4\cdot 2^m}\}}
 \equiv  \ ^hH_{1_{\{\mathcal{O}\}}}^{\Lambda(z_1,\Omega)\cap\{t<\frac{1}{1+4\cdot 2^m}\}}\nonumber\\ \equiv  \ ^hH_{1_{\{\mathcal{O}\}}}^{\Lambda(z_1,\Omega)} \Big|_{\{t<\frac{1}{1+4\cdot 2^m} \}}, \ \quad\text{for some} \ z_1\in\Omega. \label{limsupatO>1-eps}
 \end{eqnarray} 
 From \eqref{limsuphPDP},\eqref{limsupatO>1-eps} it follows that
 \begin{equation}\label{limsuphPDP1}
1-\epsilon\leq \limsup_{z\to\mathcal{O}}  \ ^hH_{1_{\{\mathcal{O}\}}}^{\Lambda(z_1,\Omega)} \leq 1,
\end{equation}
or
\begin{equation}\label{limsuphPDP2}
1-\epsilon\leq \limsup_{z\to\mathcal{O}}  \ ^hH_{1_{\{\mathcal{O}\}}}^{\Omega} \leq 1.
\end{equation}
Taking $\epsilon\to 0$, \eqref{Oirregular} follows. 

 \textbf{Proof of $h$-regularity of $\mathcal{O}$}. As before, without loss of generality we can assume $\Omega$ is contained in $\{z: 0<t<1\}$ (Lemma IV.2 assertion (2) \cite{AbdullaI}). Assume that the series \eqref{wtest} is divergent and prove that $u_*\equiv 0$. Due to \eqref{u*H1infty} it is enough to prove that $^hH_{1_{\{\mathcal{O}\}}}^\Omega\equiv 0$. The proof consists of five steps.
 
 \textbf{Step 1: Modification of $\Omega$}. The goal of this step is to replace $\Omega$ with a larger version that still preserves the divergence of the series \eqref{wtest}. The latter implies that at least one of the following four series must be divergent:
 \[ \sum\limits_l 2^{-\frac{N(4l+j)}{2}}C_h(E_{4l+j}), \ j=0,1,2,3.  \]
 Without loss of generality, we can assume that the series with $j=0$ is divergent, i.e. 
 \[  \sum\limits_{l=1}^{+\infty} 2^{-2Nl}C_h(E_{4l})=+\infty. \]
Since $h$-capacity is subadditive set function, for arbitrary $\epsilon\in (0,1)$ we have
\[ C_h(\Omega^c\cap D(2^{4l})\leq C_h(\Omega^c\cap D(2^{4l})\cap\{t\leq \epsilon\})+C_h(\Omega^c\cap D(2^{4l})\cap\{\epsilon\leq t \leq 1\})\]
and
\[ \lim\limits_{\epsilon\to 1}  C_h(\Omega^c\cap D(2^{4l})\cap\{\epsilon\leq t \leq 1\})=0, \]
it follows that for some sequence $\{\epsilon_l\}$ such that $0<\epsilon_l\uparrow 1$ we have
\begin{equation}\label{modifiedwtest}
\sum\limits_{l=1}^{+\infty}2^{-2Nl}C_h(B_l)=+\infty,
\end{equation}
where $B_l=\Omega^c\cap D(2^{4l})\cap\{t\leq \epsilon_l\})$. We define a new open set $\Omega^{'}$ such that
\begin{equation}\label{sumBl}
\Omega^{'c}\equiv \{\mathcal{O}\}\cup \sum\limits_{l=1}^{+\infty}B_l.
\end{equation}
Then we have $\Omega\subset\Omega^{'}$, $\mathcal{O}\in\partial\Omega^{'}$. We are going to prove that $^hH_{1_{\{\mathcal{O}\}}}^{\Omega^{'}} \equiv 0$. 
This will imply that $^hH_{1_{\{\mathcal{O}\}}}^{\Omega} \equiv 0$. Furthermore, without loss of generality, we assume that $\Omega= \Omega^{'}$, and its complement satisfies \eqref{sumBl}.

 \textbf{Step 2: Approximation}. Let $v=\ ^{h}\hat{R}^{\Omega^c}_1$ be a smoothed $h$-reduction of $1$ on $\Omega^c$. It is $h$-superharmonic function in $\mathbb{R}^{N+1}_+$, and
 \[ v= \ ^hH_{1-1_{\{\mathcal{O}\}}}^{\Omega} \ \quad\text{in} \ \Omega; \  \ v\equiv 1 \ \quad\text{in} \ int  \ \Omega^c. \]
 It's complementary function $u=1-v$ is $h$-subparabolic in $\mathbb{R}^{N+1}_+$, and 
 \begin{gather}\label{uhparabolicmeasure}
 \left\{
\begin{array}{l}
 u=1- \ ^hH_{1-1_{\{\mathcal{O}\}}}^{\Omega} = \ ^hH_{1_{\{\mathcal{O}\}}}^{\Omega} \quad\text{in} \ \Omega; \\
   \ u\equiv 0 \ \quad\text{in} \ int  \ \Omega^c; \ 0\leq u \leq 1 \ \quad\text{in} \ \mathbb{R}^{N+1}_+.
   \end{array} \right.
 \end{gather}
 Our goal is to prove that $u\equiv 0$. 
 
 Take an arbitrary $\epsilon>0$, and consider an open set $\Omega_\epsilon$ such that $\Omega_\epsilon^c=\sum\limits_{l=1}^{+\infty} B_l^\epsilon$, where
 \begin{equation}\label{Blepsilon}
\left\{
\begin{array}{l}
 B_l^\epsilon \ \quad\text{is a compact set in} \ \mathbb{R}^{N+1}_+;  \ B_l \subset int B_l^\epsilon; \\
 B_l^\epsilon \subseteq B(2^{4l+1})\setminus \overline{B(2^{4l-2})}.
  \end{array}\right.
\end{equation} 
Let $v^\epsilon=\ ^{h}\hat{R}^{\Omega_\epsilon^c}_1$ be a smoothed $h$-reduction of $1$ on $\Omega_\epsilon^c$, and $u^\epsilon=1-v^\epsilon$. It can be proved that 
\begin{equation}\label{uepsred1}
u^\epsilon \uparrow \tilde{u}=1-\ ^{h}R^{\Omega^c}_1 \ \quad\text{as} \ \epsilon \downarrow 0,
\end{equation}
where $\ ^{h}R^{\Omega^c}_1$ be an $h$-reduction of $1$ on $\Omega^c$ (without smoothing). 

Let $\hat{\zeta}:=\hat{\mathbb{P}}^{h_*}_{\bar{\mu}_l}$ be $h_*$-capacitary potential for $B_l$, and $\hat{\mu}=\hat{\mu}_{B_l}$ be it's $h_*$-capacitary measure. By Lemma~\ref{h*potentialvsreduction1}, (iii) we have
\begin{gather}
\mathcal{H}_*[ h_* \hat{\zeta}]=\frac{\hat{\mu}}{h} \ \text{in the sense of distributions on} \ \mathbb{R}^{N+1}_+, \label{*cappot}\\
\hat{\mu}(B_l)=C_{h_*}(B_l)=C_h(B_l). \label{*cap=cap}
\end{gather}
Next, we approximate $\hat{\mu}$ and $u_\epsilon$ by smooth functions. Let $\rho\in C^\infty( \mathbb{R}^{N+1}; [0,1])$ be a standard mollifier with compact support, and
\[ \int_{\mathbb{R}^{N+1}}\rho\,dx\,dt=1. \]
Consider convolutions
\[ u^\epsilon_\delta=\rho_\delta * u^\epsilon, \ \hat{\zeta}_\delta=\hat{\zeta} * \rho_\delta, \ \delta>0.\]
Functions $u^\epsilon_\delta$ and $\hat{\zeta}_\delta$ are smooth and $0\leq u^\epsilon_\delta, \hat{\zeta}_\delta \leq 1$, and
\begin{equation}\label{husubparabolic}
\mathcal{H}[h u^\epsilon_\delta] \leq 0.
\end{equation}
Define the measure
\[ \hat{\mu}_\delta(B)=\int_B \mathcal{H}_*[h_*\hat{\zeta}_\delta]h\,dx\,dt, \ B\in \mathcal{B}, \]
where $\mathcal{B}$ be a Borel $\sigma$-algebra of the subsets of $\mathbb{R}^{N+1}_+$. The sequence of measures $\hat{\mu}_\delta$ converges to $\hat{\mu}$ weakly as $\delta\to 0$. 
Let us fix
\begin{equation}\label{c=24l}
c=2^{4l+1}
\end{equation}
so that we have
\[ B_l\subset D\Big (\frac{c}{2}\Big ), \  B^\epsilon_l\subset B(c)\setminus \overline{B}\Big ( \frac{c}{8}\Big ).  \]
We can also arrange that
\begin{equation}\label{suptmu}
supt \ \hat{\mu}_\delta \subset int \ B^\epsilon_l \cap B\Big ( \frac{3}{4}c\Big ),
\end{equation}
if $\delta$ is small enough, depending on $\epsilon$. Recalling that $u^\epsilon=0$ in $int \ B^\epsilon_l$, from \eqref{suptmu} it follows that
\begin{equation}\label{u=0sptmu}
u^\epsilon_\delta=0 \ \quad\text{on} \ supt \ \hat{\mu}_\delta,
\end{equation}
if $\delta$ is sufficiently small. Consider a region
\[ Q(c)=\Big \{ (x,t)\in\mathbb{R}^{N+1}_+: |x|^2<2Nt(1-t)\log \frac{4ct}{1-t}, \frac{1}{1+3c}<t<1\Big \}, \]
which is the part of the $h$-heat ball $B(c)$ lying above the plane $t=1/(1+3c)$. We then write
\begin{equation}\label{split}
u^\epsilon_\delta=v^\epsilon_\delta-w^\epsilon_\delta, 
\end{equation}
where
\begin{equation}\label{vepsilondelta}
\left\{
\begin{array}{l}
\mathcal{H}[hv^\epsilon_\delta]=0, \ \quad\text{in} \ Q(2c) \\ \\
v^\epsilon_\delta=u^\epsilon_\delta \ \quad\text{in} \ \partial Q(2c)\setminus\{\bar{z}\},
  \end{array}\right.
\end{equation}
and
\begin{equation}\label{wepsilondelta}
\left\{
\begin{array}{l}
\mathcal{H}[hw^\epsilon_\delta]=-\mathcal{H}[hu^\epsilon_\delta]\geq 0, \ \quad\text{in} \ Q(2c) \\ \\
w^\epsilon_\delta=0 \ \quad\text{in} \ \partial Q(2c)\setminus\{\bar{z}\}.
  \end{array}\right.
\end{equation}

 \textbf{Step 3: Estimates}. The goal of this and the following steps is to build a connection between the divergence of the series \eqref{modifiedwtest}, and the smallness of average spatial integrals of $u$ on the time levels approaching the origin. We now estimate
 \begin{gather}
\Big( \inf_{B(3c/4)} v^\epsilon_\delta \Big )^2 \hat{\mu}_\delta\Big (B\Big (\frac{3c}{4}\Big )\Big ) \leq \int\limits_{B(3c/4)} (v^\epsilon_\delta)^2\,d\hat{\mu}_\delta \leq \int\limits_{Q(2c)} (v^\epsilon_\delta)^2\,d\hat{\mu}_\delta \nonumber\\
=\int\limits_{Q(2c)} (w^\epsilon_\delta)^2\,d\hat{\mu}_\delta \ \  \quad\text{(due to} \ \eqref{u=0sptmu},  \eqref{split}) \nonumber\\
=\int\limits_{Q(2c)} (w^\epsilon_\delta)^2 \mathcal{H}_*[h_*\hat{\zeta}_\delta]h\,dx\,dt \nonumber\\
=\int\limits_{Q(2c)} \mathcal{H}[h(w^\epsilon_\delta)^2] h_*\hat{\zeta}_\delta\,dx\,dt \ \quad\text{(since} \ w^\epsilon_\delta=0 \quad\text{on} \ \partial Q(2c)\setminus\{\bar{z}\}) \nonumber \\
=2 \int\limits_{Q(2c)} w^\epsilon_\delta \mathcal{H}[hw^\epsilon_\delta] h_*\hat{\zeta}_\delta\,dx\,dt-2\int\limits_{Q(2c)}h|\nabla w^\epsilon_\delta|^2h_*\hat{\zeta}_\delta\,dx\,dt \nonumber\\
\leq 2 \int\limits_{Q(2c)}w^\epsilon_\delta\frac{-\mathcal{H}[hu^\epsilon_\delta]}{(2t)^Nh}\,dx\,dt  \ \  \quad\text{(due to \ \eqref{wepsilondelta}, and since} \ hh_*=(2t)^{-N}) \nonumber\\
\leq C\int\limits_{Q(2c)}\frac{-\mathcal{H}[hu^\epsilon_\delta]}{(2t)^Nh}\,dx\,dt \leq C\int\limits_{B(2c)}\frac{-\mathcal{H}[hu^\epsilon_\delta]}{(2t)^Nh}\,dx\,dt. \nonumber
 \end{gather}
Hence, for sufficiently small $\delta>0$ we have
\begin{equation}\label{key}
\Big( \inf_{B(3c/4)} v^\epsilon_\delta \Big )^2 \hat{\mu}_\delta\Big (B\Big (\frac{3c}{4}\Big )\Big ) \leq C\int\limits_{B(2c)}\frac{-\mathcal{H}[hu^\epsilon_\delta]}{(2t)^Nh}\,dx\,dt.
\end{equation}
Applying the statement (ii) of Lemma~\ref{avprophparabolic} to the right-hand side, and the Harnack inequality \eqref{harnackhpar} of Lemma~\ref{harnackhparabolic} to the left-hand side of \eqref{key} we deduce
\begin{equation}\label{harnackhparproof}
\Big [\dashint_{|x|\leq \frac{(3Nc)^{\frac{1}{2}}}{1+6c}}u^\epsilon_\delta\Big ( x,\frac{1}{1+6c}\Big )\,dx\Big ]^2 \hat{\mu}_\delta\Big (B\Big (\frac{3c}{4}\Big )\Big ) c^{-\frac{N}{2}} \leq C [\phi(u^\epsilon_\delta, 8c)-\phi(u^\epsilon_\delta, 4c)],
\end{equation}
where
\[ \phi(g,c)=\frac{1}{(4c)^{\frac{N}{2}}}\int\limits_{B(c)}\frac{g(x,t)|x|^2}{t^{N+2}|t-1|^2}\,dx\,dt. \] 

 \textbf{Step 4: Passage to limits}. Passing to limit as $\delta\downarrow 0$, we have
 \[ u^\epsilon_\delta \to u^\epsilon, \quad\text{a.e. on} \ \mathbb{R}^{N+1}_+; \ \quad\text{and uniformly on compact subsets of} \ \{\bar{z}\}\cup\sum\limits_{l=1}^\infty B^\epsilon_l.\] 
 Due to weak convergence of $\hat{\mu}_\delta$ to $\hat{\mu}$ from \eqref{suptmu},\eqref{*cap=cap} it follows that
 \[ \hat{\mu}_\delta\Big (B\Big (\frac{3c}{4}\Big )\Big )=\hat{\mu}_\delta(\mathbb{R}^{N+1}_+)\to \hat{\mu}(\mathbb{R}^{N+1}_+)=C_h(B_l). \]
 Therefore, from \eqref{harnackhparproof} we derive the estimate
 \[ \Big [ \dashint_{|x|\leq \frac{(3Nc)^{\frac{1}{2}}}{1+6c}}u^\epsilon\Big ( x,\frac{1}{1+6c}\Big )\,dx\Big ]^2 \ C_h(B_l) c^{-\frac{N}{2}} \leq C [\phi(u^\epsilon, 8c)-\phi(u^\epsilon, 4c)]. \]
 Passing to limit as $\epsilon\downarrow 0$, and using \eqref{uepsred1} and the fact that $u=\tilde{u}$ in $\Omega$ we derive the estimate
 \begin{equation}\label{prewietest}
 \Big [ \dashint_{|x|\leq \frac{(3Nc)^{\frac{1}{2}}}{1+6c}}u\Big ( x,\frac{1}{1+6c}\Big )\,dx\Big ]^2 \ C_h(B_l) c^{-\frac{N}{2}} \leq C [\phi(\tilde{u}, 8c)-\phi(\tilde{u}, 4c)].
 \end{equation}
 From \eqref{husubparabolic}, and Lemma~\ref{avprophparabolic}, (ii) it follows that
 \[ \phi(u^\epsilon_\delta, 2^kc)-\phi(u^\epsilon_\delta, 2^{k-1}c) \geq 0, k=0,1,2. \]
 Passing to the limit,first as $\delta\downarrow 0$, and then as $\epsilon\downarrow 0$ we derive that
 \[ \phi(\tilde{u}, 2^kc)-\phi(\tilde{u}, 2^{k-1}c) \geq 0, k=0,1,2, \]
 and therefore, from \eqref{prewietest} it follows the estimate
  \begin{equation}\label{prewietest1}
 \Big [ \dashint_{|x|\leq \frac{(3Nc)^{\frac{1}{2}}}{1+6c}}u\Big ( x,\frac{1}{1+6c}\Big )\,dx\Big ]^2 \ C_h(B_l) c^{-\frac{N}{2}} \leq C \sum\limits_{k=0}^3 [\phi(\tilde{u}, 2^kc)-\phi(\tilde{u}, 2^{k-1}c)].
 \end{equation}
Now, we substitute $c=2^{4l+1}$, and pursue summation for $l\geq n$. Remarkably,  the right-hand side produces a telescoping series, and we deduce that
\begin{equation}\label{wtestvsu}
\sum\limits_{l=n}^{+\infty} \chi_l^2 2^{-2nl} C_h(B_l) <+\infty,
\end{equation}
where
\begin{equation}\label{chiltozero}
\chi_l \equiv \dashint_{|x|\leq \frac{(6N)^{\frac{1}{2}}2^{2l}}{1+12\cdot 2^{4l}}}u\Big ( x,\frac{1}{1+12\cdot 2^{4l}}\Big )\,dx, \ l=n,n+1,...
\end{equation}
From the divergence of the series \eqref{modifiedwtest} it follows that
\begin{equation}\label{liminf0}
\liminf_{l \to+\infty} \chi_l =0.
\end{equation}
 \textbf{Step 5: $h$-parabolic measure of $\{\mathcal{O}\}$ vanishes.} Finally, in this step, we can deduce from the condition \eqref{liminf0} that $u\equiv 0$. 
 For the reader's convenience, we divide this step into three sub-steps. 
 
 \textit{Sub-step 5.1: Transformation from $h$-parabolic to parabolic setting. }
 Let $T$ be a translation map
 \[ T(y,\tau)=(y, \tau+1/4): \mathbb{R}^{N+1}_-\cap \{\tau<-1/4\} \to \mathbb{R}^{N+1}_-. \]
 Consider a homeomorphism $T\circ A: \mathbb{R}^{N+1}_+\cup\{t<1\} \mapsto \mathbb{R}^{N+1}_-$ with
\begin{equation*}%\label{tappeltrans}
\left\{
\begin{array}{l}
(x,t)\in\mathbb{R}^{N+1}_+\cap \{t<1\} \mapsto T\circ A(x,t)=\Big (\frac{x}{2t}, -\frac{1}{4t}+\frac{1}{4}\Big )\in \mathbb{R}^{N+1}_-;\\
(y,\tau)\in\mathbb{R}^{N+1}_- \mapsto (T\circ A)^{-1}(y,\tau)\equiv A^{-1}(y,\tau-1/4) \in \mathbb{R}^{N+1}_+\cap\{t<1\}; 
\end{array}\right.
\end{equation*}
Let $\tilde{\Omega}=T\circ A\Omega$. We have
\[ \tilde{\Omega}^c\cap\mathbb{R}^{N+1}_-=\sum_{l=1}^{+\infty}\tilde{B}_l, \tilde{B}_l\subset\tilde{B}(0,2^{4l+1})\setminus \overline{\tilde{B}(0,2^{4l-2})}, \]
 where $\tilde{B}(0,c)$ be an heat ball of center $(0,0)$ and radius $c$ defined in \eqref{heatball}. 
 Consider a transformed function
 \[ \theta(y,\tau)=u(A^{-1}(y,\tau-1/4)), \ (y,\tau)\in \mathbb{R}^{N+1}_-. \]
 The map $A^{-1}(y,\tau-1/4)): \mathbb{R}^{N+1}_- \to \mathbb{R}^{N+1}_+\cap \{t<1\}$ is a composition of a translation map $T^{-1}(y,\tau)=(y, \tau-1/4)$, and a map $A^{-1}(y,\tau)$. 
 Since $u$ is $h$-subparabolic, the function $u(A^{-1}(y,\tau)$ is subparabolic in $\mathbb{R}^{N+1}_-\cap \{\tau<-1/4\}$ (Lemma IV.1, \cite{AbdullaI}). Since subparabolicity is invariant under the parallel translation, $\theta$ is subparabolic in $\mathbb{R}^{N+1}_-$. Precisely, from the Lemma IV.1, \cite{AbdullaI}) and \eqref{uhparabolicmeasure} it follows that
 \begin{equation}\label{thetaparabolicmeasure}
 \theta= H_{1_{\{\mathcal{\infty}\}}}^{\tilde{\Omega}} \quad\text{in} \ \tilde{\Omega}; 
   \ \theta\equiv 0 \ \quad\text{in} \ int  \ \tilde{\Omega}^c; \ 0\leq \theta \leq 1 \ \quad\text{in} \ \mathbb{R}^{N+1}_-.
 \end{equation} 
 From \eqref{chiltozero} it follows that
 \begin{equation}\label{chiltozerotheta}
 \chi_l \equiv \dashint_{|y|\leq (\frac{3N}{2})^{1/2}\cdot 2^{2l}}\theta( y, -3\cdot 2^{4l})\,dx, \ l=n,n+1,...
 \end{equation}
Note that due to the procedure of enlargement of $\Omega$, the integration regions in \eqref{chiltozero} are contained in $\Omega$, and therefore $u$ is a smooth function in all those sets. The property is preserved under the transformation $T\circ A$, and the integration region in \eqref{chiltozerotheta} is contained in $\tilde{\Omega}$, and therefore, $\theta$ is a smooth function on it. 
Hence, to complete the proof, we need to demonstrate that \eqref{liminf0} implies that $\theta\equiv 0$. The proof of this claim is similar to the proof given in \cite{evansgariepy} and \cite{Abdulla2}(pages 25,26). The idea of the proof is to exploit the parabolic scaling to demonstrate that the sequence of maxima of the subparabolic function $\theta$ in an expanding sequence of cylinders covering $\mathbb{R}^{N+1}_-$ will form a geometric sequence on the expense of the condition \eqref{liminf0}. 
 
 Let
 \[ S_l\equiv \Big\{z: |x|\leq \Big(\frac{3N}{2}\Big)^{1/2}\cdot 2^{2l}, \ -3\cdot 2^{4l}\leq t \leq 0\Big\} \]
 and denote
 \[ M_l \equiv \sup_{S_l} \theta, \ l= -1,0,1,2,... \]
 We will prove that $M_l \equiv 0, l=-1,0,1,2,...$. 
 
 \textit{Sub-step 5.2: Supremum norm estimate in interior sub-cylinder using Green representation for parabolic functions.}  Let $G$ be a Green function for the heat operator in a cylinder $S_0$. Assume that $w$ solves $\mathcal{H}w=0$ in $S_0$, $w$ is continuous on the parabolic boundary $\partial_pS_0$, and $0\leq w \leq M$ in $S_0$. By the Green's representation formula, we have
 \begin{gather}
 w(x,t)=\int_{|y|^2\leq 3N/2}G(x,t; y, -3)w(y,-3)\,dy  \nonumber\\
- \int_{-3}^0\int_{|y|^2=3N/2} \frac{\partial G(x,t; y,s)}{\partial n} w(y,s)\,dH^{N-1}(y)\,ds, \ (x,t)\in S_0. \label{greenf}
 \end{gather}
 In particular, we have
  \begin{gather}
 1=\int_{|y|\leq 3N/2}G(x,t; y, -3)\,dy 
- \int_{-3}^0\int_{|y|^2=3N/2} \frac{\partial G(x,t; y,s)}{\partial n}\,dH^{N-1}(y)\,ds,\nonumber\\
\frac{\partial G}{\partial n}\leq 0, G\geq 0, \min_{\substack{(x,t)\in S_{-1}\\ |y|^2\leq 3N/2}} G(x,t; y,-3) = \gamma \in (0,1). \nonumber 
 \end{gather}
 Hence, from \eqref{greenf} it follows that
 \begin{equation}\label{supnormS-1}
 \sup_{S_{-1}}w \leq \eta M + C\dashint_{|y|^2\leq 3N/2}w(y,-3)\,dy,
 \end{equation}
 where $\eta=1-\gamma \in (0,1)$, and $C>0$ is a constant depending on $N$. By using parabolic scaling $(x,t)\to (2^{2(l+1)}x,2^{4(l+1)}t)$ we derive from \eqref{supnormS-1} that
 \begin{equation}\label{supnormSl}
 \sup_{S_{l}}w \leq \eta M_{l+1} + C\dashint_{|y|^2\leq \frac{3N}{2}2^{4(l+1)}}w(y,-3\cdot 2^{4(l+1)})\,dy,
 \end{equation}
 whenever $\mathcal{H}w=0$ in $S_{l+1}$, $0\leq w \leq M_{l+1}$, $w$ continuous on $\partial_pS_{l+1}$. 
 
 \textit{Sub-step 5.3: Sequence $\{M_l\}$ has a geometric subsequence.} Now we can prove that due to \eqref{supnormSl} and \eqref{liminf0}, the sequence $\{M_l\}$ must possess a geometric subsequence with factor larger than 1. First, we fix some number $\eta<\mu<1$, and consider a cylinder $S_l$. If $M_l\neq 0$, then for some $\gamma>0$ we have $\gamma\leq M_{l'}\leq 1$ for all $l'\geq l$. Therefore, by \eqref{liminf0} we can always choose $l'>l$ such that 
 \begin{equation}\label{supnormSl'}
  \eta M_{l'+1} + C\dashint_{|y|^2\leq \frac{3N}{2}2^{4(l'+1)}}\theta(y,-3\cdot 2^{4(l'+1)})\,dy \leq \mu M_{l'+1}.
 \end{equation}
 Note that $0\leq \theta\leq 1$ on $\mathbb{R}^{N+1}_-$, and $\theta$ is smooth on the base of $S_{l'+1}$. Therefore, there is a continuous function $f$ defined on $\partial_pS_{l'+1}$ such that 
 \begin{equation}\label{supnormSl'+1}
 \theta\leq f \leq M_{l'+1} \ \text{on} \ \partial_pS_{l'+1},
 \end{equation}
  with
 \[  f=\theta \ \text{on} \ \{(y,\tau): |y|^2\leq \frac{3N}{2}2^{4(l'+1)}, \tau=-3\cdot 2^{4(l'+1)}\}. \]
Consider a Perron's solution $w=H^{S_{l'+1}}_f$. Then, $w$ is a parabolic function in $S_{l'+1}$, and from the first inequality of \eqref{supnormSl'+1} by the maximum principle it follows that $\theta\leq w$ in $S_{l'+1}$. On the other hand, from the second inequality of \eqref{supnormSl'+1} it follows that $w\leq M_{l'+1}$. Hence, by taking into account \eqref{supnormSl'}
we have
\begin{gather*}
M_l\equiv \sup_{S_l}\theta \leq \sup_{S_{l'}}\theta \equiv M_{l'} \leq \sup_{S_{l'}} w \nonumber\\
\leq  \eta \sup_{S_{l'+1}}w + C\dashint_{|y|^2\leq \frac{3N}{2}2^{4(l'+1)}}\theta(y,-3\cdot 2^{4(l'+1)})\,dy\nonumber\\
 \leq \eta M_{l'+1} + C\dashint_{|y|^2\leq \frac{3N}{2}2^{4(l'+1)}}\theta(y,-3\cdot 2^{4(l'+1)})\,dy\leq \mu M_{l'+1}. \nonumber\\
\end{gather*}
 Therefore, given any $l$, there exists a $l'\geq l$ such that
 \[ M_{l'+1}\geq \frac{1}{\mu} M_l. \]
 Since $\frac{1}{\mu} >1$ and $M_l$ is an increasing, this proves that either
 \[ M_l \equiv 0 \ \text{for all} \ l, \]
 or
 \[ \lim_{l\to +\infty} M_l=+\infty. \]
 The latter is excluded, since $M_l \leq 1$. Hence, the parabolic measure of $\mathcal{O}$ vanishes, and the regularity of $\mathcal{O}$ is proved.  \hfill{$\square$}
 
 Proof of Theorem~\ref{wienertest2} is a consequence of Lemma IV.1 of \cite{AbdullaI}, Lemma~\ref{hpotentialmeasure}, and the mapping \eqref{Appellmap}. Indeed, claim (5) of Lemma IV.1 of \cite{AbdullaI} implies that Problems $\mathcal{A}_\gamma|_\Omega$ and $\tilde{\mathcal{A}}_\gamma|_{A\Omega}$ are equivalent, i.e. $u_*\equiv 0$ if and only if $\tilde{u}_*\equiv 0$. On the other hand, claim (iii) of Lemma~\ref{hpotentialmeasure} implies that the divergence (or convergence) of the series \eqref{wtest} for $\Omega\subset\mathbb{R}^{N+1}_+$ is equivelant to the divergence (or convergence) of the series \eqref{wtest2} for $A\Omega\subset\mathbb{R}^{N+1}_-$. 
 
 \section{Proofs of Theorems~\ref{main theorem} and \ref{main theorem2}}\label{proofs}
Equivalence $(1) \Leftrightarrow (5)$ in Theorems~\ref{main theorem} and \ref{main theorem2} is proved in Theorems~\ref{wienertest} and \ref{wienertest2}. Equivalence $(1) \Leftrightarrow (2) \Leftrightarrow (3)$ in Theorems~\ref{main theorem} and \ref{main theorem2} is proved in \cite{AbdullaI} (Theorems III.1 and III.3).\\
To prove the equivalence $(1) \Leftrightarrow (4)$ in Theorem~\ref{main theorem2} with $\gamma=0$, first note that 
\begin{equation}\label{1infinity1}
H_{1_{\infty}}^{\tilde{\Omega}} = 1- H_{1-1_{\infty}}^{\tilde{\Omega}}
\end{equation}
and, therefore
\begin{equation}\label{1infinity2}
H_{1_{\infty}}^{\tilde{\Omega}} \equiv 0 \ \Leftrightarrow \ H_{1-1_{\infty}}^{\tilde{\Omega}} \equiv 1.
\end{equation}
We also have a relation
\begin{equation}\label{1infinity3}
R_1^{\tilde{\Omega}^c}(z) \equiv \bar{H}_{1-1_{\infty}}^{\tilde{\Omega}}(z) \equiv H_{1-1_{\infty}}^{\tilde{\Omega}}, \ \text{in} \ \tilde{\Omega}.
\end{equation}
The first relation of \eqref{1infinity3} is a direct consequence of the definition of the reduction and upper Perron's solution. The second relation of \eqref{1infinity3} follows from the resolutivity of the boundary function $1-1_{\infty}$ \cite{Doob}. According to the Lemma 4.4 in \cite{Abdulla2} (see also Lemma IV.2 of \cite{AbdullaI}) we have
\begin{equation}\label{1infinity4}
H_{1_{\infty}}^{\tilde{\Omega}}\equiv 0 \ \Leftrightarrow \ H_{1_{\infty}}^{\tilde{\Omega}\cap \{t\leq \delta\}}\equiv 0, \ \text{for any} \ \delta<0.
\end{equation}
Therefore, the relations \eqref{1infinity1},\eqref{1infinity2} are valid by replacing $\tilde{\Omega}$ with $\tilde{\Omega}\cap \{t\leq \delta\}, \delta<0$. 
This proves the equivalence $(1) \Leftrightarrow (4)$ in Theorem~\ref{main theorem2} with $\gamma=0$. 
Applying Lemma IV.1 of \cite{AbdullaI} and Lemma~\ref{reductionof1},
the equivalence $(1) \Leftrightarrow (4)$ in Theorem~\ref{main theorem} with $\gamma=0$ follows. Applying the translation $x\mapsto x+\gamma$, the equivalence of $(1) \Leftrightarrow (4)$ in Theorem~\ref{main theorem} with $\gamma\neq 0$ easily follows. Applying Lemma IV.1 of \cite{AbdullaI} and Lemma~\ref{reductionof1} again,
the equivalence $(1) \Leftrightarrow (4)$ in Theorem~\ref{main theorem2} with $\gamma\neq 0$ follows. \\

\section{Proofs of Theorems~\ref{geometricmaintheorem} and ~\ref{tildegeometricmaintheorem}}
Proofs of Theorems~\ref{geometricmaintheorem} and ~\ref{tildegeometricmaintheorem} is based on the following scheme:
\begin{equation}\label{geomainthm}
Thm.~\ref{tildegeometricmaintheorem}|_{\gamma=0} \stackrel{A^{-1}}{\Longrightarrow} Thm.~\ref{geometricmaintheorem}|_{\gamma=0} \stackrel{x\mapsto x+\gamma}{\Longrightarrow} Thm.~\ref{geometricmaintheorem}|_{\gamma\neq 0}\stackrel{A}{ \Longrightarrow} Thm.~\ref{tildegeometricmaintheorem}|_{\gamma\neq 0}   
\end{equation}
Precisely, the proof of Theorem~\ref{geometricmaintheorem} is pursued with the following three-step procedure:
\begin{equation}\label{thm42proof}
\Pi_\gamma:
\left\{
\begin{array}{l}
\quad\text{$h$-parabolic measure of $\mathcal{O}$ vanishes for the region $\mathcal{G}^1_\rho$}; \\
\quad\text{$\mathcal{O}$ is simultaneously $h$-regular or $h$-irregular for $\mathcal{G}^1_\rho$ and $\mathcal{G}^2_\rho$}; \\
\quad\text{$\mathcal{O}$ is $h$-regular for any open set $\Omega\subseteq\mathcal{G}^2_\rho$ for some $\delta>0$};
\end{array}\right.
\end{equation}
Equivalently, the proof of Theorem~\ref{tildegeometricmaintheorem} is pursued with a similar procedure:
\begin{equation}\label{thm48proof}
\tilde{\Pi}_\gamma:
\left\{
\begin{array}{l}
\quad\text{$\tilde{h}$-parabolic measure of $\infty$ vanishes for the region $\mathcal{L}^1_\rho$}; \\
\quad\text{$\infty$ is simultaneously $\tilde{h}$-regular or $\tilde{h}$-irregular for $\mathcal{L}^1_\rho$ and $\mathcal{L}^2_\rho$}; \\
\quad\text{$\infty$ is $\tilde{h}$-regular for any open set $\tilde{\Omega}\subseteq\mathcal{L}^2_\rho$ for some $\delta<0$};
\end{array}\right.
\end{equation}
If $\gamma=0$ ($\tilde{h}\equiv 1$), following the procedure $\tilde{\Pi}_0$, Theorem~\ref{tildegeometricmaintheorem} is proved in \cite{Abdulla2} (see Theorem 1.3, Lemma 2.5, and Lemma 2.4). To pursue the proof scheme \eqref{geomainthm}, we introduce a linear operator $\mathcal{T}:\mathcal{R} \mapsto \mathcal{R}^{'}$ as follows:
\begin{equation}\label{rhomapping}
\rho\in\mathcal{R}\ \mapsto \mathcal{T}\rho\in \mathcal{R}^{'}: \ \mathcal{T}\rho(y)= \rho\Big (-\frac{1}{y} \Big )
\end{equation}
The operator $\mathcal{T}$ is a bijection with inverse operator $\mathcal{T}^{-1}:\mathcal{R}^{'} \mapsto \mathcal{R}$ being in the same form:
\begin{equation}\label{inverserhomapping}
\rho\in\mathcal{R}^{'} \mapsto \mathcal{T}^{-1}\rho\in \mathcal{R}: \ \mathcal{T}^{-1}\rho(y)= \rho\Big (-\frac{1}{y} \Big )
\end{equation}
In particular, integral tests \eqref{tilderhoint} and \eqref{rhoint} are exchanged under the bijection $\mathcal{T}$, and so are the properties \eqref{rhodiff} at $t\downarrow -\infty$ and $t\downarrow 0$ respectively. To prove the first claim of \eqref{geomainthm}, we fix any $\rho\in\mathcal{R}$, and verify that
\begin{equation}\label{G1inL1}
\mathcal{G}^1_\rho \subset A^{-1}\mathcal{L}^1_{\mathcal{T}\rho}, \ \mathcal{G}^2_\rho \subset A^{-1}\mathcal{L}^2_{\mathcal{T}\rho},
\end{equation}
where we replace $\delta \to -(4\delta)^{-1}$ in $\mathcal{G}^1_\rho$ and $\mathcal{G}^2_\rho$. Since $\mathcal{T}\rho\in\mathcal{R}^{'}$, from Theorem~\ref{tildegeometricmaintheorem}$|_{\gamma=0}$ it follows that $\infty$ is regular for $\mathcal{L}^1_{\mathcal{T}\rho}$. By Lemma IV.1 of \cite{AbdullaI}, $\mathcal{O}$ is $h$-regular for $A^{-1}\mathcal{L}^1_{\mathcal{T}\rho}$. Therefore, from \eqref{G1inL1} and Lemma IV.2 of \cite{AbdullaI} it follows that $\mathcal{O}$ is $h$-regular for $\mathcal{G}^1_\rho$, which is the first claim of $\Pi_0$ in \eqref{thm42proof}. The second claim of $\Pi_0$ is proved with the same argument. 
By Lemma 2.5 of \cite{Abdulla2}, the $\infty$ is regular for both $\mathcal{L}^1_{\mathcal{T}\rho}$ and $\mathcal{L}^1_{\mathcal{T}\rho}$. Therefore, from \eqref{G1inL1} and Lemma IV.2 of \cite{AbdullaI} it follows that $\mathcal{O}$ is $h$-regular for both $\mathcal{G}^1_\rho$ and $\mathcal{G}^2_\rho$. Finally, the third claim of $\Pi_0$ follows from Lemma IV.2 of \cite{AbdullaI}. Noting that under the conditions of Theorem~\ref{tildegeometricmaintheorem}, $\Omega\subseteq\mathcal{G}^2_\rho$, the third claim of $\Pi_0$ follows from Lemma IV.2 of \cite{AbdullaI}. This completes the proof of the first claim of \eqref{geomainthm}. 

Second claim of \eqref{geomainthm} is straightforward by using a translation map $x \to x+\gamma$. Indeed, all three steps in $\Pi_\gamma$ are invariant with respect to the spatial translation map, and therefore $\Pi_0 \Longrightarrow \Pi_\gamma$. 

To prove the last claim of \eqref{geomainthm}, we fix any $\rho\in\mathcal{R}^{'}$, and verify that
\begin{equation}\label{L1inG1}
\mathcal{L}^1_\rho \equiv A \Omega^1, \ \mathcal{L}^2_\rho \equiv A \Omega^2.
\end{equation}
where 
\begin{gather}
\Omega^1=\Big\{z: |x_1-\gamma_1|^2<4\xi\log\rho\Big (\frac{4t^2}{\xi}\Big ), \bar{x}\in \mathbb{R}^{N-1}, 0<t<-\frac{1}{4\delta}\Big\}\label{expG1rho}\\
\Omega^2=\Big\{z: -2\Big (\xi\log\rho\Big (\frac{4t^2}{\xi}\Big )\Big )^{\frac{1}{2}}<x_1-\gamma_1<+\infty, \bar{x}\in \mathbb{R}^{N-1}, 0<t<-\frac{1}{4\delta}\Big\}\label{expG2rho}
\end{gather}
It is easy to observe that
\begin{equation}\label{G1interomega}
\Omega^1\cap \mathcal{G}^1_{\mathcal{T}^{-1}\rho} \equiv  \mathcal{G}^1_{\mathcal{T}^{-1}\rho}, \ \Omega^2\cap \mathcal{G}^2_{\mathcal{T}^{-1}\rho} \equiv \mathcal{G}^2_{\mathcal{T}^{-1}\rho}. 
\end{equation}
Since $\mathcal{T}^{-1}\rho\in\mathcal{R}$, from Theorem~\ref{geometricmaintheorem} it follows that $\mathcal{O}$ is $h$-regular for $\mathcal{G}^1_{\mathcal{T}^{-1}\rho}$ and $\mathcal{G}^2_{\mathcal{T}^{-1}\rho}$. To demonstrate that $\mathcal{O}$ is $h$-regular for $\Omega^1$ and $\Omega^2$, first note that
\begin{equation}\label{OhregOmega1}
 \ ^hH_{1_{\{\mathcal{O}\}}}^{\Omega^i} |_{\mathcal{G}^i_{\mathcal{T}^{-1}\rho}}\equiv \ ^hH_{g_i}^{\mathcal{G}^i_{\mathcal{T}^{-1}\rho}}, \ i=1,2
\end{equation}
where
\begin{equation}\label{OhregOmega1g}
g_i=
\left\{
\begin{array}{l}
1_{\{O\}} \ \quad\text{on} \ \partial\mathcal{G}^i_{\mathcal{T}^{-1}\rho}\cap\partial \Omega^i \\
^hH_{1_{\{\mathcal{O}\}}}^{\Omega^i}, \ \quad\text{on} \ \partial\mathcal{G}^i_{\mathcal{T}^{-1}\rho}\cap\Omega^i, \ i=1,2.
\end{array}\right.
\end{equation}
From \eqref{G1interomega} it follows that
\begin{equation}\label{glim}
\lim\limits_{z\to\mathcal{O}, z\in\partial\mathcal{G}^i_{\mathcal{T}^{-1}\rho}}g_i=0, \  i=1,2.
\end{equation}
Since $\mathcal{O}$ is $h$-regular for $\mathcal{G}^i_{\mathcal{T}^{-1}\rho}$, using boundary regularity counterpart (\cite{AbdullaI}, Section II C (40)) we have
\begin{equation}\label{G1rholim}
\lim\limits_{z\to\mathcal{O}, z\in\mathcal{G}^i_{\mathcal{T}^{-1}\rho}} \ ^hH_{g_i}^{\mathcal{G}^i_{\mathcal{T}^{-1}\rho}}(z) =0, \  i=1,2.
\end{equation}
Due to \eqref{G1interomega} and \eqref{OhregOmega1} it follows that
\begin{equation}\label{OhregOmega12}
\lim\limits_{z\to\mathcal{O}, z\in\mathcal{G}^i_{\mathcal{T}^{-1}\rho}} \ ^hH_{1_{\{\mathcal{O}\}}}^{\Omega^i}(z) =0, \  i=1,2,
\end{equation}
and therefore by the $h$-maximum principle we have
\begin{equation}\label{hmaxpr}
 \ ^hH_{1_{\{\mathcal{O}\}}}^{\Omega^i}\equiv 0, \ i=1,2.
\end{equation}
Hence, $\mathcal{O}$ is $h$-regular for $\Omega^i, i=1,2$. 
By Lemma IV.1 of \cite{AbdullaI} $\infty$ is $\tilde{h}$-regular for $A\Omega^i, i=1,2$. Therefore, from \eqref{L1inG1} it follows that $\infty$ is $\tilde{h}$-regular for $\mathcal{L}^i_\rho, i=1,2$, which confirms first two claims of $\tilde{\Pi}_\gamma$ in \eqref{thm48proof}. Noting that under the conditions of Theorem~\ref{tildegeometricmaintheorem}, $\Omega\subseteq\mathcal{L}^2_\rho$, the third claim of $\tilde{\Pi}_\gamma$ follows from Lemma IV.2 of \cite{AbdullaI}. This completes the proof of the last claim of \eqref{geomainthm}.  \hfill{$\square$}

\bibliographystyle{amsplain}

\end{document}